\documentstyle[amscd,amssymb,verbatim,diagrams,12pt]{amsart}
\renewcommand{\int}{{\operatorname{int}}}

\newcommand{\Tot}{{\operatorname{Tot}}}

\newcommand{\bC}{{\bf C}}
\newcommand{\bc}{{\bf c}}
\newcommand{\bn}{{\bf n}}

\newcommand{\SA}{\operatorname{SA}}

\newcommand{\fG}{{\frak G}}

\renewcommand{\mod}{\operatorname{mod}}

\newcommand{\gr}{\operatorname{gr}}
\newcommand{\cusp}{\operatorname{cusp}}

\newcommand{\OO}{{\cal O}}
\newcommand{\RR}{{\cal R}}

\newcommand{\Det}{\operatorname{Det}}

\newcommand{\Polar}{\operatorname{Polar}}

\newcommand{\coker}{\operatorname{coker}}
\newcommand{\DD}{{\cal D}}
\newcommand{\NN}{{\cal N}}
\newcommand{\KK}{{\cal K}}

\newcommand{\be}{{\bf e}}

\newcommand{\BB}{{\cal B}}

\newcommand{\unit}{{\bf 1}}
\newcommand{\G}{{\Bbb G}}

\newcommand{\mg}{{\frak m}}

\newcommand{\lan}{\langle}
\newcommand{\ran}{\rangle}

\newcommand{\CC}{{\cal C}}
\newcommand{\UU}{{\cal U}}
\newcommand{\WW}{{\cal W}}

\newcommand{\supp}{\operatorname{supp}}
\newcommand{\Spec}{\operatorname{Spec}}

\newcommand{\Proj}{\operatorname{Proj}}
\renewcommand{\P}{{\Bbb P}}
\newcommand{\sV}{{\Bbb V}}

\newcommand{\si}{\sigma}

\newcommand{\Pic}{\operatorname{Pic}}

\newcommand{\de}{\delta}

\renewcommand{\ker}{\operatorname{ker}}
\newcommand{\im}{\operatorname{im}}

\newcommand{\A}{{\Bbb A}}

\numberwithin{equation}{subsection}
\newcommand{\bm}{{\bf m}}
\newcommand{\bp}{{\bf p}}
\newcommand{\bD}{{\bf D}}
\newcommand{\br}{{\bf r}}

\newtheorem{thm}{Theorem}[subsection]
\newtheorem{prop}[thm]{Proposition}
\newtheorem{lem}[thm]{Lemma}
\newtheorem{cor}[thm]{Corollary}
{  \theoremstyle{definition}
\newtheorem{defi}[thm]{Definition}
\newtheorem{ex}[thm]{Example}
\newtheorem{exs}[thm]{Examples}
\newtheorem{rem}[thm]{Remark}
\newtheorem{rems}[thm]{Remarks}
}

\newcommand{\Pf}{\noindent {\it Proof}}
\newcommand{\id}{\operatorname{id}}

\newcommand{\ov}{\overline}
\newcommand{\we}{\wedge}

\newcommand{\rk}{\operatorname{rk}}

\newcommand{\MM}{{\cal M}}

\newcommand{\XX}{{\cal X}}
\newcommand{\HH}{{\cal H}}

\newcommand{\VV}{{\cal V}}
\newcommand{\SS}{{\cal S}}
\newcommand{\LL}{{\cal L}}

\newcommand{\forg}{\operatorname{for}}
\newcommand{\Om}{\Omega}

\newcommand{\Hom}{\operatorname{Hom}}

\renewcommand{\a}{\alpha}

\newcommand{\om}{\omega}

\newcommand{\la}{\lambda}

\newcommand{\N}{{\Bbb N}}
\newcommand{\R}{{\Bbb R}}
\newcommand{\Z}{{\Bbb Z}}
\newcommand{\Q}{{\Bbb Q}}

\newcommand{\wt}{\widetilde}
\newcommand{\ot}{\otimes}

\newcommand{\sub}{\subset}
\newcommand{\ed}{\qed\vspace{3mm}}

\newcommand{\Kr}{\operatorname{Kr}}
\newcommand{\WMA}{\operatorname{WMA}}
\newcommand{\MA}{\operatorname{MA}}

\newcommand{\bw}{{\bf w}}
\newcommand{\ba}{{\bf a}}
\newcommand{\bb}{{\bf b}}
\newcommand{\bS}{{\bf S}}

\newcommand{\Si}{\Sigma}

\newcommand{\sslash}{\mathbin{/\mkern-6mu/}}

\title{Moduli of curves, Gr\"obner bases, and the Krichever map}
\author{Alexander Polishchuk}
\thanks{Supported in part by the NSF grant DMS-1400390}

\begin{document}

\maketitle
\begin{abstract}
We study moduli spaces of (possibly non-nodal)
curves $(C,p_1,\ldots,p_n)$ of arithmetic genus $g$ with $n$ smooth marked points,
equipped with nonzero tangent vectors, such that $\OO_C(p_1+\ldots+p_n)$ is ample and
$H^1(\OO_C(a_1p_1+\ldots+a_np_n))=0$ for given weights $\ba=(a_1,\ldots,a_n)$ such that $a_i\ge 0$ and $\sum a_i=g$.
We show that each such moduli space $\wt{\UU}^{ns}_{g,n}(\ba)$
is an affine scheme of finite type, and the Krichever map identifies it with the quotient of
an explicit locally closed subscheme of the Sato Grassmannian by the free action of the group of changes of formal
parameters. We study the GIT quotients of $\wt{\UU}^{ns}_{g,n}(\ba)$ by the natural torus action and show
that some of the corresponding stack quotients give modular compactifications of $\MM_{g,n}$ with projective coarse
moduli spaces. More generally, using similar techniques,
we construct moduli spaces of curves with chains of divisors supported at marked points, with prescribed number of sections,
which in the case $n=1$ corresponds to specifying the Weierstrass gap sequence at the marked point. 
\end{abstract}

\section*{Introduction}

Recently, modular compactifications of $\MM_{g,n}$ parametrizing curves with not necessarily nodal singularities 
attracted a lot of attention, in particular,
in connection with Hassett-Keel program to study log-canonical models of $\ov{M}_{g,n}$ (surveyed
in \cite{FS}). In the present paper, generalizing \cite{P-ainf}, we construct new examples of such modular compactifications 
that arise in connection with the Krichever map.
The idea of the Krichever map is to study a projective curve $C$ with a smooth point $p$ via the subspace of all Laurent expansions of regular functions on $C\setminus \{p\}$ with respect to a formal parameter at $p$, viewed as a point
in the Sato Grassmannian. In the case when $C$ is reducible, it is better to use an analog of this construction with several marked points, such that there is at least one marked point on each component. The corresponding morphism from the moduli space of curves with fixed formal parameters at marked points to the Sato Grassmannians has been used to study the geometry of the moduli spaces of curves (see \cite{ADKP}, \cite{BS}, \cite{Kontsevich}, \cite{Schwarz-L}, \cite{SW}).
In the present paper we study situations when there is a canonical choice of a formal parameter (with a given
1st jet) at the marked points. This leads to an identification of certain moduli stacks of curves with quotients
of explicit affine varieties by the torus actions.                         

The main objects of our study are the moduli stacks of curves with nonspecial divisors, described as follows.
For every $g\ge 1$ let $X(g,n)$ denote the subset $\Z_{\ge 0}^n$ consisting of
$\ba=(a_1,\ldots,a_n)$ such that $a_1+\ldots+a_n=g$. For each $\ba\in X(g,n)$ we consider
(not necessarily smooth) reduced connected projective curves $C$ of arithmetic genus $g$ with smooth marked points
$p_1,\ldots,p_n$, such that $h^1(\OO_C(a_1p_1+\ldots+a_np_n))=0$ and the line bundle $\OO_C(p_1+\ldots+p_n)$
is ample. We denote the moduli stack of such data by $\UU^{ns}_{g,n}(\ba)$.
We rigidify these data by adding a choice of nonzero tangent vectors at every marked point and denote by
$\wt{\UU}^{ns}_{g,n}(\ba)\to \UU^{ns}_{g,n}(\ba)$ the corresponding $\G_m^n$-torsor.

In the case when $n=g$ and all the weights $a_i$ are equal to $1$, we recover the moduli stacks $\UU^{ns,a}_{g,g}$ and
$\wt{\UU}^{ns,a}_{g,g}$ considered in \cite{P-ainf}. 
We proved in \cite{P-ainf} that if we restrict the base to 
$\Spec(\Z[1/6])$ then $\wt{\UU}^{ns}_{g,g}(1,\ldots,1)$ is in fact an affine scheme of finite type over $\Z[1/6]$.
One of the main results of this paper is the following generalization of this to $\wt{\UU}^{ns}_{g,n}(\ba)$.

\medskip

\noindent
{\bf Theorem A.} {\it (i) There exists an integer $N$ depending only on $\ba$, such that working
over $\Z[1/N]$ we get that $\wt{\UU}^{ns}_{g,n}(\ba)$ is an affine scheme of finite type over $\Z[1/N]$.
It is equipped with a $\G_m^n$-action corresponding to rescaling of the tangent vectors at the marked points.

\noindent
(ii) For each $\ba\in X(g,n)$ there is a natural ``forgetting last point" map 
$$\forg_{n+1}:\wt{\UU}^{ns}_{g,n+1}(\ba,0)\to \wt{\UU}^{ns}_{g,n}(\ba)$$ 
and a compatible map of universal affine curves
$$C_{g,n+1}(\ba,0)\setminus\{p_1,\ldots,p_{n+1}\}\to C_{g,n}(\ba)\setminus\{p_1,\ldots,p_n\}.$$

\noindent
(iii) For a collection $\ba_1,\ldots,\ba_r\in X(g,n)$ let us define $\wt{\UU}^{ns}_{g,n}(\ba_1,\ldots\ba_r)$
as the intersection $\wt{\UU}^{ns}_{g,n}(\ba_1)\cap\ldots\cap\wt{\UU}^{ns}_{g,n}(\ba_r)$ inside the moduli stack of
curves. Then $\wt{\UU}^{ns}_{g,n}(\ba_1,\ldots\ba_r)$ is a principal affine open subscheme in
$\wt{\UU}^{ns}_{g,n}(\ba_1)$.}

The key feature in \cite{P-ainf} that led to an explicit embedding of $\wt{\UU}^{ns}_{g,g}$ into an affine space is the existence
of a certain canonical basis of $H^0(C\setminus\{p_1,\ldots,p_g\})$ for $(C,p_1,\ldots,p_g)$ in this moduli space.
In characteristic zero this is complemented by the construction of the canonical formal parameters at each marked point.
In the present paper we explain both these phenomena in terms of the Krichever map to the Sato Grassmannian
and generalize them to $\wt{\UU}^{ns}_{g,n}(\ba)$.
Namely, we associate with each $\ba\in X(g,n)$ a cell $SG^\ba$ in the Sato Grassmannian of subspaces in
$\HH=\bigoplus_{i=1}^n k((t_i))$, such that the Krichever map associated with $n$ marked points and formal parameters at them
lands in $SG^\ba$ if and only if $h^1(\OO_C(a_1p_1+\ldots+a_np_n))=0$
(for the precise definition of $SG^\ba$ see Section \ref{cell-sec}).
The existence of a canonical basis (constructed for the universal curve) in $H^0(C\setminus\{p_1,\ldots,p_g\})$ has to do with a section for the action of the pro-unipotent group $\fG$ of formal changes of variables $t_1,\ldots,t_n$, 
trivial modulo $(t_i^2)$, on $SG^\ba$ (see Theorem B below).

Let $ASG$ be the closed subscheme in the Sato Grassmannian corresponding to subspaces $W\sub\HH$
such that $1\in W$ and $W\cdot W\sub W$, and let $ASG^\ba=ASG\cap SG^\ba$. 
Note that this subscheme is preserved by the $\fG$-action.
Our next result is that this action is free and the quotient is isomorphic to the moduli space of curves considered above. 

\medskip

\noindent
{\bf Theorem B.} {\it Let us work over $\Q$. For any $\ba\in X(g,n)$
the natural action of the group $\fG$ of formal changes of variables on $SG^\ba$ admits
a section $\Si^{\ba,i_0}\sub SG^{\ba}$ (depending on a choice of $i_0$ such that $a_{i_0}>0$, in the case when not
all $a_i$ are positive), isomorphic to an infinite-dimensional affine space. We have a
commutative diagram of affine schemes
\begin{diagram}
\wt{\UU}^{ns,(\infty)}_{g,n}(\ba)&\rTo^{\Kr}_{\sim}& ASG^{\ba}\\
\dTo{}&&\dTo{}\\
\wt{\UU}^{ns}_{g,n}(\ba)&\rTo^{\ov{\Kr}}_{\sim}& ASG^{\ba}/\fG
\end{diagram}
where the horizontal arrows are given by the Krichever map and the vertical arrows are quotients by the free action of
$\fG$. 
}

Note that the conditions $1\in W$, $W\cdot W\sub W$ were known to characterize the image of the Krichever map in related
contexts (see \cite[Sec.\ 6]{SW}, \cite{Mulase}). However, as far as we know, this map was usually considered at the set-theoretic
level, not at the level of moduli functors (an exception to this is the construction of \cite[Prop.\ 6.3]{MP-PM} 
for the moduli of integral curves). 
In fact, it seems that even to define the Krichever map as a map from some
moduli stack of not necessarily irreducible curves $C$ with fixed formal parameters at the marked points $p_1,\ldots,p_n$, one needs to impose some restrictions on $(C,p_1,\ldots,p_n)$. For example, we define such a map for the moduli of $(C,p_1,\ldots,p_n)$ such that
a sufficiently high multiple of the divisor $p_1+\ldots+p_n$ has vanishing $h^1$ 
%(note that this is automatic for irreducible curves),
(see Proposition \ref{Kr-map-prop}).

%Proof of isomorphism of the moduli space with the Krichever data: construct isomorphism of the enhanced moduli,
%obtained by fixing formal parameters, the actions are compatible, use the section.

The proof of Theorems A and B starts with an explicit construction of the section $\Sigma^{\ba,i_0}\sub SG^\ba$ of the $\fG$-action.
Then, generalizing \cite[Lem.\ 1.2.2]{P-ainf},
we introduce an affine scheme $S_{GB}$ of finite type which parametrizes
commutative algebras with Gr\"obner bases of special type. The proof is achieved by studying natural morphisms
$$\wt{\UU}^{ns}_{g,n}(\ba)\rTo^{\ov{\Kr}} \Sigma^{\ba,i_0}\cap ASG^\ba\to S_{GB}\to \wt{\UU}^{ns}_{g,n}(\ba)$$
(see Section \ref{proof-A-B-sec}).

We also study the GIT quotients of the schemes $\wt{\UU}^{ns}_{g,n}(\ba)$ by the natural $\G_m^n$-action.
%rescaling the tangent vectors at the marked points. 
In our next result, generalizing \cite[Prop.\ 2.4.2]{P-ainf},
we show that these quotients, which depend on a character $\chi$ of $\G_m^n$, are projective, and 
identify an explicit region of $\chi$ for which the generic curve (resp., every smooth curve) is stable. 

\medskip

\noindent
{\bf Theorem C (see Theorem \ref{GIT-thm}(i),(ii) and Corollary \ref{GIT-cor}).} 
{\it Let us work over an algebraically closed field $k$ 
of characteristic zero. For any $\chi\in\Z^n$ the GIT quotient
$\wt{\UU}^{ns}_{g,n}(\ba)\sslash_\chi \G_m^n$ is a projective scheme over $k$. Let $\bC_\ba\sub\R^n$ denote the closed cone
generated by the vectors $(a_i+1)\be_i-a_j\be_j$, with $i\neq j$, and let $\bC_0$ be the subcone generated by
the standard basis vectors $\be_i$. 

\noindent
(i) For $\chi\not\in\bC_\ba$ one has $\wt{\UU}^{ns}_{g,n}(\ba)\sslash_\chi \G_m^n=\emptyset$.
For $\chi$ in the interior of  $\bC_\ba$ and for any smooth curve $C$ of genus $g$, a point 
$(C,p_1,\ldots,p_n,v_1,\ldots,v_n)\in \wt{\UU}^{ns}_{g,n}(\ba)$ corresponding to generic points 
$p_1,\ldots,p_n$, is $\chi$-stable.

\noindent
(ii) For $\chi$ in the interior of $\bC_0$ every smooth curve is $\chi$-stable.
For such $\chi$ every semistable point is stable and the notion of semistability does not depend on $\chi$.
The corresponding quotient stack $\wt{\UU}^{ns}_{g,n}(\ba)^{ss}_\chi/\G_m^n$ is proper.
}

Note that Theorem C(ii) implies that for $\chi$ in the interior of $\bC_0$ the irreducible component of the
quotient stack $\wt{\UU}^{ns}_{g,n}(\ba)^{ss}_\chi/\G_m^n$, containing smooth curves, is a
modular compactifications of $\MM_{g,n}$ in the sense of \cite{Smyth}. The same is true for other generic characters
in $\bC_\ba$, i.e., for $\chi$ in the complement to a finite number of codimension $1$ walls, however, 
for arbitrary $\ba$ we do not know explicitly the set of walls (cf.\ Remark \ref{GIT-rems}.1).

From Theorem B we see that $ASG\cap\Sigma^{\ba,i_0}$ is a section for the free action of $\fG$ on $ASG^\ba$, in particular, 
$$ASG^{\ba}/\fG \simeq ASG\cap\Sigma^{\ba,i_0}.$$
Combining this isomorphism with $\ov{\Kr}$ we get an explicit $\G_m^n$-equivariant
embedding of $\wt{\UU}^{ns}_{g,n}(\ba)$ into the infinite-dimensional affine space $\Sigma^{\ba,i_0}$. 
The corresponding coordinates $\a_{ij}[p,q]$ on $\wt{\UU}^{ns}_{g,n}(\ba)$ (see \eqref{coordinates-a-ij-eq}) descend
to sections $\ov{\a}_{ij}[p,q]$ of certain line bundles on $\UU^{ns}_{g,n}(\ba)$ (linear combinations of $\psi$-classes). 
In particular, we can define similar sections $\ov{\a}_{ij}[p,q]$ of line bundles on the open substack 
$\ov{\MM}_{g,n}(\ba)$ of stable curves $(C,p_1,\ldots,p_n)$ such  that $H^1(C,\OO(p_1+\ldots+p_n))=0$.
Generalizing \cite[Prop.\ 3.1.1]{FP}, we estimate the poles of these sections along the complement to 
$\ov{\MM}_{g,n}(\ba)$ in $\ov{\MM}_{g,n}$ (see Theorem \ref{pole-thm}). This can be considered as the first step
towards studying the birational maps from $\ov{\MM}_{g,n}$ to the GIT quotients of $\wt{\UU}^{ns}_{g,n}(\ba)$ by
$\G_m^n$.

We apply similar ideas to analyze more general moduli spaces of curves. Namely,
we consider pointed curves $(C,p_1,\ldots,p_n)$, where $p_1,\ldots,p_n$ are smooth and distinct, with 
the following additional structure: an effective divisor $D$, supported at $\{p_1,\ldots,p_n\}$,
such that $h^1(D)=0$, and a chain of effective divisors between $0$ and $D$, such that 
$h^0$ is stipulated along the whole chain. We leave the precise formulation for Section \ref{special-curves-sec},
and state here our result in the case $n=1$. In this case we study the moduli space 
$\wt{\UU}_{g,1}[\ell_1,\ldots,\ell_g]$ of
irreducible pointed curves $(C,p)$ of arithmetic genus $g$, 
with a choice of a nonzero tangent vector at $p$, with a given Weierstrass gap
sequence $1=\ell_1<\ell_2<\ldots<\ell_g$ at $p$. This means that 
the function $n\mapsto h^1(np)$ jumps precisely at these values, and is constant on the intervals $[\ell_i,\ell_{i+1}-1]$,
$i=0,\ldots,g-1$ (where we set $\ell_0=0$).
The precise moduli problem can be formulated by considering families of
irreducible pointed curves $(\pi:C\to S,p:S\to C)$ with the requirement that $R^1\pi_*(mp)$ are locally free
of given ranks (see Corollary \ref{gap-moduli-cor}). We also define similar loci $ASG[\ell_1,\ldots,\ell_g]$ in the Sato Grassmannian, corresponding to subalgebas $W\sub\HH=k((t))$ with given dimensions of intersections with 
$t^{\ge -m}k[[t]]$, determined by $\ell_1,\ldots,\ell_g$ (see Section \ref{special-divisors-setup-sec}).

\medskip

\noindent
{\bf Theorem D.} {\it Let us work over $\Q$. For each gap sequence $1=\ell_1<\ldots<\ell_g$
the moduli stack $\wt{\UU}_{g,1}[\ell_1,\ldots,\ell_g]$ is an affine scheme of finite type over $\Q$. The action of $\fG$
on $ASG[\ell_1,\ldots,\ell_g]$ is free and admits a section. The Krichever map induces an isomorphism
$$\wt{\UU}_{g,1}[\ell_1,\ldots,\ell_g]\simeq ASG[\ell_1,\ldots,\ell_g]/\fG.$$
}

Note that moduli spaces similar to $\UU_{g,1}[\ell_1,\ldots,\ell_g]$ were constructed by St\"ohr \cite{Stohr} in
the special case when $\ell_g=2g-1$ (the maximal possible value for $\ell_g$), $\ell_2=2$ and $\ell_{g-1}\ge g$.
His approach was to work with Gorenstein curves and to use the canonical embedding, whereas we consider not necessarily
Gorenstein curves and use the Krichever map. As a byproduct of our more general approach we
are able to deduce that any projective irreducible and reduced curve of arithmetic genus $g\ge 1$, that has a smooth point $p$ such that $h^1((2g-2)p)=1$, is Gorenstein and has $\om_C\simeq \OO((2g-2)p)$ (see Corollary \ref{Gorenstein-cor}).
Thus, it seems plausible that St\"ohr's moduli spaces are special cases of our moduli spaces 
(see Remark \ref{Stohr-rem}.1 for more discussion).
%(see also Remark \ref{Pinkham-rem}).

%A similar picture exists for curves with $n>1$ points (see Theorem \ref{special-divisors-thm}).

%The second motivation is the problem of describing explicitly the equations of $\wt{\UU}^{ns}_{g,g}$ in the above
%affine embedding. 

In \cite{P-ainf} we also proved that $\wt{\UU}^{ns}_{g,g}(1,\ldots,1)$ 
can be interpreted as the moduli space of minimal $A_\infty$-structures on
a certain finite-dimensional algebra. Elsewhere we will discuss a similar interpretation of
the moduli schemes $\wt{\UU}^{ns}_{g,n}(\ba)$.

%gluing maps.

The paper is organized as follows.

In Section \ref{moduli-krich-sec}, after introducing the basic framework for working with the Sato Grassmannian and the Krichever map, we study the moduli spaces $\wt{\UU}^{ns}_{g,n}(\ba)$. In addition to
proving Theorems A and B, we describe them as moduli spaces of {\it marked algebras of type $\ba$}, which
are commutative algebras equipped with some special filtrations---this notion captures the properties of the filtration
on $H^0(C\setminus\{p_1,\ldots,p_n\},\OO)$ by order of poles at the marked points (see Proposition \ref{marked-alg-prop}).
In comparison to the case $g=n$, $\ba=(1,\ldots,1)$, considered in \cite{P-ainf}, 
the use of Gr\"obner bases is slightly more technical for general $\ba$. In Section \ref{Groebner-sec} we identify a class of
Gr\"obner bases for commutative algebras arising from the Krichever map.

In Section \ref{more-moduli-sec} we further study the schemes $\wt{\UU}^{ns}_{g,n}(\ba)$ and the $\G_m^n$-action on them.
In \ref{special-curves-subsec} we construct some special curves corresponding to points in $\wt{\UU}^{ns}_{g,n}(\ba)$.
Then in \ref{subgroup-sec} 
we find a subgroup $\G_m$ in $\G_m^n$, with respect to which the coordinates on $\wt{\UU}^{ns}_{g,n}(\ba)$
have positive weight. In particular, we get that the action of this subgroup degenerates any curve in $\wt{\UU}^{ns}_{g,n}(\ba)$
to a certain cuspidal curve $C^{\cusp}(\ba)$. In \ref{generators-sec} we find a convenient set of generators in
the algebra of functions on $\wt{\UU}^{ns}_{g,n}(\ba)$, which is helpful in the analysis of the GIT stability conditions.
In \ref{GIT-sec} we study the GIT quotients of $\wt{\UU}^{ns}_{g,n}(\ba)$ by the $\G_m^n$-action, in particular, proving
Theorem C.  In \ref{extending-sec} we analyze the poles of the coordinates $\ov{\a}_{ij}[p,q]$, 
viewed as rational sections of line bundles on $\ov{\MM}_{g,n}$.

In Section \ref{examples-sec} we consider examples of the moduli schemes $\wt{\UU}^{ns}_{g,n}(\ba)$ with $g=1$. 
In particular, we explain the connection to the moduli spaces of genus $1$ curves studied in \cite{LP}.

Finally, in Section \ref{special-curves-sec} we study the moduli of curves with chains of divisors supported
at marked points and with prescribed $h^0$, proving Theorem D and its generalization to the case $n>1$,
Theorem \ref{special-divisors-thm}.

\medskip

\noindent
{\it Notation and conventions.}
All curves (over algebraically closed fields) are assumed to be reduced and connected.
We always assume that $g\ge 1$.
We denote by $\be_1,\ldots,\be_n$ the standard basis vectors in $\Z^n$.
By the infinite-dimensional affine space over $R$ we mean $\Spec(R[x_1,x_2,\ldots])$.

\medskip

\noindent
{\it Acknowledgment.}
 I am grateful to the Institut des Hautes \'Etudes Scientifiques 
where some of this work was done, for the hospitality and the excellent working conditions.
%This research is supported in part by the NSF grant DMS-1400390.

\section{Moduli spaces of curves via the Krichever map}\label{moduli-krich-sec}

\subsection{The Sato Grassmannian and the Krichever map}

Let $k$ be a field. Set
$$\HH=\HH_k=\bigoplus_{i=1}^n k((t_i)),$$
where $t_1,\ldots,t_n$ are variables, and for $N\in\Z$,
$$\HH_{\ge N}=\bigoplus_{i=1}^n k[[t_i]]t_i^N \sub \HH.$$
The Sato Grassmannian $SG=SG(\HH)$ parametrizes subspaces $W\sub\HH$,
such that the operator 
\begin{equation}\label{d-V-operator-eq}
d_W:W\oplus \HH_{\ge 0}\to \HH
\end{equation}
has finite-dimensional kernel and cokernel.
We denote these kernel and cokernel by
$$H^0(W):=W\cap\HH_{\ge 0}=1, \ \ H^1(W):= \HH/(W+\HH_{\ge 0}).$$

For each $v\in\HH$ we call the component of $v$ in $k((t_i))$ the {\it expansion of $v$ in $t_i$}.

The Sato Grassmannian can be defined as a scheme by gluing {\it open cells} that are identified with 
infinite-dimensional affine spaces (see \cite[Sec.\ 4.3]{BS},
\cite[Sec.\ 2]{SW}). Namely, for every collection $\bS$ of subsets $S_i\sub\Z$, $i=1,\ldots,n$, 
such that $S_i\setminus \N$ and $\N\setminus S_i$ are finite, there is an open subset $U_{\bS}\sub SG$ parametrizing subspaces $W$ such that 
$W\oplus \HH_{\bS}=\HH$, where
$$\HH_{\bS}=\bigoplus_{i=1}^n\prod_{j\in S_i}k t_i^j.$$
All such $W$ can be represented as graphs of 
linear maps $\phi:\HH_{\bS^c}\to \HH_\bS$, where 
$$\HH_{\bS^c}=\bigoplus_{i, j\not\in S_i}k t_i^j.$$
The components of $\phi$ give coordinates on $U_{\bS}$ identifying it with the infinite-dimensional affine space.
The transition maps are algebraic and defined over $\Z$, so in fact, $SG$ can be defined as a scheme over $\Z$.
More generally, for every lattice $L\sub\HH$, i.e., a subspace commensurable with $\bigoplus_{i=1}^n k[[t_i]]$,
there is an open subset $U_L\sub SG$ (defined over $k$), isomorphic to an infinite-dimensional affine space, 
consisting of $W$ such that $W\oplus L=\HH$.

For simplicity of exposition we will discuss below various subschemes of $SG$ in terms of $k$-points, where $k$ is a field.
However, there are natural analogs of these subschemes defined over $\Z$ using the above description
of $SG$ as a scheme glued from the open subsets $U_{\bS}$. One can also explicitly work with
$R$-points of $SG$, where $R$ is a commutative ring (see e.g. \cite[Sec.\ 2.A]{MP-PM}).

For an integer $g$ we denote by $SG(g)$ 
the connected component of $SG$ consisting of $W$ such that
$\dim H^0(W)-\dim H^1(W)=1-g$. 

Let us denote by $e_i\in \HH$, $i=1,\ldots,g$, the natural idempotents, and let
$$\unit:=e_1+\ldots+e_g\in \HH.$$ 
Let $SG'(g)$ denote the closed subscheme of $SG(g)$ consisting of $W$ such that $\unit\in W$.
Note that $SG'(g)$ can be naturally identified with a connected component of the Sato Grassmannian 
$SG(\HH/\lan\unit\ran)$, where $\lan\unit\ran$ is the subspace generated by $\unit$. 
We are interested in the open subset $SG_1(g)\sub SG'(g)$ given by
$$SG_1(g):=\{W\in SG'(g)\ |\ H^0(W)=\lan 1\ran\}.$$
Equivalently, $W\in SG_1(g)$ if and only if $H^1(W)$ 
has minimal possible (for $W$ containing $\unit$) rank $g$.

In fact, $SG_1(g)$ is the union of the open subsets, which we call {\it open cells}, defined by
\begin{equation}\label{standard-open-eq}
SG_1(g)_L:=SG'(g)\cap U_{L'}
\end{equation} 
where $L=\lan\unit\ran\oplus L'\sub\HH$, is a subspace containing $\HH_{\ge 0}$ and such that $\dim L/\HH_{\ge 0}=g$.
The following lemma implies that this open subset does not depend on a choice of the complement $L'\sub L$ to $\lan\unit\ran$.

\begin{lem}
One has
\begin{equation}\label{standard-open2-eq}
SG_1(g)_L=\{W\in SG_1(g) \ |\ W+L=\HH\}=\{W\in SG'(g) \ |\ L/\HH_{\ge 0}\rTo{\sim} H^1(W)\}.
\end{equation}
\end{lem}

\Pf . By definition, $W\in SG'(g)$ is in $SG_1(g)_L$ iff the map $L'\simeq L/\lan\unit\ran \to \HH/W$ is an isomorphism.
This implies that the map 
\begin{equation}\label{quotient-H1-map-eq}
L/\HH_{\ge 0}\to \HH/(\HH_{\ge 0}+W)=H^1(W)
\end{equation}
is an isomorphism, or equivalently, $W\in SG_1(g)$ and the map $L/\HH_{\ge 0}\to H^1(W)$ is surjective.
The latter surjectivity is equivalent to $W+L=\HH$. Conversely, if the map \eqref{quotient-H1-map-eq}
is an isomorphism then $\dim H^1(W)=g$, so $\HH_{\ge 0}\cap W=\lan\unit\ran$, hence $\HH/W$ is an extension of $H^1(W)$ by
$\HH_{\ge 0}/\lan\unit\ran$, and so the map $L/\lan\unit\ran\to \HH/W$ is an isomorphism.
\ed

%Let $W\sub\HH$ be a complement to $L'$, such that $\unit\in W$, then we 
%We will refer to the open subsets to the 

In particular, for $\bS=\sqcup_{i=1}^nS_i$, where the subsets $S_i\sub\Z$, $i=1,\ldots,r$, are such that $\Z_{\ge 0}\sub S_i$ and $\sum_i |S_i\setminus \Z_{\ge 0}|=g$, we set
\begin{equation}\label{U-S-1-eq}
U_{\bS,1}=U_{\bS,1}(g):=SG_1(g)_{\HH_\bS}.
\end{equation}
These cells form an open covering of $SG_1(g)$ defined over $\Z$.
Note that for any subspace $W\sub \HH_R$ corresponding to a point in $U_{\bS,1}(R)$ (where $R$ is a commutative ring), 
the quotient $W/R\cdot 1$ has a canonical basis of the form
\begin{equation}\label{cell-basis-eq}
(t_i^j+v_{i,j})_{i=1,\ldots,n, j\not\in S_i} \ \ \text{ with } \ v_{i,j}\in \HH_\bS.
\end{equation} 

\begin{defi}\label{Lambda-defi}
We denote by $\sV$ the vector bundle of rank $g$ over $SG_1(g)$, 
whose fiber over $W$ is $H^1(W)$ (below we will show that $\sV$ is related to the similarly defined bundle
on a certain moduli space of curves). By definition, this bundle is trivialized over each open subset $SG_1(g)_L$ 
with $L$ as above via the natural morphism 
$$(L/\HH_{\ge 0})\ot\OO\to \sV$$
which is an isomorphism over $SG_1(g)_L$. 
Similarly, for every $N\ge 0$ we denote by $\sV_N$ the vector bundle over $SG_1(g)$ whose fiber over
$W$ is $\HH/(W+\HH_{\ge N})$. More precisely, for $N\ge 1$ it has rank $g+Nn-1$ and is defined using 
the natural trivializations 
$$(L/(\HH_{\ge N}+\lan 1\ran))\ot\OO\to \sV_N$$
over $SG_1(g)_L$.
\end{defi}

\begin{rem}\label{Det-rem}
Note that the line bundle $\det(\sV)$ is isomorphic to the inverse of the determinant bundle $\Det$ on $SG$, restricted to $SG_1(g)$. Indeed, by definition, on the open subset of $W$ such that $W\cap \HH_{\ge 1}=0$ (which contains $SG_1(g)$),
we have $\Det^{-1}_W\simeq \det(\HH/(W+\HH_{\ge 1}))$.
Since on $SG_1(g)$ we have $W\cap H_{\ge 0}=\lan 1\ran$, there is a canonical isomorphism 
$$\det(\HH/(W+\HH_{\ge 1}))\simeq \det(\HH/W+\HH_{\ge 0})=\det(\sV).$$
\end{rem}

%The Krichever map is usually defined for smooth curves (or set-theoretically for some singular curves, see 
%\cite[Sec.\  6]{SW}).
Now we are we are going to define a version of the Krichever map for the moduli spaces of curves we are interested in.

Let $\UU^{(\infty)}_{g,n}$ denote the moduli stack of projective curves $C$ with $h^0(C,\OO)=1$, of arithmetic genus $g$,
equipped with distinct smooth marked points $p_1,\ldots,p_n$ and formal parameters $(t_1,\ldots,t_n)$ at these points.
Let us consider the open substack in the stack of all curves 
$$\UU^r_{g,n}\sub \UU_{g,n}$$ 
consisting of $(C,p_1,\ldots,p_n)$ such that $h^1(C,\OO(N(p_1+\ldots+p_n)))=0$ for $N\gg 0$,
and let $\UU^{r,(\infty)}_{g,n}\sub \UU^{(\infty)}_{g,n}$ denote the preimage of $\UU^r_{g,n}$ under the projection
$\UU^{(\infty)}_{g,n}\to \UU_{g,n}$.

Below we use the extension of the standard ``cohomology and base change" results to algebraic stacks proved
in \cite{Hall-BC}.

\begin{lem}\label{R1pi*O-lem} 
Let $\pi:C\to\UU^r_{g,n}$ be the universal curve. Then $\sV:=R^1\pi_*\OO$ is a vector bundle on 
$\UU^r_{g,n}$.
\end{lem}

\Pf . Let $\UU_{g,n}(N)\sub \UU_{g,n}$ be the open substack defined by the condition
$H^1(C,\OO(N(p_1+\ldots+p_n)))=0$. 
By definition, $\UU^r_{g,n}=\cup_N \UU_{g,n}(N)$, so it is enough to prove the assertion over $\UU_{g,n}(N)$.
By the base change we know that $R^1\pi_*(\OO(N(p_1+\ldots+p_n)))=0$, 
%and the base change morphism
%$$\varphi^0(x):\pi_*(\OO(N(p_1+\ldots+p_n)))\to H^0(C_x,\OO((N(p_1+\ldots+p_n))))$$
%is an isomorphism for any point $x\in \UU_{g,n}(N)$.
so we have an exact sequence
\begin{equation}\label{R1pi*O-seq}
0\to \pi_*\OO(N(p_1+\ldots+p_n))/\pi_*\OO\to \pi_*(\OO(N(p_1+\ldots+p_n))/\OO)\rTo{\de} R^1\pi_*\OO\to 0.
\end{equation}
Furthermore, for any $x\in \UU_{g,n}(N)$ we have a commutative square
\begin{diagram}
\pi_*(\OO(N(p_1+\ldots+p_n))/\OO)\otimes k(x)&\rTo{}& R^1\pi_*\OO\otimes k(x)\\
\dTo{\varphi^0(x)}&&\dTo{\varphi^1(x)}\\
H^0(C_x,\OO(N(p_1+\ldots+p_n))/\OO)&\rTo{}& H^1(C_x,\OO)
\end{diagram}
where the vertical arrows are the base change maps and the horizontal arrows are surjective.
Since the left vertical arrow is an isomorphism, this implies that 
$$\varphi^1(x):R^1\pi_*\OO\otimes k(x)\to H^1(C_x,\OO)$$ 
is surjective.
On the other hand, since $H^0(C_x,\OO)$ is spanned by $1$, the base change map
$$\varphi^0(x):\pi_*\OO\otimes k(x)\to H^0(C_x,\OO)$$ 
is surjective for any $x$. By \cite[Thm.\ A]{Hall-BC}, this implies that $R^1\pi_*\OO$ is a vector bundle.
\ed

We have a natural action of the group $\fG=\prod_{i=1}^n\fG_i$ on $\HH$, where $\fG_i$ acts by changes of the parameter 
$t_i$ of the form
$$t_i\mapsto t_i+c_1t_i^2+c_2t_i^3+\ldots.$$
Note that the group $\fG_i$ is the projective limit of the groups $\fG_i(p)$, which only track
the changes of $t_i\mod (t_i^{p+1})$. We have $\fG_i(1)=1$ and each $\fG_i(p)$ is a unipotent algebraic group
(obtained by successive extensions of $\G_a$).
In this way we can view $\fG_i$ as a pro-unipotent group scheme defined over $\Z$.

We consider the induced action of $\fG$ on $SG_1(g)\sub SG$, which is an algebraic action of a group scheme on a scheme.
The group $\fG$ also acts naturally on $\UU^{(\infty)}_{g,n}$ by changing the formal parameters at the marked points.

\begin{prop}\label{Kr-map-prop}
There is a $\fG$-equivariant
morphism ({\it the Krichever map})
\begin{equation}\label{Krichever-map}
\Kr:\UU^{r,(\infty)}_{g,n}\to SG_1(g):
(C,p_1,\ldots,p_n, t_1,\ldots,t_n)\mapsto H^0(C\setminus\{p_1\ldots,p_n\})\sub \HH.
\end{equation}
Here the embedding $H^0(C\setminus\{p_1\ldots,p_n\})\sub \HH$ is given 
by the expansions of functions in Laurent series with respect to the parameters $t_i$.
\end{prop}

\Pf . For $N\ge 1$ let $\UU^{(\infty)}_{g,n}(N)\sub \UU^{(\infty)}_{g,n}$ be the preimage of $\UU_{g,n}(N)\sub \UU_{g,n}$.
It is enough to define compatible morphisms on each $\UU^{(\infty)}_{g,n}(N)$.
The main point is that we can define open subsets in $\UU^{(\infty)}_{g,n}(N)$ which will become the preimages
of the open cells $U_{\bS,1}\sub SG_1(g)$ 
under the map $\Kr$. Namely, taking $\bS$ to be the collection of subsets  $S_i\sub \Z$, 
$i=1,\ldots,n$, such that $\N\sub S_i\sub -N+\N$ (our convention is that $0\in\N$), we define the open
subset 
$$\UU^{(\infty)}_{g,n}(N,\bS)\sub \UU^{(\infty)}_{g,n}(N)$$
by the condition that the composition
$$\bigoplus_{i=1}^n\bigoplus_{j\in [-N,-1]\cap S_i} \OO\cdot t_i^j\to
\bigoplus_{i=1}^n\bigoplus_{j\in [-N,-1]} \OO\cdot t_i^j\simeq
\pi_*(\OO(N(p_1+\ldots+p_n))/\OO)\rTo{\de} R^1\pi_*\OO$$
%$$H^0(C,\OO(N(p_1+\ldots+p_n)))/H^0(C,\OO)\to H^0(C,\OO(N(p_1+\ldots+p_n))/\OO)\to \bigoplus_{i=1}^n
%\bigoplus_{j\in [1,N]\setminus S_i} k\cdot t_i^{-j}$$
is an isomorphism, where the isomorphism in the middle uses the formal parameters at the
marked points. Note that since $\de$ is surjective and since $R^1\pi_*\OO$ is a vector bundle by
Lemma \ref{R1pi*O-lem}, these open subsets cover $\UU^{(\infty)}_{g,n}(N)$.
%Since the spaces $H^0(C,\OO(N(p_1+\ldots+p_n)))/H^0(C,\OO)$ form a vector bundle
%over $\UU^{(\infty)}_{g,n}(N)$, one can define the above condition for any family of curves in $\UU^{(\infty)}_{g,n}(N)$.
Now we claim that over $\UU^{(\infty)}_{g,n}(N,\bS)$ the subspaces
$H^0(C\setminus \{p_1,\ldots,p_n\})\sub\HH$ give rise to points of $U_{\bS,1}$, and this defines a $\fG$-equivariant morphism
$$\UU^{(\infty)}_{g,n}(N,\bS)\to U_{\bS,1}.$$
Indeed, this follows from the fact that
for each $(C,p_1,\ldots,p_n)$ underlying a point in $\UU^{(\infty)}_{g,n}(N,\bS)$, the image of the map
$$H^0(C,\OO(N'(p_1+\ldots+p_n)))\to H^0(C,\OO(N'(p_1+\ldots+p_n))/\OO)\simeq 
\bigoplus_{i=1}^n\bigoplus_{j\in [-N',-1]} k\cdot t_i^j$$
for $N'\ge N$, is complementary to the subspace $\bigoplus_{i=1}^n\bigoplus_{j\in [-N',-1]\cap S_i} k\cdot t_i^j$. 
Furthermore, the exact sequence \eqref{R1pi*O-seq} shows that over $\UU^{(\infty)}_{g,n}(N,\bS)$ the composed map
$$\pi_*\OO(N'(p_1+\ldots+p_n))/\pi_*\OO\to \pi_*(\OO(N'(p_1+\ldots+p_n))/\OO)\to \bigoplus_{i=1}^n
\bigoplus_{j\in [-N',-1]\setminus S_i} \OO\cdot t_i^j$$
is an isomorphism for $N'\ge N$, so using the inverse map we get a map
$$\bigoplus_{i=1}^n\bigoplus_{j\in [-N',-1]\setminus S_i} \OO\cdot t_i^j\to
\pi_*\OO(N'(p_1+\ldots+p_n))/\pi_*\OO\to\bigoplus_{i=1}^n\prod_{j\in S_i\setminus\{0\}} \OO\cdot t_i^j,$$
defined over $\UU^{(\infty)}_{g,n}(N,\bS)$.
The collection of the resulting functions on $\UU^{(\infty)}_{g,n}(N,\bS)$ gives the required morphism to $U_{\bS,1}$. 
%These morphisms for different $\bS$ and $N$ are compatible, so we get the required morphism. 
\ed

\begin{rems}\label{Krich-map-rems}
1. The pull-back of the vector bundle $\sV$ (see Definition \ref{Lambda-defi})
under $\Kr$, is isomorphic to the pull-back of the similarly denoted bundle $\sV$ on $\UU^r_{g,n}$ 
(see Lemma \ref{R1pi*O-lem}).
Over the locus of stable curves the latter bundle is isomorphic to the dual of the Hodge bundle.

\noindent 2. In the case of moduli of integral curves the morphism $\Kr$ was defined in \cite[Prop.\ 6.3]{MP-PM}.
An example of $(C,p_1,\ldots,p_n)$, which does not belong to the substack $\UU^{r,(\infty)}_{g,n}$, is a curve
that has a nonrational smooth component $Z$, joined with the union of other components in a single node,
such that there are no marked points on $Z$. In fact, for such a curve 
the subspace $H^0(C\setminus \{p_1,\ldots,p_n\})\sub\HH$
does not belong to $SG(g)$.
\end{rems}

%Since all the subspaces of the form $V=H^0(C\setminus D)$ contain $\unit$,
%it is convenient to set 
%and to view $V'=V/k\unit$ as elements of the Sato Grassmannian
%of $\HH'$. More precisely, we define the component $SG(\HH',g)$ as those $V'\sub \HH'$ for which
%$V'\cap \HH'_{\ge 0}=0$ and $H^1_{V'}:=\HH'/(V'+\HH'_{\ge 0})$ is $g$-dimensional,
%where $\HH'_{\ge 0}=\HH'_{\ge 0}/k\unit$.

\subsection{Formal divisors and cohomology}\label{form-div-sec}

Let us denote for any integer vector $\bb=(b_1,\ldots,b_n)\in\Z^n$,
$$\HH_{\ge \bb}:=\bigoplus_{i=1}^n t_i^{b_i} k[[t_i]]\sub\HH, \ \ \HH_{\le \bb}:=\bigoplus_{i=1}^n t_i^{b_i} k[t_i^{-1}]\sub \HH.
$$
We define $\HH_{>\bb}$ and $\HH_{<\bb}$ in a similar fashion.
%Also, if all $b_i\le 0$ then we set $\HH'_{\ge \bb}:=\HH_{\ge \bb}/k\unit$.

We will also use formal divisors, which are elements of the free abelian group with
the basis $\bp_1,\ldots,\bp_n$ (formal points). For such divisors we write $\sum a_i\bp_i\le \sum b_i\bp_i$
if $a_i\le b_i$ for every $i$. The support $\supp(\bD)$ of $\bD=\sum a_i\bp_i$ is the set of $i$ such that $a_i\neq 0$.
We set $\deg(\bD)=\sum a_i$. For $\bD=\sum a_i\bp_i$ we set
$$\HH(\bD):=\HH_{\ge (-a_1,\ldots,-a_n)}.$$
For a point $W$ of the Sato Grassmannian we set
\begin{equation}\label{H-W-D-def-eq}
H^0(W(\bD)):=W\cap \HH(\bD), \ \ H^1(W(\bD)):=\HH/(W+\HH(\bD)).
\end{equation}

For example, the fiber of the bundle $\sV_N$ at $W$ is $H^1(W(-N(p_1+\ldots+p_n)))$. For 
$W=H^0(C\setminus\{p_1,\ldots,p_n\})$, where $(C,p_1,\ldots,p_n,t_1,\ldots,t_n)\in \UU^{r,(\infty)}_{g,n}$
one has 
$$H^i(W(a_1\bp_1+\ldots+a_n\bp_n))=H^i(C,\OO_C(a_1p_1+\ldots+a_np_n)).$$

Note that for every $\bD\ge 0$ the subset $SG(\bD)\sub SG_1(g)$ consisting of $W$ such that 
$H^1(W(\bD))=0$ is open. 
%and these open subsets cover $SG_1(g)$. 
In the case when $\deg(\bD)=g$ we get one of the open cells
\eqref{standard-open-eq}:
$$SG(\bD)=SG_1(g)_{\HH(\bD)}.$$

For every $\bD=\sum a_i\bp_i\ge 0$ we have a natural morphism of vector bundles over $SG_1(g)$,
\begin{equation}\label{pi-ba-eq}
\pi_\bD=\pi_\ba: (\HH_{\ge -\ba}/\HH_{\ge 0})\ot\OO\to \sV,
\end{equation}
induced by the embedding $\HH_{\ge -\ba}/\HH_{\ge 0}\to \HH/\HH_{\ge 0}$ and by the natural projection.

\begin{defi}\label{K-C-defi}
For a scheme $X$ with a morphism $f:X\to SG_1(g)$ and a formal divisor $\bD=\sum a_i\bp_i\ge 0$ we define
the coherent sheaves on $X$,
\begin{equation}\label{K-f-C-f-def-eq}
\KK(f,\ba)=\KK(f,\bD):=\ker(f^*\pi_\bD), \ \ \CC(f,\ba)=\CC(f,\bD):=\coker(f^*\pi_\bD).
\end{equation}
In the case when $f=\id$ we set $\KK(\bD)=\KK(\id,\bD)$, $\CC(\bD)=\CC(\id,\bD)$.
\end{defi}

For example, $\CC(0)=\sV$.

It is easy to see that for the embedding $i:\{W\}\to SG_1(g)$ of a $k$-point in $SG_1(g)$ we have
$$\KK(i,\bD)=H^0(W(\bD))/k\cdot 1, \ \ \CC(i,\bD)=H^1(W(\bD)).$$

\begin{lem}\label{H-0-H-1-exact-seq-lem} 
For every pair of formal divisors $0\le \bD'\le \bD$, and a morphism $f:X\to SG_1(g)$ one has
an exact sequence
$$0\to \KK(f,\bD')\to \KK(f,\bD)\to (\HH(\bD)/\HH(\bD'))\ot \OO_X\to \CC(f,\bD')\to \CC(f,\bD)\to 0.$$
\end{lem}

\Pf . This is the long exact sequence of cohomology associated
with the exact sequence of two-term complexes 
%\begin{equation}\label{ex-seq-2-term-com}
\begin{diagram}
0 &\rTo{}& (\HH(\bD')/\HH_{\ge 0})\ot \OO_X  &\rTo{}&  (\HH(\bD)/\HH_{\ge 0})\ot \OO_X &\rTo{}& 
(\HH(\bD)/\HH(\bD'))\ot\OO_X &\rTo{}& 0\\
&&\dTo{f^*\pi_{\bD'}}&&\dTo{f^*\pi_\bD}&&\dTo{}\\
0 &\rTo{}& f^*\sV&\rTo{\id}&  f^*\sV&\rTo{}& 0 &\rTo{}& 0
\end{diagram}
%\end{equation}
\ed

\begin{rem}\label{coherent-rem} 
It is well known that the ring of polynomials in infinitely many variables (with coefficients in $\Z$)
is coherent. It follows that the structure sheaf $\OO$ is coherent on $SG_1(g)$. Hence,
%for each $\pi$ 
%locally there exists a projection $SG_1(g)\to \A^N$ given by some of the coordinates
%on $SG_1(g)\simeq \A^\infty$ such that the morphism $\pi_D$ is the pull-back of a morphism of
%vector bundles on $\A^N$. It follows that 
for any $\bD\ge 0$ the $\OO$-modules $\KK(\bD)$ and $\CC(\bD)$ are in fact coherent sheaves on $SG_1(g)$.
\end{rem}

\begin{cor}\label{C-f-isom-cor} 
For any $0\le \bD'\le \bD$, over the open subset $SG(\bD)\sub SG_1(g)$  one has a natural isomorphism of
$\CC(\bD')$ with the cokernel of the morphism of vector bundles
$$\KK(\bD)\to (\HH(\bD)/\HH(\bD'))\ot \OO.$$
\end{cor}

\Pf . Over $SG(\bD)$ the map $\pi_\bD$ is surjective, hence $\KK(\bD)$ is a vector bundle. Now the assertion
follows from the exact sequence of Lemma \eqref{H-0-H-1-exact-seq-lem}.
\ed

\subsection{The open cell associated with $\ba=(a_1,\ldots,a_n)$}\label{cell-sec}

For $\ba=(a_1,\ldots,a_n)\in X(g,n)$ we consider the open subset
$SG^{\ba}\sub SG_1(g)$ defined by
$$SG^{\ba}:=SG(a_1\bp_1+\ldots+a_n\bp_n)=
\{W\in SG_1(g) \ |\ W+\HH_{\ge -\ba}=\HH\},$$
where $-\ba=(-a_1,\ldots,-a_n)$ (see \eqref{standard-open2-eq}).
%This is equivalent to the condition 
%that the projection $\HH_{\ge -\ba}/\HH_{\ge 0}\to \HH/W=\HH^1_V$ is surjective, hence
%an isomorphism (note that the dimension of both spaces is equal to $g$).
In other words, $SG^{\ba}$ is the locus in $SG_1(g)$, where $H^1(W(a_1\bp_1+\ldots+a_n\bp_n))=0$, or where 
%the natural morphism of vector bundles of rank $g$ over $SG_1(g)$,
$\pi_\ba$ is an isomorphism.

%Let us define the open subset $SG^{ns}=SG(\HH',g)^{ns}\sub SG(\HH',g)$ by
%$$SG^{ns}:=\{V'\in SG(\HH',g)\ |\ V'+\HH'_{\ge -1}=\HH'\}.$$
%In other words, $V'\in\UU$ if and only if $t_1^{-1},\ldots,t_g^{-1}$ project to a basis in the $g$-dimensional
%space $H^1_{V'}$. 

Note that the 
%subspace $H^0(C\setminus\{p_1,\ldots,p_g\},\OO)/k \sub\HH'$, associated with
%$(C,p_1,\ldots,p_g,t_1,\ldots,t_g)$ belongs to $SG^{ns}$ if and only if 
%$h^1(p_1+\ldots+p_g)=0$, or equivalently $h^0(p_1+\ldots+p_g)=1$.??? 
preimage of $SG^{\ba}$ under the Krichever map
\eqref{Krichever-map} is the open substack
$\UU^{(\infty)}_{g,n}(\ba)\sub \UU^{(\infty)}_{g,n}$
given by the condition $H^1(C,\OO(a_1p_1+\ldots+a_gp_g))=0$.

We identify $SG^{\ba}$ with the infinite-dimensional affine space
$\Hom_k(\HH_{<-\ba},\HH_{\ge -\ba}/\lan\unit\ran)$, by associating with $W\in SG^{\ba}$
the unique
%Let $\HH'_{\ge -\ba}\sub\HH_{\ge -\ba}$ be a complement to the one-dimensional subspace
%$k\unit\sub \HH_{\ge -\ba}$. Then $V\in SG^{\ba}$, if and only if $V$ is the graph of
linear map 
$$\phi':\HH_{<-\ba}\to \HH_{\ge -\ba}/\lan\unit\ran $$
such that $W/\lan\unit\ran$ is the graph of $\phi'$. 
We can lift $\phi'$ to a linear map 
\begin{equation}\label{phi-main-map}
\phi:\HH_{<-\ba}\to \HH_{\ge -\ba},
\end{equation}
%such that the expansion of $\phi(t_i^p)$ in $t_i$ has no constant term.???
defined up to adding a linear map with values in $\lan 1\ran\sub\HH_{\ge -\ba}$.
Then the subspace $W$ corresponding to $\phi'$ is spanned by elements $\unit, (f_i[p])_{i=1,\ldots,n; p<-a_i}$, where
\begin{equation}\label{f-i-p-eq}
f_i[p]=t_i^p+\phi(t_i^p),
\end{equation}
%Let us consider the coefficients of $\phi(t_i^p)$, for $i=1,\ldots,n$, $p<-a_i$:
with
\begin{equation}\label{coordinates-a-ij-eq}
\phi(t_i^p)=\sum_{j=1}^n \sum_{q\ge -a_j}\a_{ij}[p,q] t_j^q,
\end{equation}
%where $\wt{\a}_{ij}[p,0]$ are defined up to $\wt{\a}_{ij}[p,0]\mapsto \wt{\a}_{ij}[p,0]+c_p$.
%To get rid of this ambiguity we set
%\begin{equation}
%\a_{ij}[p,q]=\begin{cases} \wt{\a}_{ij}[p,q], & q\neq 0,\\ \wt{\a}_{ij}[p,0]-\wt{\a}_{ii}[p,0], & q=0.\end{cases}
%\end{equation}
%Note that by our choice of $\phi$ we have $\a_{ii}[p,0]=0$.
Note that the elements $f_i[p]$ are defined uniquely up to adding a constant. Thus,
the coefficients $(\a_{ij}[p,q])$, where $p<-a_i, q\ge -a_j$, are well defined for $q\neq 0$ (i.e., do not depend
on a choice of lifting $\phi$ of $\phi'$),
whereas for $q=0$ only the differences of the form $\a_{ij}[p,0]-\a_{ij'}[p,0]$ are well defined.
%The remaining coefficients $(\a_{ij}[p,q])$ (with $p<-a_i, q\ge -a_j$; $i\neq j$ if $q=0$) form a set of
We can normalize the coefficients $\a_{ij}[p,0]$ (and $f_i[p]$)
by requiring that $\a_{ii}[p,0]=0$. This is the normalization used in this
Section, as well as in Section \ref{generators-sec}.
% and \ref{g-2-n-2-a-1-1-sec}.
Another normalization, which is used in Section \ref{extending-sec}, is obtained by
choosing $j_0\in [1,n]$ and requring that $\a_{i,j_0}[p,0]=0$.
With either choice the remaining nonzero functions $(\a_{ij}[p,q])$ 
form coordinates on the affine space $SG^{\ba}$.

Using the above coordinates we identify $SG^\ba$ with the infinite-dimensional affine space over $\Z$.
For a commutative ring $R$, the $R$-points of $SG^\ba$ can be identified with $R$-submodules of
$$\HH_R:=\bigoplus_{i=1}^nR((t_i))$$
for which there exists a basis $\unit, (f_i[p])$ of the form specified above.

%Having identified the cells $SG^{\ba}$ with the affine spaces, let us look at the intersection of two cells.

%\begin{rem}
%For every subset of monomials $\bS\sub \{t_i^n \ |\ n\in\Z, 1\le i\le g\}$??? such that 
%$|S\setminus S_0|=|S_0\setminus S|<\infty$, we can define the corresponding Pl\"ucker coordinate $\si_S$ on $SG_1(g)$
%which is a section of the dual of the determinantal bundle $\Det$.
%On $SG^{ns}=U_{\bS_0,1}$ we have a trivialization of $\Det$ such that $\si_{\bS_0}=1$.
%Thus, $\si_\bS/\si_{\bS_0}$ become functions of the canonical coordinates on $SG^{ns}$, etc.
%\end{rem}

Recall that for every $N\ge 0$ we introduced a vector bundle $\sV_N$ over $SG_1(g)$ with the 
fiber $H^1(W(-N\bp_1-\ldots-N\bp_n))$ over $W$. Furthermore, over $SG^\ba$ we have a trivialization
$$\HH_{\ge -\ba}/(\HH_{\ge N}+\lan 1\ran)\ot \OO\rTo{\sim} \sV_N$$
(see Definition \ref{Lambda-defi}).

%\begin{defi}
With every set $\bS=\sqcup_{i=1}^n S_i$, 
such that  $\Z_{\ge 0}\sub S_i\sub\Z$ for $i=1,\ldots,n$ and
$\sum_{i=1}^n |S_i\setminus\Z_{\ge 0}|=g$, 
we associate a morphism of vector bundles of rank $g$ on
$SG_1(g)$,
\begin{equation}\label{pi-S-eq}
\pi_{\bS}:\HH_{\bS}/\HH_{\ge 0}\ot \OO\to \sV,
\end{equation}
and hence a global section
\begin{equation}\label{s-S-det-eq}
s_{\bS}:=\det(\pi_\bS)
\end{equation}
of the line bundle $\det(\sV)$, where we use the trivialization of $\det(\HH_{\bS}/\HH_{\ge 0})$ associated
with a fixed basis of $\HH$, given by $(t_i^m)$ in some order.

Similarly, for $N\ge 1$, we consider sets $\bS=\sqcup_{i=1}^n S_i$, 
where $\Z_{\ge N}\sub S_i\sub\Z$, such that 
$\sum_{i=1}|S_i\setminus\Z_{\ge N}|=g+Nn-1$,
and the corresponding morphism of vector bundles of rank $g+Nn-1$ on
$SG_1(g)$,
\begin{equation}\label{pi-S-N-eq}
\pi_{\bS,N}:\HH_{\bS}/\HH_{\ge N}\ot \OO\to \sV_N.
\end{equation}
We set 
$$s_{\bS,N}:=\det(\pi_{\bS,N})\in H^0(SG_1(g),\det(\sV_N)).$$
%where we use the trivialization of ???

\begin{lem}\label{Plucker-lem}
(i) There is a natural isomorphism $\kappa_N:\det(\sV_N)\to\det(\sV)$, such that
$$s_{\bS,1}:=\kappa_N(s_{\bS,N})$$ 
does not depend on $N\ge 1$, 
for appropriate ordering of the standard bases of $\HH_{\bS}/\HH_{\ge N}$.

\noindent
(ii) Let $\bS=\sqcup_{i=1}^n S_i$ be such that $\Z_{\ge 0}\sub S_i$ and 
$\sum_{i=1}^n |S_i\setminus \Z_{\ge 0}|=g$.
For fixed $i$ consider $\bS'=\sqcup_{j=1}^n S'_j$, where
$S'_i=S_i\setminus\{0\}$ and $S'_j=S_j$ for $j\neq i$.
Then
$$s_{\bS}=s_{\bS',1},$$
for appropriate ordering of the standard basis of $\HH_{\bS}/\HH_{\ge 0}$.
%where the left hand side is given by \eqref{s-S-det-eq}.
\end{lem}

\Pf . (i) For $N\ge 1$ we have an exact sequence of vector bundles on $SG_1(g)$
\begin{equation}\label{La-N-La-seq}
0\to \HH_{\ge 0}/(\HH_{\ge N}+\lan 1\ran)\ot\OO\to \sV_N\to \sV\to 0,
\end{equation}
Passing to determinants we get the isomorphism $\kappa_N$. If $\Z_{\ge N}\sub S_i$ for all $i$ then
we have a morphism of exact sequences
\begin{diagram}
0\to &\HH_{\ge N}/\HH_{\ge N+1}\ot \OO &\rTo{}&\HH_{\bS}/\HH_{\ge N+1}\ot \OO &\rTo{}& \HH_{\bS}/\HH_{\ge N}\ot \OO &\to 0\\
&\dTo{}&&\dTo{\pi_{\bS,N+1}}&&\dTo{\pi_{\bS,N}}\\
0\to &\HH_{\ge N}/\HH_{\ge N+1}\ot \OO &\rTo{}&\sV_{N+1}&\rTo{}&\sV_N&\to 0
\end{diagram}
The lower exact sequence gives an isomorphism $\kappa_{N,N+1}:\det(\sV_N)\to\det(\sV_{N+1})$, so that
$$\kappa_{N,N+1}(\det(\pi_{\bS,N}))=\det(\pi_{\bS,N+1}).$$
It is easy to see that $\kappa_{N+1}\circ\kappa_{N,N+1}=\kappa_N$, so we deduce that
$$\kappa_N(\det(\pi_{\bS,N}))=\kappa_{N+1}(\det(\pi_{\bS,N+1})).$$

(ii) This is proved similarly to (i), using the exact sequence 
$$0\to \HH_{\ge 0}/(\HH_{\ge 1}+\lan 1\ran)\to \sV_1\to \sV_0\to 0.$$
\ed

We refer to $(s_\bS)$ and $(s_{\bS,1})$ as {\it Pl\"ucker coordinates} on $SG_1(g)$.
Since we do not fix an ordering of the standard bases of $\HH_{\bS}/\HH_{\ge N}$,
there is a sign ambiguity in the definition of these sections, hence the appearance of $\pm$ in the formulas below.

%\end{defi}

\begin{rem}
Under the identification of $\det(\sV)$ with $\Det^{-1}$ (see Remark \ref{Det-rem}) our sections $(s_\bS)$ 
correspond to the restrictions of the Pl\"ucker coordinates on $SG$ (see \cite[Sec.\ 2.D]{MP-PM}).
\end{rem}

Note that for $\ba=(a_1,\ldots,a_n)\in X(g,n)$, if we set 
$S_i=[-a_i,+\infty)$, $i=1,\ldots,n$, then the morphism $\pi_\bS$ is exactly the morphism 
%we have a canonical morphism of vector bundles of rank $g$,
$\pi_\ba$ (see \eqref{pi-ba-eq}), such that
$SG^\ba$ is the locus where $\pi_\ba$ is an isomorphism. In other words,
$SG^\ba$ is the complement in $SG_1(g)$
to the closed subscheme $Z_\ba\sub SG_1(g)$ given as the zero locus of
$$s_\ba:=\det(\pi_\ba)$$ 
which is a section of the line bundle $\det(\sV)$.

\begin{prop}\label{Plucker-prop} Let us normalize $(\a_{ij}[p,0])$ by setting $\a_{ii}[p,0]=0$.

\noindent
(i) The coordinate $\a_{ij}[p,q]$, where $p\le -a_i-1$, $q\ge -a_j$ (see \eqref{coordinates-a-ij-eq}) on 
$SG^\ba$ has the following expression in terms of the Pl\"ucker coordinates:
%(and hence on $\wt{\UU}^{ns}_{g,n}(\ba)\simeq ASG^\ba\cap\Si^{\ba,j}$)
$$\a_{ij}[p,q]=\pm\frac{s_{\bS,1}}{s_\ba},$$ 
where $\bS=\sqcup S_k$ is defined by
%as follows. 
%If $q<0$ then we set 
%$S_j=[-a_j,+\infty)\setminus \{q\}$, $S_i=[-a_i,+\infty)\cup\{p\}$, and
%$S_k=[-a_k,+\infty)$ for $k\neq i,j$. 
%If $0\le q<N$ then 
$$S_i=[-a_i,+\infty)\setminus\{0\}\cup\{p\}, \ \ S_j=[-a_j,+\infty)\setminus \{q\}, \ \ S_k=[-a_k,+\infty) \text{ for }k\neq i,j.$$
%on the Sato Grassmannian (note that the determinant bundle is naturally trivialized on $SG^\ba$).
If in addition $q<0$, then we have
$$\a_{ij}[p,q]=\pm\frac{s_{\bS'}}{s_\ba},$$ 
where $\bS'=\sqcup S'_k$ with $S'_i=[-a_i,+\infty)\cup\{p\}$ and $S'_k=S_k$ for $k\neq i$. 

\noindent
(ii) For another element $\ba'\in X(g,n)$, 
the intersection $SG^{\ba}\cap SG^{\ba'}$ is the principal affine open in $SG^{\ba}$
given by the nonvanishing of the function
$$\frac{s_{\ba'}}{s_{\ba}}=\det(A_{\ba,\ba'})$$
on $SG^\ba$, where $A_{\ba,\ba'}=\pi_\ba^{-1}\circ \pi_{\ba'}$ can be viewed as
the $\Hom(\HH_{\ge -\ba'}/\HH_{\ge 0},\HH_{\ge -\ba}/\HH_{\ge 0})$-valued function
on $SG^{\ba}$ given by
\begin{equation}\label{A-matrix-formula-eq}
A_{\ba,\ba'}(\phi)(t_i^j)=\begin{cases} t_i^j \mod\HH_{\ge 0}, & j\ge -a_i, \\ -\phi(t_i^j) \mod\HH_{\ge 0} & j<-a_i,\end{cases}
\end{equation}
where we take $(t_i^j)_{i=1,\ldots,g; -a'_i\le j<0}$ as a basis of $\HH_{\ge -\ba'}/\HH_{\ge 0}$ and identify
$SG^{\ba}$ with the affine space of maps $\phi$ as in \eqref{phi-main-map}.
\end{prop}

\Pf . (i) Set $N=\max(1,q+1)$, and let
$\HH'_{\ge -\ba}\sub\HH_{\ge -\ba}$ (resp., $\HH'_{\ge 0}\sub \HH_{\ge 0}$)
be the subspace obtained by omitting the element $t_i^0$ in the standard basis
of $\HH_{\ge -\ba}$ (resp., $\HH_{\ge 0}$). 
As in Lemma \ref{Plucker-lem}, using the exact sequence \eqref{La-N-La-seq}, and the morphism
of exact sequences
\begin{diagram}
0\rTo{} &\HH'_{\ge 0}/\HH_{\ge N}\ot \OO &\rTo{}&\HH'_{\ge -\ba}/\HH_{\ge N}\ot \OO &\rTo{}& \HH_{\ge -\ba}/\HH_{\ge 0}\ot \OO &
\rTo{} 0\\
&\dTo{\simeq }&&\dTo{\pi'_{\ba,N}}&&\dTo{\pi_{\ba}}\\
0\rTo{} &\HH_{\ge 0}/(\HH_{\ge N}+\lan 1\ran)\ot \OO &\rTo{}&\sV_{N}&\rTo{}&\sV&\rTo{} 0
\end{diagram}
we see that the natural map
%$$\HH_{\ge -\ba}/(\HH_{\ge N}+\lan 1\ran)\to \sV_N$$
%is an isomorphism, and its determinant is given by $s_\ba$ (see ???).
%Then the above map can be identified with the map
$$\pi'_{\ba,N}:\HH'_{\ge -\ba}/\HH_{\ge N}\to \sV_N$$ 
is an isomorphism on $SG^\ba$, and its determinant is equal to $s_\ba$.
Thus, $s_{\bS,1}/s_\ba$ is equal to the determinant of the map
$$A=(\pi'_{\ba,N})^{-1}\circ\pi_{\bS,N}:\HH_{\bS}/\HH_{\ge N}\ot\OO\to \HH'_{\ge -\ba}/\HH_{\ge N}\ot\OO.$$
For any element $b$ of the standard basis of $\HH_{\bS}/\HH_{\ge N}$, except
for $t_i^p$, we have $A(b)=b$, viewed as an element of $\HH'_{\ge -\ba}/\HH_{\ge N}$.
As for $A(t_i^p)$, we have to find an element $v(W)\in \HH'_{\ge -\ba}/\HH_{\ge N}$ such that
$$t_i^p\equiv v(W) \mod W+\HH_{\ge N}.$$
The element $f_i[p]\in W$ (see \eqref{f-i-p-eq}, \eqref{coordinates-a-ij-eq}), normalized by $\a_{ii}[p,0]=0$, satisfies
$$f_i[p]\equiv t_i^p \mod \HH'_{\ge -\ba}.$$
Thus, we have
$$A(t_i^p)=v(W)=t_i^p-f_i[p] \mod \HH_{\ge N}.$$
Now computing the determinant reduces to taking the coefficient of $t_j^q$ in $A(t_i^p)$.

The second formula in the case $q<0$ follows from Lemma \ref{Plucker-lem}(ii).

\noindent
(ii) The first assertion is an immediate consequence of the fact that $SG^{\ba'}$ is the nonvanishing locus of $s_{\ba'}$.
To calculate
$$A_{\ba,\ba'}: (\HH_{\ge -\ba'}/\HH_{\ge 0})\ot\OO\rTo{\pi_{\ba'}} \sV\rTo{\pi^{-1}_\ba} (\HH_{\ge -\ba}/\HH_{\ge 0})\ot\OO$$
we note that $A(t_i^j)$ at $W\in SG^{\ba}$ 
is the projection of $t_i^j\mod\HH_{\ge 0}$ to
$\HH_{\ge -\ba}/\HH_{\ge 0}$ along $W/\lan\unit\ran$. Taking $W$ to be the graph of $\phi$ we get the formula 
\eqref{A-matrix-formula-eq}.
\ed

\begin{cor}\label{a-simple-for-cor} 
If $a_j>0$ then one has 
$$\a_{ij}[-a_i-1,-a_j]=\pm \frac{s_{\ba+\be_i-\be_j}}{s_\ba}.$$
\end{cor} 

%Show that all the coordinates on $SG^\ba$ are of the form $\frac{s}{s_\ba}$ for some sections $s\in \det(\sV)$.

\subsection{Action by changes of parameters}\label{changes-parameter-sec}

Recall that we denote by $\fG$ the product over $i=1,\ldots,n$ of the groups $\fG_i$ of formal changes of parameters
of the form $t_i\mapsto t_i+c_1t_i^2+c_2t_i^3+\ldots$.
Note that the action of $\fG$ on $\HH$ preserves all the subspaces $\HH_{\ge -\bb}$ (where $\bb\in\Z^n$).
Hence, the induced action on $SG_1(g)$ preserves each open cell $SG^{\ba}$.

\begin{defi}
Let $G$ be a group scheme acting on a scheme $X$.
We say that a closed subscheme $S\sub X$ is {\it a section for the action of $G$} if the map
$G\times S\to X$, given by the action, is an isomorphism. In this case the action of $G$ on $X$ is free.
\end{defi}

The following simple observation will be useful for us.

\begin{lem}\label{section-lem} 
Let $S\sub X$ be a section for an action of $G$ on $X$ and let $f:Y\to X$ be a $G$-equivariant morphism.
Then the schematic preimage $f^{-1}(S)$ is a section for the action of $G$ on $Y$.
\end{lem}

\Pf . Set $T=f^{-1}(S)$.              Consider the commutative diagram
\begin{diagram}
G\times T &\rTo{}& G\times Y &\rTo{a_Y}& Y\\
\dTo{}&&\dTo{}&&\dTo{}\\
G\times S&\rTo{}& G\times X &\rTo{a_X}& X\\
\end{diagram}
with the vertical arrows induced by $f$ and $a_X, a_Y$ the action morphisms.
The left square is cartesian by the construction of $T$. We claim that the right square is also cartesian.
Indeed, this immediately follows from the fact that the composition of $a_X$ with the automorphism $(g,x)\mapsto (g,g^{-1}x)$
of $G\times X$ is simply the projection $G\times X\to X$ (and similarly for $Y$).
Hence, the big rectangle in the above diagram is cartesian. Since the map $G\times S\to X$ is an isomorphism,
we deduce that the map $G\times T\to Y$ is also an isomorphism.
\ed

The following result is a generalization of \cite[Lem.\ 4.1.3]{FP}. We use the notation from Section \ref{form-div-sec}.

\begin{prop}\label{section-prop} 
Fix $i$, $1\le i\le n$, and $\bb=(b_1,\ldots,b_n)\in\Z_{\ge 0}^n$, such that $b_i>0$.
Let $X$ be a scheme over $\Q$ equipped with an action of $\fG_i$, $f:X\to SG_1(g)$ a $\fG_i$-equivariant morphism, 
such that $f^*\pi_\bb$ is surjective. Assume also that 
the $\LL=\CC(f,\bb-\be_i)$ is locally free of rank $1$.
Then there is a closed subscheme $\Sigma_i^X\sub X$ parametrizing
$x\in X$ such that for each $p\ge 1$, there exists an element $v_p\in W_x:=f(x)$ with 
$$v_p\equiv t_i^{-b_i-p} \mod \HH_{\ge -\bb+\be_i},$$ 
such that $\Sigma_i^X$ is a section for the $\fG_i$-action on $X$. 
\end{prop}\

\Pf . Let us set for each $p\ge -1$,
$$\KK_p:=\KK(f,\bb+p\be_i)=\ker(f^*\pi_{\bb+p\be_i}).$$
By assumption, the morphism $f^*\pi_{\bb}$ is surjective, so from Lemma \ref{H-0-H-1-exact-seq-lem} we get
an isomorphism
\begin{equation}\label{Kp-K0-quot-isom}
\KK_p/\KK_0\simeq (\HH_{\ge -\bb-p\be_i}/\HH_{\ge -\bb})\ot\OO_X.
\end{equation}
On the other hand, from the same Lemma we get an exact sequence
%$$0\to [(\HH_{\ge -\bb+\be_i}/\HH_{\ge 0})\ot \OO_T\to f^*\sV]\to [(\HH_{\ge -\bb}/\HH_{\ge 0})\ot \OO_T\to f^*\sV]\to
%[(\HH_{\ge -\bb}/\HH_{\ge -\bb+\be_i})\ot\OO_T\to 0]\to 0$$
\begin{equation}\label{K-L-long-ex-seq}
0\to \KK_{-1}\to \KK_0\to (\HH_{\ge -\bb}/\HH_{\ge -\bb+\be_i})\ot\OO_X\to\LL\to 0.
\end{equation}
Since $\LL$ is locally free of rank $1$, this implies that the arrow $\KK_{-1}\to \KK_0$ is an isomorphism.
Thus, using \eqref{Kp-K0-quot-isom} we obtain a morphism
$$\si_p:(\HH_{\ge -\bb-p\be_i}/\HH_{\ge -\bb})\ot\OO_X\simeq \KK_p/\KK_0\simeq \KK_p/\KK_{-1}\to 
(\HH_{\ge -\bb-p\be_i}/\HH_{\ge -\bb+\be_i})\ot\OO_X,$$
where the last arrow is induced by the embedding $\KK_p\to (\HH_{\ge -\bb-p\be_i}/\HH_{\ge 0})\ot\OO_X$. 
Using the natural basis of $\HH$ we define the composed map
$$t_i^{-b_i-p}\ot\OO_X\to (\HH_{\ge -\bb-p\be_i}/\HH_{\ge -\bb})\ot\OO_X\rTo{\si_p}
(\HH_{\ge -\bb-p\be_i}/\HH_{\ge -\bb+\be_i})\ot\OO_X\to t_i^{-b_i}\ot \OO_X,$$
and we denote by $a_i(p)$ the corresponding function on $X$.
%\a_{ii}[-b_i-p,-b_i]

The meaning of this function is the following: $a_i(p)(x)$ is the coefficient of $t_i^{-b_i}$ in the expansion of an element
$v_p\in W_x$ such that $v_p\equiv t_i^{-b_i-p}\mod \HH_{\ge -\bb}$. The point is that such an element $v_p$
is defined up to adding an element in $W_x\cap \HH_{\ge -\bb}=W_x\cap \HH_{\ge -\bb+\be_i}$, so this coefficient is
well defined.

We define the subscheme $\Sigma=\Sigma_i^X\sub X$ as the common zero locus of the functions $a_i(p)$ for $p\ge 1$.
We are going to check that it is a section for the action of $\fG_i$.

Let $\fG_i(>p):=\ker(\fG_i\to \fG_i(p))$. Also, for $p\ge 1$ let
$\Sigma_p\sub X$ denote the common zero locus of the functions $a_i(p')$ for $1\le p'\le p-1$,
%=\{t \ |\ p_{ii}[m](V)=0, \ 2\le m\le n \},$$
so that $\Sigma_1=X$ and $\Sigma=\cap_p \Sigma_p$.

We have a decomposition into a semi-direct product
$$\fG_i(>p)=\fG_i(>p+1)\rtimes F_{p+1},$$
where $F_{p+1}=\G_a$ acts by $t_i\mapsto t_i+ct_i^{p+1}$.  

We claim that $\fG_i(>p)(\Sigma_p)\sub \Sigma_p$,  
the action of $F_{p+1}$ on $\Sigma_p$ is free and $\Sigma_{p+1}$ is a section for this $F_{p+1}$-action.
Indeed, note that under the change of variable
$$g:t_i\mapsto t_i+ct_i^{p+1}+\ldots$$
we have
$$t_i^{-b_i-p'}\mapsto t_i^{-b_i-p'}-(b_i+p')ct_i^{-b_i+p-p'}+\ldots.$$
Thus, under the action of the above element of $\fG_i(>p)$ the functions $a_i(p')$ with $p'<p$ do not change, while
the function $a_i(p)$ transforms by 
\begin{equation}\label{a-i-p-transformation-eq}
a_i(p)(gx)=a_i(p)(x)-(b_i+p)c.
\end{equation}
This immediately implies our claim. 
%$$p_{ii}[n+1]\mapsto p_{ii}[n+1]-(n+1)c_i.$$

Let us define for $p\ge 1$ the morphisms $g_p:X\to F_p$ (where $F_1=0$) and $\rho_p:X\to \Sigma_p$ inductively as follows.
We set $g_1=1$ and $\rho_1=\id_X$. Assuming that $g_{p-1}$ and $\rho_{p-1}$ are already defined, we define 
$g_p(x)$ and $\rho_p$ by 
$$g_p(x): t_i\mapsto t_i+\frac{a_i(p-1)(\rho_{p-1}(x))}{b_i+p-1}t_i^p,$$ 
%t_i\mapsto t_i+\frac{p_{ii}[-n](V_{n-1})}{n}t_i^n, \ \ i=1,\dots,g,$$
$$\rho_p:X\to X: x\mapsto g_p(x)\cdot \rho_{p-1}(x).$$
%$V_n=g_n(V_{n-1})$. 
Let us check that $\rho_p$ factors through the subscheme $\Sigma_p\sub X$.
By induction we can assume that $\rho_{p-1}$ factors through $\Sigma_{p-1}$. Since $F_p$ preserves $\Sigma_{p-1}$,
this implies that $\rho_p$ factors through $\Sigma_{p-1}$. Now using \eqref{a-i-p-transformation-eq} we see
that 
$$a_i(p-1)(\rho_p(x))=a_i(p-1)(g_p(x)\cdot \rho_{p-1}(x))=0,$$ 
hence, $\rho_p$ factors through $\Sigma_p$.

We claim that any $g\in\fG$, such that $g\cdot x\in\Sigma_p$, satisfies 
$$g\equiv g_p(x)\ldots g_2(x)\mod\fG_i(>p).$$
To see this we use the induction on $p$. Suppose this is true for $p$, and assume $g$ is such that $g\cdot x\in\Sigma_{p+1}$.
Replacing $x$ by $\rho_p(x)=g_p(x)\ldots g_2(x)\cdot x$ and $g$ by $gg_2(x)^{-1}\ldots g_p(x)^{-1}$ we can assume that
$x\in\Sigma_p$. By the induction assumption, this implies that $g\in \fG_i(>p)$, and the assertion follows from the fact that
$\Sigma_{p+1}$ is a section for the free action of $F_{p+1}$ on $\Sigma_p$.

Now we observe that the infinite product $g(x)=\ldots g_4(x)g_3(x)g_2(x)$ in $\fG_i$ is well defined. Furthermore, for
each $p\ge 1$ we have
$$g(x)\cdot x=\ldots g_{p+2}(x)g_{p+1}(x)\cdot \rho_p(x)\in \Sigma_p,$$
hence, $g(x)\cdot x\in\Sigma$. 
Together with the above claim this implies that $\Sigma$ is a section for the $\fG_i$-action on $X$.
\ed

%\begin{prop} The action of $\fG$ on $SG^{ns}$ is free, and the subvariety
%$\Sigma\sub SG^{ns}$ defined by equations $p_{ii}[-n]=0$ for $i=1,\ldots,g$, $n\ge 2$,
%is a section for this action.
%\end{prop}
%\Pf . Given a subspace $V=\lan f_i[-n], n\ge 2, i=1,\ldots,g\ran\in SG^{ns}$ where $f_i[-n]$ are given by ???
%Let us write $p_{ii}[-n]=p_{ii}[-n](V)$ to remember the dependence on $V$.

We can apply the above construction with $X$ being the open cell $SG^\ba$ (and $f$ the natural inclusion).

\begin{cor}\label{section-open-cor} Let us work over $\Q$.
For $\ba=(a_1,\ldots,a_n)\in X(g,n)$ let $I_\ba\sub\{1,\ldots,n\}$
be the set of $i$ such that $a_i>0$. Define
$$\Sigma^{\ba}\sub SG^\ba$$ 
as the closed subscheme cut out by the equations 
\begin{equation}\label{Sigma-a-eq}
\a_{ii}[-a_i-p,-a_i]=0 
\end{equation}
for all $p\ge 1, i\in I_\ba$.
Then
$\Sigma^\ba$, which is isomorphic to an infinite-dimensional affine space, 
is a section for the action of the group $\prod_{i\in I_\ba}\fG_i$ on $SG^\ba$.
\end{cor}

\Pf . We take $\bb=\ba$ in Proposition \ref{section-prop}. 
Note that $\pi_\ba$ is an isomorphism  over $SG^\ba$ by definition. Now the exact sequence
\eqref{K-L-long-ex-seq} for $i\in I_\ba$ gives an isomorphism
$$\HH_{\ge -\ba}/\HH_{\ge -\ba+\be_i}\ot\OO_X\simeq\coker(\pi_{\ba-\be_i})$$
over $SG^\ba$. Thus, we see that Proposition \ref{section-prop} is applicable to the 
action of $\fG_i$ on $SG^\ba$, and we get that the subscheme
$\Sigma^\ba_i\sub SG^\ba$ cut out by the equations \eqref{Sigma-a-eq} for all $p\ge 1$ and given $i\in I_\ba$
is a section for the action of $\fG_i$ on $SG^\ba$. Since each $\Sigma^\ba_i$ is invariant with respect
to the action of $\fG_j$, for $j\neq i$, using Lemma \ref{section-lem} we derive that
$$\Sigma^\ba=\cap_{i\in I_a} \Sigma^\ba_i$$
is a section for the action of $\prod_{i\in I_\ba}\fG_i$ on $SG^\ba$.
\ed

In the case when some of coordinates of $\ba$ are zero (but not all, since we assume that $g\ge 1$), 
we still get a section for the action of the full group $\fG$
on $SG^\ba$ as follows.

\begin{prop}\label{section-open-prop} 
Keep the assumptions and the notations of Corollary \ref{section-open-cor}, and let us fix some $i_0\in I_\ba$. 
%Assume $n\ge 1$ and let us fix $j_0$
Then the closed subscheme
$\Sigma^{\ba,i_0}\sub \Sigma^\ba$, cut out by the equations 
$$\a_{i,i_0}[-1,0]-\a_{ii}[-1,0]=\a_{ii}[-1,p]=0 \text{ for all } i\not\in I_\ba, p\ge 1,$$
is a section for the $\fG$-action on $SG^\ba$.
\end{prop}

\Pf . Note that for each $i\not\in I_\ba$ the change of the parameter $t_i$ of the form
$$t_i\mapsto t_i+ct_i^{p+1}+\ldots$$
leads to the transformation
$$t_i^{-1}\mapsto t_i^{-1}-ct_i^{-1+p}+\ldots.$$
Hence, for $p>1$ it does not change the coordinates $\a_{ij}[-1,0]$, $\a_{ii}[-1,-1+p']$, $1<p'<p$ and
shifts $\a_{ii}[-1,-1+p]$ by $c$. Also, for $p=1$ this transformation shifts $\a_{ij}[-1,0]$ by $-c$.
Now we can use the same argument as in the proof of Proposition \ref{section-prop}. 
%For each $i$ such that $a_i=0$, the long exact sequence associated with 
%the exact sequence of complexes \eqref{ex-seq-2-term-com} for $\bb=\ba+\be_i$
%implies that the composed arrow
%$$\ker(f^*\pi_{\ba+\be_i})\to (\HH_{\ge -\ba-\be_i}/\HH_{\ge 0})\ot \OO_X\to (\HH_{\ge -\ba-\be_i}/\HH_{\ge -\ba})\ot\OO_X$$
%is an isomorphism. In other words, for each such $i$ we have a canonical section
%$v(i)\in \ker(f^*\pi_{\ba+\be_i})$ such that $v(i)\equiv t_i^{-1}\mod \HH_{\ge -\ba}$.
%Let us first explain the remainder of the argument ``pointwise".
%For each $V\in \wt{\Sigma}^\ba$ we can lift $v(i)(V)$ to an element $\wt{v(i)}(V)\in V$, defined uniquely up to adding a constant.
%We fix this ambiguity by requiring that the expansion of $\wt{v(i)}(V)$ in $t_j$ has no constant term.
%Now there is a unique $g_i\in \fG_i$ such that the
%change of parameter $t_i\mapsto g_i(t_i)$ leads to the transformation $t_i^{-1}\mapsto \wt{v(i)}(V)(t_i)$
%???
\ed

Let $ASG=ASG(\HH,g)$ be the closed subscheme of $SG_1(g)$ consisting of $W$ such that
$W\cdot W\sub W$ (so that $W$ is a subalgebra of $\HH$). Note that the Krichever map lands in $ASG$.
For $\ba=(a_1,\ldots,a_n)\in X(g,n)$ let us set 
$$ASG^{\ba}:=SG^{\ba}\cap ASG.$$
Note that the subschemes $ASG$ and $ASG^{\ba}$ are preserved by the action of the group $\fG$. Hence,
we derive the following Corollary from Proposition \ref{section-open-prop}

\begin{cor}\label{ASG-sec-cor} 
The closed subscheme $\Sigma^{\ba,i_0}\cap ASG\sub ASG^\ba$ is a section for the $\fG$-action on
$ASG^\ba$ (over $\Q$).
\end{cor}

The Krichever map \eqref{Krichever-map} induces a morphism
$$\wt{\UU}^{ns,(\infty)}_{g,n}(\ba)\rTo{\Kr} ASG^\ba,$$
compatible with the action of the group $\fG$. 
Hence, using Lemma \ref{section-lem} we derive the following result, generalizing \cite[Lem.\ 2.1.1]{P-ainf}.

\begin{cor}\label{curves-sec-cor}
Let us work over $\Q$.
The closed subscheme $\Kr^{-1}(ASG\cap\Sigma^{\ba,i_0})$ is a section for the $\fG$-action
on $\wt{\UU}^{ns,(\infty)}_{g,n}$, i.e., it is a section for
the natural projection $\wt{\UU}^{ns,(\infty)}_{g,n}(\ba)\to \wt{\UU}^{ns}_{g,n}(\ba)$.
Hence, 
%$$\Kr^{-1}(ASG\cap\Sigma^{\ba,i_0})\simeq \wt{\UU}^{ns,(\infty)}_{g,n}(\ba)/\fG\simeq \wt{\UU}^{ns}_{g,n}(\ba)$$
we have canonical formal parameters at the marked points on the universal curve $C$
over $\wt{\UU}^{ns}_{g,n}(\ba)$.
\end{cor}

%The projection $\pi:SG^{ns}\to \Sigma$ is given by some algebraic map of coordinates. Explicit partial formula??? 

\subsection{Marked and weakly marked algebras}

In \cite{P-ainf} we have shown that $\wt{\UU}^{ns}_{g,g}(1,\ldots,1)$ is isomorphic to the moduli space of
{\it marked algebras}, which are algebras equipped with a filtration with some special properties.
Here we generalize the concept of marked algebras to the case of any weight $\ba\in X(g,n)$
and show that their moduli space has a natural interpretation in terms of the Sato Grassmannian.

Let $M_\ba\sub (\Z_{\ge 0})^n$ be the submonoid consisting of $\bm\in(\Z_{\ge 0})^n$ such that either 
$\bm\ge\ba$ or $\bm=(0,\ldots,0)$.
%and let us consider
%the submonoid $M_\ba^+$ of $(\Z_{\ge 0})^n$ given by
%$$M_\ba=M_\ba^+\cup\{(0,\ldots,0)\}.$$ 
%Sometimes we denote the elements of $M_\ba$ as formal divisors $m_1\bp_1+\ldots+m_n\bp_n$.

For a commutative ring $R$ we define an $R$-subalgebra
$$C_N\sub R[u_1]\oplus\ldots\oplus R[u_n]$$
as follows. Let $e_1,\ldots,e_n$ be the natural idempotents in $R[u_1]\oplus\ldots\oplus R[u_n]$.
Then $C_N$ is spanned as an $R$-module by $1=e_1+\ldots+e_n$ and by the elements
$u_i^{m}\in R[u_i]$ for $i=1,\ldots,n$, $m>N$.
We view $C_N$ as a graded $R$-algebra, where $\deg(u_i)=1$ for each $i$.

%$$C_\ba\sub R[t_1^{-1}]\oplus\ldots\oplus R[t_n^{-1}]$$
%as follows. Let $e_1,\ldots,e_n$ be the natural idempotents in $R[t_1^{-1}]\oplus\ldots\oplus R[t_n^{-1}]$.
%Then $C_\ba$ is spanned as $R$-module by $e_1+\ldots+e_n$ and by the elements
%$t_i^{-m}e_i$ for $i=1,\ldots,n$, $m>a_i$.

\begin{defi}\label{marked-alg-def}
(i) A {\it marked algebra of type $\ba$ over $R$} is a commutative $R$-algebra $A$ equipped with an exhaustive 
$M_\ba$-valued algebra filtration $(F_{\bm}A)$, such that $F_\ba=F_0=R$ and the map
$\bm\mapsto F_\bm A$ is a morphism of lattices, i.e.,
\begin{equation}\label{lattice-eq}
F_{\min(\bm,\bm')}A=F_\bm A\cap F_{\bm'} A, \ F_{\max(\bm,\bm')}A=F_{\bm} A+F_{\bm'} A.
\end{equation}
In addition, for each $\bn\ge\ba$ and each $i=1,\ldots,n$, there should be a fixed isomorphism
$$\varphi_{\bn,i}:F_{\bn+\be_i}A/F_{\bn}A\rTo{\sim} R,$$
%where $(\be_i)$ is the standard basis in $\Z^n$,
such that the maps
$$F_{\bn+\be_i}A/F_{\bn}A\ot F_{\bn'+\be_i}A/F_{\bn'}A\to F_{\bn+\bn'+2\be_i}A/F_{\bn+\bn'+\be_i}A,$$\
induced by the multiplication on $A$, get identified with the multiplication on $R$. 
%for the $\N$-valued filtration 
%$$F_mA=\sum_{\bm\in S_\ba, \bm\le (m,\ldots,m)} F_\bm A$$
%one has an isomorphism 
%$$\gr^\bullet_{F} A\simeq C_\ba$$
%of graded $R$-algebras.

\noindent (ii) Let us set $\bD_1:=\be_1+\ldots+\be_n\in\Z^n$, and for $N\ge 0$ let us consider the subsemigroup
$M_N\sub (\Z_{\ge 0})^n$ given by
$$M_N:=\{k\bD_1+\sum_{i=1}^n m_i\be_i \ |\ k\ge N, \ 0\le m_i\le \frac{k}{N} \text{ for } i=1,\ldots,n\}.$$
A {\it weakly $N$-marked algebra} over $R$ is a commutative $R$-algebra $A$
equipped with $R$-submodules $F_{\bm}A\sub A$, for $\bm\in M_N$,
such that $F_\bm A\cdot F_{\bm'}A\sub F_{\bm+\bm'}A$ and the lattice condition \eqref{lattice-eq} is satisfied.
We denote by $(F_m A)_{m\ge N}$ the induced filtration $F_m A:=F_{m\bD_1}A$. In addition, there should be
fixed an isomorphism of graded non-unital $R$-algebras
\begin{equation}\label{marking-2-eq}
\gr^F_{>N} A=\bigoplus_{m>N} F_mA/F_{m-1}A\simeq C_N.
\end{equation}
We impose the following two conditions on these data.
First, we require that for each $m\ge N$ the image of $F_{m\bD_1+\be_i}A/F_{m}A$ in 
$F_{m+1}A/F_mA$
gets identified with $Ru_i^{m+1}\sub (C_N)_{m+1}$.
Secondly, for each $m>N$ let us denote by $S_m$
the set of elements $s\in F_mA$, such that 
$$\ov{s}=s\mod F_{m-1}A\in \gr^F_m A\sub R^n$$
has components that are invertible in $R$. Then we require that if for some $m_0>N$ and some
$x\in A$, for each $m>m_0$ there exists $s\in S_m$ such that $xs\in F_{m_0}A$, then $x=0$.
Using the fact that $S_m\cdot S_{m'}\sub S_{m+m'}$, one can easily deduce from 
this that each $S_m$ consists of non-zero-divisors in $A$. 

A {\it weakly marked algebra} is a weakly $N$-marked algebra for some $N$.

%together with a fixed injective homomorphism
%\begin{equation}\label{marking-2-eq}
%\gr^\bullet_{F} A\to \bigoplus_{i=1}^n R[u_i],
%\end{equation}
%of graded $R$-algebras, where $\deg(u_i)=1$, such that the image has finite codimension.
An {\it isomorphism  of weakly marked algebras} is an isomorphism $A\simeq A'$
compatible with structures of weakly $N$-marked algebras for sufficiently large $N$.
\end{defi}

%Fix $N$ such that $N>\max(a_1,\ldots,a_n)$.

\begin{lem}\label{marked-weakly-marked-lem} A marked algebra $A$ of type $\ba$ over $R$,
is naturally a weakly $N$-marked algebra, where $N=\max(a_1,\ldots,a_n)$.
\end{lem}

\Pf . First, using the lattice condition in the definition of a marked algebra, we get that for $m>N$ one has
$$F_mA/F_{m-1}A\simeq \bigoplus_{i=1}^n F_{(m-1)\bD_1+\be_i}A/F_{m-1}A\simeq R^n.$$
This gives an identification of $\gr^F_{>N}A$ with $C_N$ as $R$-modules. The compatibility of the isomorphisms
$(\varphi_{\bn,i})$ with the product implies that this is an isomorphism of $R$-algebras. Next,
suppose $x\in F_\bn A$ is such that for each $m>m_0$ there exists $s\in S_m$ such that $xs\in F_{m_0}A$. 
If $\bn=\ba$ then $x\in F_\ba A=R$, so the condition $xs\in F_{m-1}A$ implies that $x=0$. 
Suppose now that $\bn>\ba$. Then there exists $i$ such that $\bn-\be_i\ge \ba$.
Let $\ov{x}= x \mod F_{\bn-\be_i}A$. Then considering for $s\in S_m$ with sufficiently large $m$ the condition
$$xs\in F_{m_0}A\sub F_{\bn+m\bD_1-\be_i}A,$$ 
we deduce that
$\ov{s}_i \cdot \ov{x}=0$,
hence, $\ov{x}=0$. Thus, we get $x\in F_{\bn-\be_i}$, and we can apply the induction to deduce that $x=0$.
\ed

Given a weakly marked algebra $A$ let us fix $N$ such that $A$ is weakly $N$-marked.
%the image of \eqref{marking-2-eq} contains $\bigoplus_{i=1}^nu_i^NR[u_i]$.
Let us consider the localization $K(A):=S^{-1}A$, where 
$$S=\sqcup_{m=0 \text{ or }m> N} S_m$$ with 
$S_0=\{1\}$ and $S_m$ for $m> N$ is as in Definition \ref{marked-alg-def}.
It is clear that $S$ is multiplicative.
The elements of $S$ are non-zero-divisors in $A$, so
the natural homomorphism $A\to K(A)$ is injective.
Note also that $S_m\neq\emptyset$ for all $m\gg 0$, so
the localization $K(A)$ does not depend on a choice of $N$.

Next, we define a decreasing $\Z^n$-valued algebra filtration 
 $(G_{\br}K(A))_{\br\in\Z^n}$ on $K(A)$ by
$$G_\br K(A):=\{ \frac{a}{s} \ |\ a\in F_{m\bD_1-\br}A, s\in S_m \text{ for some } m\gg 0\}$$
(note that for any $\br$ one has $m\bD_1-\br\in M_N$ for sufficiently large $m$).
We also set $G_r K(A):=G_{r\bD_1} K(A)$.

%Note that since $S_m\neq\emptyset$ for all $m> N$, any element $\frac{a}{s}$ with $a\in F_{m-r}A$, $s\in S_m$
%can be rewritten as a similar fraction with the denominator in $S_{m'}$ for any $m'\ge m+N$.

%Similarly, if $A$ is an  $\ba$-marked algebra then we define a $\Z^n$-valued filtration

%Let us extend the filtration $(F_\bm A)$ on a marked algebra $A$ to a $\Z^n$-valued filtration by setting
%$F_\bm A=0$ for $\bm\not\in (\Z_{\ge 0})^n$ and $F_\bm A=F_0=R$ for $

%Note also that 
%$$\OO_{A,\infty}:=G_0K(A)$$
%is an $R$-subalgebra of $K(A)$, and $G_rB:=G_rK(A)$ is the decreasing algebra filtration on $B$.

At this point we are going to make a formal change variables $t_i=u_i^{-1}$. Thus, we view each $R[u_i]$ as a subring
in $R[t_i,t_i^{-1}]$.

\begin{lem}\label{marked-infty-lem} Let $A$ be a weakly $N$-marked algebra.

\noindent
(i) Let $\br\in\Z^n$. For $\frac{a}{s}\in K(A)$, where $s\in S_m$, $m>N$ and $m\bD_1-\br\in M_N$, 
one has $\frac{a}{s}\in G_\br K(A)$ if and only if $a\in F_{m\bD_1-\br}A$.
%If in addition $A$ is an $\ba$-marked algebra then for $s\in S_m$ and $\br$, such that 
%$m(\be_1+\ldots+\be_n)-\br\in M^+_\ba$, one has
%$\frac{a}{s}\in G_\br K(A)$ if and only if $a\in F_{m(\be_1+\ldots+\be_n)-\br}A$.
%For any $M\ge 0$ one has 
%$$G_rK(A)=\{ \frac{a}{s} \ |\ a\in F_{m-r}A, s\in S_m \text{ for some } m\ge M\}.$$

\noindent
(ii) The filtration $(G_rK(A))$ on $K(A)$ is exhaustive and Hausdorff. 
%(i.e., $\cap_r G_rK(A)=0$). 
Also, for $\bm\in M_N$ we have
\begin{equation}\label{G-inters-eq}
G_{-\bm}K(A)\cap A= F_\bm A.
\end{equation}
If in addition $A$ is an $\ba$-marked algebra then \eqref{G-inters-eq}
holds for $\bm\in M_\ba$. 
 
\noindent
(iii) There a is natural isomorphism of graded $R$-algebras
\begin{equation}\label{gr-K-A-eq}
\gr_\bullet^G K(A)\simeq \bigoplus_{i=1}^n R[t_i,t_i^{-1}],
\end{equation}
where $\deg(t_i)=1$.
\end{lem}

\Pf . (i) It is enough to prove that if $s\in S_m$ and $sa\in F_{\bn+m\bD_1} A$ for $\bn\in M_N$ then $a\in F_{\bn}A$. 
We can use induction on $k$ such that $a\in F_k A$. Since the leading term $\ov{s}$ is not a zero divisor in $C_N$,
from the condition $sa\in F_{\min(k\bD_1,\bn)+m\bD_1}A$ we derive the existence of
$a'\in F_{\min(k\bD_1,\bn)}$ such that $a\equiv a'\mod F_{k-1} A$. Now we replace $a$ by $a-a'$
and use the induction assumption.
%The proof of the second assertion is similar.
%Hence, if $\frac{a}{s}\in G_rK(A)$, where $s\in S_m$, then $a\in F_{m-r}A$. 

\noindent
(ii) The fact that $F_\bullet A$ is exhaustive immediately implies that $\cup_r G_rK(A)=K(A)$. 
Next, \eqref{G-inters-eq} for $\bm\in M_N$ (resp., for $\bm\in M_\ba$ if $A$ is $\ba$-marked)
follows easily from (i). Now to see that $\cap_r G_r K(A)=0$, it is enough to prove that if $x\in A$
belongs to $G_r K(A)$ for all $r\ge 0$ then $x=0$. From (i) we see that such $a$ satisfies
$xs\in F_N A$ for any $s\in S_m$ with $m>N$. By the definition of a weakly $N$-marked algebra
this implies that $x=0$.

\noindent
(iii) For a given $r\in\Z$ let us pick $m>N$ such that $m-r>N$.
For any $s\in S_m$ and $a\in F_{m-r}A$ let us consider the elements 
$\ov{s}$ and $\ov{a}$ of $\bigoplus_{i=1}^n R[t_i,t_i^{-1}]$ defined by
$$\ov{s}:=s\mod F_{m-1}A\in \bigoplus_{i=1}^n Rt_i^{-m}, \ \
\ov{a}:=a\mod F_{m-r-1}A \in \bigoplus_{i=1}^n Rt_i^{-m+r}.$$
By the definition of $S_m$, $\ov{s}$ has invertible coefficients with each $t_i^{-m}$. Hence, 
$\ov{s}$ is invertible in $\bigoplus_{i=1}^n R[t_i,t_i^{-1}]$, so we can consider
the fraction $\frac{\ov{a}}{\ov{s}}\in \bigoplus_{i=1}^n Rt_i^r$.
It is easy to check that this construction gives
a well defined map of $R$-modules
$$G_rK(A)/G_{r+1}K(A)\to \bigoplus_{i=1}^n Rt_i^r,$$
which is in fact an isomorphism, and these maps are compatible with multiplication.
\ed

\begin{cor}\label{marked-alg-filtr-cor} 
(i) Let $A$ and $A'$ be weakly marked algebras over $R$, such that $A$ is weakly
$N$-marked and $A'$ is weakly $N'$-marked.
Then any isomorphism $\phi:A\to A'$ of weakly marked algebras satisfies
$\phi(F_\bm A)=F_\bm A'$ for $\bm\in M_N\cap M_{N'}$, and the induced isomorphism 
$$\gr^F_{\ge \max(N,N')}A\simeq \gr^F_{\ge \max(N,N')}A'$$ 
is compatible with the identifications \eqref{marking-2-eq}.

\noindent
(ii) Let $A$ and $A'$ be $\ba$-marked algebras over $R$. Then any isomoprhism $A\simeq A'$ of
weakly marked algebras is an isomorphism of $\ba$-marked algebras.
\end{cor}

\Pf . (i) Note that the localization $K(A)$ and the filtration $G_rK(A)$ depend only on $F_\bm A$ for $\bm\gg 0$ (see
Lemma \ref{marked-infty-lem}(i)). Now using Lemma \ref{marked-infty-lem}(ii) we recover $F_\bm A$ for all $\bm\in M_N$
as the intersection of $G_{-\bm}K(A)$ with $A$. This shows that $\phi(F_\bm A)=F_\bm A'$ for all 
$\bm\in M_N\cap M_{N'}$.
The isomorphism induced by $\phi$ on $\gr^F_{\ge \max(N,N')}$ acts as identity in
degrees $\gg 0$. This easily implies that it acts as identity in all degrees.

\noindent
(ii) An isomorphism $A\simeq A'$ of weakly marked algebras induces an isomorphism $K(A)\simeq K(A')$ compatible
with the filtrations $G_\br$. Since the filtrations $F_\bm$ are recovered from $G_\br$ via \eqref{G-inters-eq},
the assertion follows.
\ed

For a weakly marked algebra $A$ let us consider the completion
$$\widehat{K(A)}:=\varprojlim_r K(A)/G_rK(A)$$
% \ \text{ and } \hat{B}=\varprojlim_r B/G_rB\sub \hat{K(A)}.$$
Then Lemma \ref{marked-infty-lem}(iii) easily implies that there exists a noncanonical isomorphism 
%$\hat{B}\simeq \bigoplus_{i=1}^n R[[t_i]]$,
$$\widehat{K(A)}\simeq \bigoplus_{i=1}^n R((t_i)),$$ 
compatible with the $\Z^n$-filtrations $\widehat{G_\br K(A)}$ and $(\bigoplus_{i=1}^n t_i^{r_i}R[[t_i]])$,
and with the isomorphism \eqref{gr-K-A-eq}.

Let  $\MA^\ba(R)$ (resp., $\WMA(R)$) denote the groupoid of marked algebras of type $\ba$ (resp.,
weakly marked algebras) over $R$ and their isomorphisms. Note that by Corollary
\ref{marked-alg-filtr-cor}, $\MA^\ba(R)$ is a full subgroupoid in
$\WMA(R)$.

\begin{prop}\label{marked-alg-prop}
(i) Let $\SA(\HH_R)$ be the set of $R$-subalgebras $W\sub \HH_R$, such that for some $N\ge 0$,
one has $H^1(W(N\bp_1+\ldots+N\bp_n))=0$, i.e., $\HH_R=W+\HH_{R,\ge -N}$.
Let $[\SA(\HH_R)/\fG(R)]$ be the groupoid with objects $\SA(\HH_R)$ and morphisms given
by elements of $\fG(R)$. Then the natural functor 
$$[\SA(\HH_R)/\fG(R)]\to \WMA(R)$$
 is fully faithful.

\noindent
(ii) The natural functor of groupoids
$$\Phi:[ASG^\ba(R)/\fG(R)]\to \MA^\ba(R),$$
%$$[ASG^\ba(R)/\fG(R)]\to \MA^\ba(R)$$ 
induced by the functor in part (i), is an equivalence.
\end{prop}

\Pf . (i) Given a subalgebra $W\sub \HH_R=\bigoplus_{i=1}^n R((t_i))$ in $\SA(\HH_R)$, 
we define the structure of a weakly marked algebra on $W$ by considering the filtration
\begin{equation}\label{F-bm-W-eq}
F_\bm W:=W\cap \bigoplus_{i=1}^n t_i^{-m_i}R[[t_i]],
\end{equation}
where each $m_i\ge N$, with $N$ large enough (note that for a weakly marked algebra structure we only need
$F_\bm W$ for $\bm\in M_N$). By assumption, $F_mW/F_{m-1}W\simeq R^n$,
and this gives an isomorphism of $\gr^F_{>N}W$ with $C_N$. 
To see that this defines a weakly marked algebra structure on $W$ we need to check that if for some $m_0>N$ and some
$x\in W$, for each $m>m_0$ there exists $s\in S_m$ such that $xs\in F_{m_0}W$, then $x=0$.
But this follows immediately by considering the leading terms of the components of $x$ in $R((t_i))$.

%Furthermore, any element $f\in F_mW$, such that
%$f\mod F_{m-1}W$ has invertible components, is invertible in $\HH_R=\bigoplus_{i=1}^n R((t_i))$, so we have a structure of
%a weakly marked algebra on $W$. 

We can view an action of $g\in\fG(R)$ as an automorphism of $\HH_R$, compatible with the filtration
by the submodules $F_{\bm}=\bigoplus_{i=1}^n t_i^{-m_i}R[[t_i]]$ and
acting as identity on the associated graded algebra of the filtration $(F_m)$. 
Thus, if $gW=W'$ then we have an induced isomorphism of
weakly marked algebras $W\to W'$.

Note that the localization $K(W)$ can be viewed as an $R$-subalgebra in $\HH_R$,
since all the elements of $S$ are invertible in the latter ring. Furthermore, it is easy to see
that 
$$G_\br K(W)=K(W)\cap \bigoplus_{i=1}^n t_i^{r_i}R[[t_i]].$$
Now Lemma \ref{marked-infty-lem} implies that the embedding of $K(W)$ into $\bigoplus_{i=1}^n R((t_i))$
induces an isomorphism 
\begin{equation}\label{K-W-compl-eq}
\widehat{K(W)}\rTo{\sim}\bigoplus_{i=1}^n R((t_i))=\HH_R,
\end{equation}
compatible with the embedding of $W$ into $K(W)$ and into $\HH_R$.  

Given a subalgebra $W\sub \HH_R$ and an element $g\in \fG(R)$, let us consider the corresponding isomorphism 
$\phi_g:W\to W'=gW$ of weakly marked algebras. Passing to completions we get an isomorphism
$$\hat{\phi}_g:\widehat{K(W)}\to \widehat{K(W')},$$
fitting into a commutative diagram
\begin{diagram}
\widehat{K(W)}&\rTo{\sim}& \HH_R\\
\dTo{\hat{\phi}_g}&&\dTo{g}\\
\widehat{K(W')}&\rTo{\sim}& \HH_R\\
\end{diagram}
This shows that $g$ can be recovered from $\phi_g$.

Conversely, given an isomorphism of weakly marked algebras $\phi:W\to W'$, for subalgebras $W,W'\sub \HH_R$, 
we get an induced isomorphism $\hat{\phi}:\hat{K(W)}\to \hat{K(W')}$. It corresponds via the isomorphisms
\eqref{K-W-compl-eq} for $W$ and for $W'$ to an automorphism of $\HH_R$,
compatible with the $\Z^n$-filtration $(\bigoplus_{i=1}^n t_i^{r_i}R[[t_i]])$, and inducing the identity on the associated
graded with respect to the $\Z$-filtration $(\bigoplus_{i=1}^n t_i^mR[[t_i]])$. Hence, this automorphism of $\HH_R$
is given by the action of an element of $g\in \fG$. Since $\hat{\phi}|_W=\phi$, we see that $g(W)=W'$.

\noindent 
(ii) 
%Since the action of $\fG$ on $ASG^\ba$ admits a section (see Corollary \ref{ASG-sec-cor}), we can identify
%$R$-points of $ASG^\ba/\fG$ with the orbits of the free $\fG(R)$-action on $ASG^\ba(R)$.
%For each $R$, we have a natural functor of groupoids
%compatible with the functor considered in (i). 
%Namely, this 
The functor $\Phi$ associates with $W\sub \HH_R$,
the same $W$ viewed as an algebra, with the filtration $(F_\bm W)_{\bm\in M_\ba}$
given by \eqref{F-bm-W-eq}
and with the isomorphisms $\varphi_{\bn,i}$ coming from the standard basis of $\HH_R$.
By part (i) and Corollary \ref{marked-alg-filtr-cor}(ii), the functor $\Phi$ is fully faithful. It remains to check that it is
essentially surjective. Let $A$ be a marked algebra of type $\ba$ over $R$. Let us fix some isomorphism 
$\widehat{K(A)}\simeq \HH_R$, compatible with the filtrations $G_\br$ and with \eqref{gr-K-A-eq}.
Then we obtain an embedding 
$$A\to K(A)\to \widehat{K(A)}\simeq \HH_R.$$
It is easy to see that $A$, viewed as a subagebra of $\HH_R$, gives an $R$-point of $ASG^\ba$, that maps to $A$ under 
$\Phi$.
\ed

\begin{rem} Note that for any $\bS=\sqcup_{i=1}^n S_i$, where $\Z_{\ge 0}\sub S_i\sub\Z$ with 
$|S_i\setminus \Z_{\ge 0}|<\infty$, we have an inclusion
$$U_{\bS,1}\cap ASG(R)\sub \SA(\HH_R),$$
where $U_{\bS,1}$ is the open cell \eqref{U-S-1-eq}.
\end{rem}

\subsection{Gr\"obner bases}\label{Groebner-sec}

Recall that Gr\"obner bases can be defined with respect to any complete order $<$ on the monomials, which is
{\it admissible}, i.e., satisfies
\begin{equation}\label{groebner-order}
a<b \implies ac<bc, \ \ a<ab.
\end{equation}
If $<$ is such a complete order and $\deg$ is a nonnegative grading on the variables 
then we can define a new admissible complete order $<_{\deg}$ by
$$a<_{\deg} b \ \text{ if either } \deg(a)<\deg(b) \ \text{ or } \deg(a)=\deg(b) \text{ and } a<b.$$
%For example, the degree reverse lexicographical order is obtained like this with $<$ being the lexicographical order.
If we have two nonnegative gradings $\deg_1$, $\deg_2$ then iterating this construction we
obtain the admissible complete order $<_{\deg_1,\deg_2}$, where we first check whether $\deg_1(a)<\deg_1(a)$,
then whether $\deg_2(a)<\deg_2(a)$ and finally in the case of two ties, whether $a<b$.

Let $R$ be a commutative ring, and let 
$W\sub \HH_R=\bigoplus_{i=1}^n R((t_i))$ be a subalgebra, which as a subspace comes
from an $R$-point of one of the open cells of $SG$.
Let us set 
$$\bD_1:=\bp_1+\ldots+\bp_n.$$

We make the following additional assumptions on $W$:

\noindent
($\star$) Assume that for some $N>0$ one has $H^1(W(N\bD_1))=0$ and $H^0(W(N\bD_1))/R\cdot 1$ is a free $R$-module of finite rank.

%We start with a point $W\in ASG^\ba(R)$, which we view as  
%$(C,p_1,\ldots,p_n)$ be a pointed curve such that $h^1(\OO_C(ND))=0$ for some $N>0$, where
%$D=p_1+\ldots+p_n$. 

Then for each $i=1,\ldots,n$ and each $j=0,\ldots,N$, we can choose an element
$$h_i(j)\in H^0(W((1+j)\bp_i+N\bD_1)) \text{ such that } h_i(j)\equiv t_i^{-N-1-j}\mod H^0(W((j\bp_i+N\bD_1))).$$
In other words, $h_i(j)$ is required to have the pole of order exactly $N+1+j$ in $t_i$ and a pole of order at most $N$
in all other $t_{i'}$'s. Let us set $f_i=h_i(0)$. 
Also, let $1,g_1,\ldots,g_r$ be an $R$-basis of $H^0(W(N\bD_1))$ which exists by the assumption ($\star$).

\begin{prop}\label{Groebner-prop} 
The elements 
$$(f_i,h_i(j),g_k)_{i=1,\ldots,n, j=1,\ldots,N,k=1,\ldots,r}$$ 
generate the algebra $W$.
Let us view $W$ as the quotient of the polynomial algebra in $f_i,h_i(j),g_k$, by an ideal $J$.
Consider the ordering $<_{\deg_1,\deg_2}$ on monomials in $f_i,h_i(j),g_k$,
%viewed as independent var
where 
$$\deg_1(h_i(j))=N+1+j, \ \deg_2(h_i(j))=0, \ \deg_1(g_k)=N, \ \deg_2(g_k)=1$$
(where $j\ge 0$ and $f_i=h_i(0)$),
and $<$ is the reverse lexicographical order for any ordering of the variables $<$ 
such that $\deg_1(a)<\deg_1(b)$ implies $a<b$.
Then the corresponding normal monomials are 
\begin{equation}\label{normal-monomials-eq}
f_i^m, \ f_i^mh_i(j), \ g_k,
\end{equation}
with $m\ge 0$, $i=1,\ldots,n$, $j=1,\ldots,N$. The 
corresponding Gr\"obner basis of $J$ consists of elements of the form
(written starting from the leading term)
\begin{equation}\label{Groebner-leading-eq}
\begin{array}{l}
h_i(j)h_i(j')-f_ih_i(j+j')+ \text{terms}(\deg_1\le 2N+1+j+j'), \ \text{ for } j+j'\le N,\\
h_i(j)h_i(j')-f_i^2h_i(j+j'-N-1)+\text{terms}(\deg_1\le 2N+1+j+j'), \ \text{ for } j+j'\ge N+1,\\
h_i(j)h_{i'}(j')+P_i(\le 2N+1+j)+P_{i'}(\le 2N+1+j')+\text{terms}(\deg_1\le 2N), \\ 
h_i(j)f_{i'}+Q_i(\le 2N+1+j)+Q_{i'}(\le 2N+1)+\text{terms}(\deg_1\le 2N), \\ 
f_if_{i'}+R_i(\le 2N+1)+R_{i'}(\le 2N+1)+ \text{terms}(\deg_1\le 2N), \\ 
h_i(j)g_k+\ldots, \ f_ig_k+\ldots, \ g_kg_l+\ldots,   
\end{array}
\end{equation}
where $i\neq i'$, $j\ge 1$, $j'\ge 1$, the expression $\text{terms}(\deg_1<a)$ stands for a linear combination
of normal monomials of $\deg_1<a$, and $P_i(\le a), Q_i(\le a), R_i(\le a)$ are some linear combinations
of normal monomials of the form $f_ih_i(j)$, $h_i(j)$, and $f_i$ with $\deg_1\le a$.
\end{prop}

\Pf . First, we claim that the monomials \eqref{normal-monomials-eq} form a basis $W$.
Indeed, for each $m\ge N$ let $\BB_m$ be the subset of these monomials that have $\deg_1$ equal to $m$.
Since $\BB_0\cup\BB_N=\{1,g_1,\ldots,g_r\}$ is a basis of $H^0(W(N\bD_1))$, it is enough to prove that
for each $m\ge N+1$ the elements of $\BB_m$ project to a basis of $H^0(W(m\bD_1))/H^0(W((m-1)\bD_1))$.
Note that $\BB_m$ consists of elements $f_i^ph_i(j)$, where $i=1,\ldots,g$, such that
$m=p(N+1)+N+1+j$ with $p\ge 0$ and $0\le j\le N$ (where we use $h_i(0)=f_i$ for $j=0$).
Thus, for each $i=1,\ldots,g$ there exists a unique monomial in $\BB_m$ having the pole of order exactly $m$ in $t_i$ and 
poles of order $<m$ in other $t_{i'}$'s. This immediately implies our claim.

It is easy to see that \eqref{normal-monomials-eq} are precisely all the monomials that are not divisible by any of the initial monomials in \eqref{Groebner-leading-eq}. 
Thus, it remains to prove that the ideal $J$ contains elements of the form \eqref{Groebner-leading-eq}, where 
the leading terms are bigger than the subsequent terms. 
In other words, we have to check that any initial monomial $M$
in \eqref{Groebner-leading-eq}, viewed as an element of $W$, is a linear combination of  smaller monomials
among \eqref{normal-monomials-eq} of prescribed form.

For the monomials of the form $h_i(j)h_{i'}(j')$ (including the possibilities $i=i'$ and $j=0$ or $j'=0$) this
is obtained by analyzing polar parts of the expansion of this monomial in $t_k$. 
Note that for $k\neq i,i'$ the order of pole in $t_k$ is $\le 2N$, while for $i\neq i'$ the order of pole in $t_i$ is $\le 2N+1+j$.
In the case $i=i'$ we also use the fact that the expansion of $h_i(j)h_i(j')$ in $t_i$ starts with $t_i^{-2N-2-j-j'}$ and
find the matching normal monomial.

In the remaining case when $M$ is one of the monomials $h_i(j)g_k$, $f_ig_k$ or 
$g_kg_l$, we have $\deg_1(M)>N$, so all elements of
$\BB_m$ have $\deg_2=0$, hence are smaller than $M$. 
\ed

%\subsection{Affine moduli spaces}

\subsection{Proof of Theorems A and B}\label{proof-A-B-sec}
First, we will work over $\Q$. Then we will explain what modifications have to be made to work over $\Z[1/N]$.

Recall that by Corollary \ref{curves-sec-cor}, we have a morphism
$$\wt{\UU}^{ns}_{g,n}(\ba)\simeq \Kr^{-1}(ASG\cap\Sigma^{\ba,i_0})
\rTo{\ov{\Kr}} ASG\cap\Sigma^{\ba,i_0},$$
where $\Sigma^{\ba,i_0}$ is a section for the $\fG$-action on $SG^\ba$ defined in Proposition \ref{section-open-prop}.

Note that $ASG\cap\Sigma^{\ba,i_0}$ is an affine scheme, as a closed subscheme of the affine scheme $\Sigma^{\ba,i_0}$,
although at the moment we do not know that it is of finite type. Let $R$ be the corresponding commutative ring, so
$ASG\cap\Sigma^{\ba,i_0}=\Spec(R)$. 
Let $W\sub\HH_R$ be the subalgebra corresponding to the universal $R$-point of $ASG\cap\Sigma^{\ba,i_0}$.
We are going to use the procedure of Section \ref{Groebner-sec} to 
construct a canonical set of generators of the $R$-algebra $W$ and an $R$-basis of $W$ consisting of normal monomials.
%$H^0(C\setminus D,\OO)$,  where $D=p_1+\ldots+p_n$, as in Proposition \ref{Groebner-prop}.

Namely, let $N$ be the maximal of the numbers $a_1,\ldots,a_n$. 
Then the condition $(\star)$ is satisfied for $W$ and the elements 
\begin{equation}\label{g-k-f-i-p-eq}
f_i[-p], \text{ for } i=1,\ldots,n, a_i<p\le N,
\end{equation}
where $f_i[-p]$ are given by \eqref{f-i-p-eq}, \eqref{coordinates-a-ij-eq} (normalized by $\a_{ii}[-p,0]=0$),
form a basis of $H^0(W(N\bD_1))/R\cdot 1$.
Similarly, for $i=1,\ldots,n$ and $j=0,\ldots,N$ we set
$h_i(j)=f_i[-N-j-1]$ (and $f_i=h_i(0)=f_i[-N-1]$). Thus, we are in the setting of Proposition \ref{Groebner-prop},
with the elements \eqref{g-k-f-i-p-eq} playing the role of $(g_k)$.

Therefore, the algebra $W$ has the basis \eqref{normal-monomials-eq} as an $R$-module,
and the Gr\"obner basis of the ideal of relations between $(f_i,h_i(j),g_k)$ has form
\eqref{Groebner-leading-eq}. Note that by Buchberger's Criterion (see \cite[Thm.\ 15.8]{Eisenbud}), 
if we write general relations of the
form \eqref{Groebner-leading-eq}
and treat coefficients as indeterminate variables, then the condition that the normal monomials
form a basis will give a system of polynomial equations on the coefficients. 
Thus, we have an affine scheme $\wt{S}_{GB}$ (where ``GB" stands for Gr\"obner bases) of
finite type over $\Q$, which parametrizes Gr\"obner bases of this form.
%Let $S_{GB}=\Spec(R_{GB})$ 
We are going to define a  closed subscheme $S_{GB}\sub \wt{S}_{GB}$ by imposing additional equations
on the coefficients in the relations \eqref{Groebner-leading-eq}
(in the case when all $a_i=1$ we considered such a scheme in \cite[Lem.\ 1.2.2.]{P-ainf}). 
Namely, we define $S_{GB}$ inside $\wt{S}_{GB}$ by imposing the additional
requirement for the relations with the leading terms $f_i[-p]f_i$ and $f_i[-p]f_{i'}$ to have form
%For a monomial $M$ in the left-hand side of
%\eqref{Groebner-leading-eq} and for a normal monomial $M_0$, let us denote by $c(M,M_0)$ the coefficient of $M_0$
%in the expression of $M$ in terms of normal monomials.
%Then we define $S_{GB}$ inside $\wt{S}_{GB}$ by the equations
\begin{equation}\label{S-GB-eq}
\begin{array}{l}
f_i[-p]f_i=h_i(p)+A_i(\le N+p)+\sum_{k\neq i} A_k(\le 2N)+\text{terms}(\deg_1\le N),\\
%\text{terms}(\deg_1<N+1-p)???,\\
%c(f_i[p]f_i,h_i(-p))=1, \\
f_i[-p]f_{i'}=B_{i'}(\le N+1+a_{i'})+\sum_{k\neq i'} B_k(\le 2N)+\text{terms}(\deg_1\le N),
%c(f_i[p]f_{i'},h_{i'}(j))=0 \ \text{ for } j>a_{i'},
\end{array}
\end{equation}
where $i\neq i'$, $a_i<p\le N$, $A_i(\le a)$ and $B_i(\le a)$ are some linear combinations of the elements
$h_i(l)$ and $f_i$ of $\deg_1\le a$.
%where other terms in the right-hand side have smaller $\deg_1$

The above construction defines a morphism of affine schemes
$ASG\cap\Sigma^{\ba,i_0}\to \wt{S}_{GB}$, and we claim that it factors through a morhism
$$i:ASG\cap\Sigma^{\ba,i_0}\to S_{GB}.$$
%(the fact that \eqref{h-h-leading-equations} holds follows by considering the expansions in $t_i$).
Indeed, the equations \eqref{S-GB-eq} follow from the fact that the expansions of $f_i[-p]f_i$ (resp., $f_i[-p]f_{i'}$)
in $t_k$ with $k\neq i$ (resp., $k\neq i'$) have poles of order $\le 2N$, while the expansion of $f_i[-p]f_i$ in $t_i$
starts with $t_i^{-N-1-p}$ (resp., the expansion of $f_i[-p]f_{i'}$ in $t_{i'}$ has pole of order $\le N+1+a_{i'}$).

Next, we are going to construct a morphism 
\begin{equation}\label{r-map}
r:S_{GB}\to \wt{\UU}^{ns}_{g,n}(\ba)
\end{equation}
such that the compositions
\begin{equation}\label{Kr-i-r}
\wt{\UU}^{ns}_{g,n}(\ba)\rTo{\ov{\Kr}} ASG\cap\Sigma^{\ba,i_0}\rTo{i} S_{GB}\rTo{r}  \wt{\UU}^{ns}_{g,n}(\ba) \ \text{ and}
\end{equation}
\begin{equation}\label{i-r-Kr}
ASG\cap\Sigma^{\ba,i_0}\rTo{i} S_{GB}\rTo{r}  \wt{\UU}^{ns}_{g,n}(\ba)\rTo{\ov{\Kr}} ASG\cap\Sigma^{\ba,i_0}
\end{equation}
are the identity morphisms. Since $S_{GB}$ is an affine scheme of finite type, using Lemma \ref{closed-subscheme-lem} below, this would imply that both maps
$$i\circ \ov{\Kr}:\wt{\UU}^{ns}_{g,n}(\ba)\to S_{GB} \ \text{ and }\ i:ASG\cap\Sigma^{\ba,i_0}\to S_{GB}$$ 
are 
closed embeddings that factor through each other. It would follow that
$$\wt{\UU}^{ns}_{g,n}\simeq ASG\cap\Sigma^{\ba,i_0}$$
and that this is an affine scheme of finite type.

To construct the morphism \eqref{r-map}
we use the universal family over $S_{GB}=\Spec(R_{GB})$, i.e., the algebra $A_{GB}$ over $R_{GB}$, obtained
as a quotient by the ideal generated by the universal Gr\"obner basis. Note that 
the basis of normal monomials \eqref{normal-monomials-eq} in our case gives the following $R_{GB}$-basis
of $A_{GB}$:
$$f_i^m, \ f_i^mh_i(j), \ f_i[-p],
$$
where $m\ge 0$, $i=1,\ldots,n$, $j=1,\ldots,N$, $a_i<p\le N$.
Let us define the increasing filtration $(F_m A_{GB})$ on $A_{GB}$, 
by letting $F_m A_{GB}$ to be the $R_{GB}$-submodule spanned by the normal monomials with $\deg_1\le m$.
Any leading term among \eqref{Groebner-leading-eq}
is expressed in terms of elements with smaller or equal $\deg_1$, hence $F_\bullet A_{GB}$ is an algebra filtration.
We also have a $M_N$-valued filtration defined by
\begin{equation}\label{weak-marked-A-GB-eq}
F_{k\bD_1+\sum_i m_i\be_i} A_{GB}=F_k+\sum_{k<(p+1)(N+1)+j\le k+m_i}R\cdot f_i^ph_i(j),
\end{equation}
where $0\le m_i\le k/N$. We will see later that it defines a structure of a weakly $N$-marked algebra on $A_{GB}$.

Let $\RR(A_{GB})=\bigoplus_{m\ge 0} F_m A_{GB}$ be the corresponding Rees algebra, and let $T\in \RR(A_{GB})$
be the element of degree $1$ corresponding to $1\in F_1 A_{GB}$. Note that $T$ is a non-zero-divisor and
$\RR(A_{GB})/(T)\simeq \gr^F A_{GB}$, the associated graded algebra for the filtration $F_\bullet A_{GB}$.

We claim that 
$$C_{GB}:=\Proj(\RR(A_{GB})$$ 
is a flat family of curves over $R_{GB}$, equipped with the natural marked points $p_1,\ldots,p_n$,
making it into an $R_{GB}$-point of $\wt{\UU}^{ns}_{g,n}(\ba)$. 

Note that for $m>N$ the quotient $F_m A_{GB}/F_{m-1} A_{GB}$ has the basis $(f_i^ph_i(j))_{i=1,\ldots,n}$, where
$p\ge 0$ and $j$, $0\le j\le N$, are unique such that $m=p(N+1)+N+1+j$. 
Note also that in $\gr^F A_{GB}$ we have the relations 
$$h_i(j)h_i(j')=f_ih_i(j+j'), \ \text{ for } j+j'\le N,$$
$$h_i(j)h_i(j')=f_i^2h_i(j+j'-N-1), \ \text{ for } j+j'\ge N+1,$$
$$h_i(j)h_{i'}(j')=0, \ h_i(j)f_{i'}=0, \ f_if_{i'}=0, \ \text{ for } i\neq i'.$$ 
This easily implies the following identification of the truncated associated graded algebra:
\begin{equation}\label{gr-A-GB-eq}
\gr^F_{>N} A_{GB}\simeq \bigoplus_{i=1}^n R_{GB}[u_i]_{>N},
\end{equation}
where the variables $u_i$ have degree $1$. 
Therefore, the divisor $D:=(T=0)$ in $C_{GB}$ is isomorphic
to $\Proj(\bigoplus_{i=1}^n R_{GB}[u_i])$, i.e., it is the disjoint union of $n$ copies of the base $\Spec(R_{GB})$.
Note that the algebra of functions on the complement to $D$ is isomorphic to
the degree $0$ part of the localization $\RR(A_{GB})[T^{-1}]$,
which is $A_{GB}$, so we get an identification of $R_{GB}$-algebras
\begin{equation}\label{H0-A-GB-eq}
H^0(C_{GB}\setminus D,\OO)\simeq A_{GB}.
\end{equation}

Let $F_i, H_i(j)$ be elements $f_i, h_i(j)$ viewed as elements of $F_m A_{GB}=\RR_m(A_{GB})$, 
where $m$ is equal to $\deg_1$ of the corresponding element of $A_{GB}$. We define the marked points
$p_1,\ldots,p_n$, so that $D=p_1\sqcup\ldots\sqcup p_n$, and $F_i\neq 0$, $H_i(j)\neq 0$ at $p_i$.

As in the proof of  \cite[Thm.\ 1.2.4]{P-ainf}, we check that $C_{GB}$ is flat over $R_{GB}$,
the sheaf $\OO(1)$ on $C_{GB}$ is locally free, and the divisor $D=(T=0)$ in $C_{GB}$
is ample. 

%Define valuation associated with $p_1$ as follows. 
As in the proof of \cite[Thm.\ 1.2.4]{P-ainf}, 
we consider the affine open neighborhood of $p_i$ in which $F_i$ and all $H_i(j)$ are invertible.
Then $p_i$ is the intersection of $D$ with this open neighborhood, so
\begin{equation}\label{t-i-T-F-H-1}
t_i:=TF_i/H_i(1)
\end{equation}
generates the ideal of $p_i$ over this neighborhood.
Since $t_i$ is a non-zero-divisor, we get that $C_{GB}$ is smooth over $R_{GB}$ near $p_i$. 

Let us observe also that \eqref{gr-A-GB-eq} implies that for any $j\ge 0$
the function $\frac{H_i(j)F_i^{N+1+j}}{H_i(1)^{N+1+j}}\in A(p_i)$, 
has value $1$ at $p_1$.

We claim that considering polar conditions near $p_1,\ldots,p_n$ we recover the filtration $F_\bm A_{GB}$ on
$A_{GB}=H^0(C_{GB}\setminus D,\OO)$, for $\bm\in M_N$ (see \eqref{weak-marked-A-GB-eq}).
For this we need to estimate the orders of poles of $h_i(j)$, $j\ge 0$, and of $f_i[p]$
at all the marked points. By symmetry, it is enough to consider the marked point $p_1$.
We use $t_1=TF_1/H_1(1)$ as a local parameter at $p_1$.
We have for $j\ge 0$,
$$h_1(j)=H_1(j)/T^{N+1+j}=\frac{H_1(j)F_1^{N+1+j}}{H_1(1)^{N+1+j}}\cdot t_1^{-N-1-j}.
$$
As we observed before, the fraction is invertible at $p_1$, so we get that the order of pole is exactly $N+1+j$.
%This argument also works for $f_1=h_1(0)$.

If $G$ is a monomial of $\deg_1=N$, viewed as an element of $F_N A_{GB}=\RR_N(A_{GB})$, then 
$GF_1^N/H_1(1)^N\in A(p_1)$ is regular at $p_1$. This implies that 
$$\frac{G}{T^N}=\frac{GF_1^N}{H_1(1)^N}\cdot t_1^{-N}$$
has pole of order $\le N$ at $p_1$. Thus, any element $f_i[-p]\in A_{GB}$ has pole of order $\le N$ at $p_1$.

Note that for $i\neq 1$ the function $f_i/f_1=F_i/F_1$ vanishes at $p_1$, 
hence, $f_i$ has pole of order $\le N$ at $p_1$. 
Similarly, $h_i(j)/h_1(j)$ vanishes at $p_1$, so $h_i(j)$ has pole of order $\le N+j\le 2N$ at $p_1$.
Thus, we deduce that all normal monomials of $\deg_1\le 2N$, have poles of order $\le 2N$ at $p_1$.
Indeed, these are just the elements $f_i[-p]$ and $h_i(j)$ for $j\le N-1$ (with possibly $i=1$).

Now let us check by induction on $j$ that for $i\neq 1$,
$h_i(j)$ has pole of order $\le N$ at $p_1$. 
By \eqref{Groebner-leading-eq}, we have 
$$h_i(j)f_1=Q_i(\le 2N+1+j)+Q_1(\le 2N+1)+\text{terms}(\deg_1\le 2N),$$
where $Q_i(\le 2N+1+j)$ is a linear combination of $f_i$, $h_i(j')$ and $f_ih_i(k)$ for $0\le k<j$, while
$Q_1(\le 2N+1)$ is a linear combination of $h_1(j')$ and $f_1$. Using the induction assumption we see
that $Q_i(\le 2N+1+j)$ has pole of order $\le 2N$ at $p_1$. Hence, we derive that
$h_i(j)f_1$ has pole of order $\le 2N+1$ at $p_1$, which implies that $h_i(j)$ has pole of order $\le N$ at $p_1$.

The above information on the poles implies that the elements $(f_i^ph_i(j))$ with $k<(p+1)(N+1)+j\le k+m_i$
project to a set of generators of $H^0(C_{GB},\OO(kD+\sum_i m_ip_i))$, provided $0\le m_i\le k/N$ for each $i$.
Thus, we deduce that for each $\bm=k\bD_1+\sum_i m_i\be_i\in M_N$ one has 
\begin{equation}\label{weak-marked-A-C-GB-eq}
F_\bm A_{GB}=H^0(C_{GB},\OO(kD+\sum_i m_ip_i)),
\end{equation}
where the left-hand side is given by \eqref{weak-marked-A-GB-eq}.

We need a more precise information on the poles of $f_i[-p]$.
Let us consider the first of the relations \eqref{S-GB-eq} for $i=1$. As we have shown, the 
terms of $\deg_1\le N$ have poles of order $\le N$ at $p_1$. The same is true for linear combinations of 
$h_k(j)$ with $k\neq 1$, so we derive that
$f_1[-p]f_1-h_1(p)$ has pole of order $\le N+p$ at $p_1$. Thus, $f_1[-p]f_1$ has pole of order exactly $N+1+p$ at $p_1$,
which implies that $f_1[-p]$ has pole of order exactly $p$ at $p_1$.
Similarly, considering the second of the relations \eqref{S-GB-eq} for $i\neq 1$, $i'=1$, we see that 
$f_i[-p]f_1$ has pole of order $\le N+1+a_1$ at $p_1$, which implies that $f_i[-p]$ has pole of order $\le a_1$ at $p_1$.

The above information on the poles implies that the elements $(f_i[-p])$ project to an $R_{GB}$-basis of
$H^0(C_{GB},\OO(ND))/H^0(C_{GB},\OO(a_1p_1+\ldots+a_np_n))$. Hence, we also
deduce that $H^0(C_{GB},\OO(a_1p_1+\ldots+a_np_n)=R_{GB}$.

Also, for $m\ge N$ the isomorphism $F_mA_{GB}\simeq H^0(C_{GB},\OO(mD))$ shows that $H^0(C_{GB},\OO(mD))$
 is a free $R_{GB}$-module of rank $mn-g+1$. Hence, we
 have a family of curves of arithmetic genus $g$ with marked points satisfying the condition 
 $h^0(\OO(a_1p_1+\ldots+a_np_n))=1$.

We define the family of tangent vectors at the marked points using the local parameters $t_1,\ldots,t_n$
given by \eqref{t-i-T-F-H-1}.
Hence, we obtained an $R_{GB}$-point of  $\wt{\UU}^{ns}_{g,n}(\ba)$,
i.e., defined a morphism $r:S_{GB}\to \wt{\UU}^{ns}_{g,n}(\ba)$.

%Claim: for $i\neq 1$, $j>0$ $H_i(j)\in \II^2$, $F_i\equiv c_i\cdot TF_1/H_1(1)\mod \II^2$.
%Define a homomorphism $A(p_i)\to k[[t_i]]$ identifying the completion of the local rings at $p_i$ ith $k[[t_i]]$??? 
%Check flatness, smoothness of points, etc.???

Since $ASG\cap \Sigma^{\ba,i_0}\simeq ASG^\ba/\fG$,
to check that $\ov{\Kr}\circ r\circ i=\id$, by Proposition \ref{marked-alg-prop}(ii), we need to show
that the isomorphism \eqref{H0-A-GB-eq} is compatible with the marked algebra structures on both sides.
By Proposition \ref{marked-alg-prop}(i), it is enough to show the compatibility with the weakly marked algebra
structures on both sides. But this follows from the identification \eqref{weak-marked-A-C-GB-eq}.

Next, we need to check that $r\circ i\circ\ov{\Kr}=\id$. Thus, starting with a family of curves $(C,p_1,\ldots,p_n)$ over
a commutative ring $R$, with fixed tangent vectors at $p_i$, defining an $R$-point of
$\wt{\UU}^{ns}_{g,n}(\ba)$, we apply the morphism $r$ to the algebra $A=H^0(C\setminus D,\OO)$, viewed
as a subalgebra of $\bigoplus_i R((t_i))$, where $t_i$ are the canonical parameters at the marked points 
(see Corollary \ref{curves-sec-cor}).
Note that the corresponding Rees algebra $\RR(A)$ (truncated in degree $\ge N$) can be identified with 
$\bigoplus_{m\ge N} H^0(C,\OO(mD))$.  
Since $D$ is ample, we have a natural isomoprhism 
$$C\simeq \Proj(\RR(A)),$$
compatible with the marked points $p_i$. Since the expansion of 
$f_i/h_i(1)=f_i[-N-1]/f_i[-N-2]$ at $p_i$ has form $f_i/h_i(1)=t_i+\ldots$, we get
the compatibility with the tangent vectors at the marked points.

This finishes the proof of Theorem A(i) and of Theorem B, working over $\Q$. Let us now explain how to pass
to $\Z[1/N]$.
Recall that we needed to work over $\Q$ so that we could use the section $\Sigma^{\ba,i_0}$ for the $\fG$-action
on $SG^\ba$ and the corresponding canonical formal parameters
at the marked points over $\wt{\UU}^{ns}_{g,n}(\ba)$ (see Proposition \ref{section-open-prop} and Corollary
\ref{curves-sec-cor}). However, if we only keep track of finite jets of formal parameters and work
with the Grassmannian of subspaces of $\HH_{\ge -N}/\HH_{\ge N}$ for large $N$, then it would be enough to
work over $\Z[1/N']$. Now we observe that since the affine scheme $S_{GB}$ is of finite type, we can define
versions of the morphisms in \eqref{Kr-i-r}, replacing $ASG$ 
with the following truncated version for large enough $N$: 
we consider subspaces $\ov{W}\sub \HH_{\ge -N}/\HH_{\ge 2N}$ such that
$1\in \ov{W}$ and 
$$\ov{W}\cdot\ov{W}\sub (\ov{W}+\HH_{\ge N})/\HH_{\ge N}\sub \HH_{\ge -N}/\HH_{\ge N}.$$
Then the entire argument above still works and we only need 
$\Spec(\Z[1/N'])$ as a base.

Next, we define the forgetting map
$$\forg_{n+1}:\wt{\UU}^{ns}_{g,n+1}(\ba,0)\to \wt{\UU}^{ns}_{g,n}(\ba)$$
by associating with the universal curve $(C,p_1,\ldots,p_{n+1})$ over $\wt{\UU}^{ns}_{g,n+1}(\ba,0)$
the algebra $H^0(C\setminus \{p_1,\ldots,p_n\},\OO)$, viewed as a marked algebra of type $\ba$. By Proposition
\ref{marked-alg-prop}(ii) and the first part of the proof, this gives a family in $\wt{\UU}^{ns}_{g,n}(\ba)$, hence
the required morphism. Furthermore, we have a compatible map between universal affine curves coming from
the natural homomorphism of filtered algebras
$$H^0(C\setminus\{p_1,\ldots,p_n\})\to H^0(C\setminus\{p_1,\ldots,p_{n+1}\}).$$
This proves Theorem A(ii).

Finally, let us check that for any $\ba'\in X(g,n)$ the open subset
$\wt{\UU}^{ns}_{g,n}(\ba,\ba')$ is a distinguished open affine in $\wt{\UU}^{ns}_{g,n}(\ba)$.
%Pick $N$ bigger than all $a_i$ and $a'_i$. Then we have a natural morphism of $R$-modules,
%where $\wt{\UU}^{ns}_{g,n}(\ba)=\Spec(R)$,
%$$H^0(C,\OO_C(ND))/H^0(C,\OO_C)\to H^0(C,\OO_C(ND)/\OO_C(a'_1p_1+\ldots+a'_np_n),$$
%where $C$ is the universal curve. Furthermore, both $R$-modules are free of rank $Nn-g$, so taking the determinant
%we get an element $f\in R$, such that $f\neq 0$ is precisely the open subset $\wt{\UU}^{ns}_{g,n}(\ba,\ba')$.
Indeed, $\wt{\UU}^{ns}_{g,n}(\ba,\ba')$ is the preimage of $SG^{\ba'}$
under the Krichever map $\wt{\UU}^{ns}_{g,n}(\ba)\to SG_1(g)$.
But $SG^{\ba'}\cap SG^\ba$ is the distinghuished open affine in $SG^\ba$ by Proposition \ref{Plucker-prop}(ii),
and the assertion follows.
%By definition $SG^{\ba'}$ is given by the inequality $\det(\pi_{\ba'})\neq 0$ (see ???),
%where $\det(\pi_{\ba'})$ is a section of the line bundle $\det(\sV)$. Since on $SG^\ba$
%we have a trivialization of $\det(\sV)$ given by the nonvanishing section $\det(\pi_\ba)\neq 0$, the above inequality is %equivalent to the nonvanishing of the function,
%\begin{equation}\label{double-det-eq}
%\frac{\det(\pi_{\ba'})}{\det(\pi_\ba)}=\det(\pi_\ba^{-1}\circ \pi_{\ba'})\neq 0.
%\end{equation}
\ed

We have used the following result.

\begin{lem}\label{closed-subscheme-lem} Let $S$ be a base scheme, $X$ a stack over $S$, $Y$ a separated scheme over $S$. Suppose we are given morphisms $i:X\to Y$ and $p:Y\to X$ over $S$, such that $p\circ i\simeq \id_X$.
Then $i$ induces an isomorphism of $X$ with a closed subscheme of $Y$.
\end{lem} 

\Pf . Consider the closed subscheme $Z\sub Y$ given by the equation $y=ip(y)$ (here we use separatedness of $Y$).
It is easy to check that the morphisms $i$ and $p$ induce mutually inverse isomorphisms between $X$ and $Z$.
\ed

%This gives the map
%$$\phi:ASG\cap\Sigma\to \UU^{ns,a}_{g,g}.$$
%We have to show that it is the inverse of the Krichever map $\ov{\Kr}$. The fact that $\phi\circ\ov{\Kr}=\id$ follows from 
%the proof of Theorem ???. It remains to prove that $\ov{\Kr}\circ\phi=\id$. 

%If $A,A'\in ASG(\HH)\cap\Sigma$ are isomorphic as marked algebras of genus $g$ then
%in fact $A=A'$.

%Coordinates on $\Sigma$ correspond to restrictions of Plucker coordinates.

%Can we use only $p_{ij}[-n]$ as coordinates???
%Have $\a_{ij}$, $\b_{ij}$, $\eta_{ij}$, $p_{ij}[5]$, $p_{ij}[6]$, etc.
%NO since cannot get $\ga_{ij}$.

%Can write equations for the image of the moduli of curves in terms of these coordinates
%(equations that we have plus express higher $p_{ii}[-n]$ in terms of the smaller ones).

%Homogenizing get equations on Plucker coordinates satisfied by the image of the Krichever map
%(with $g$ distinct points).

\begin{rems} 1. The Gr\"obner basis of the ideal defining the algebra $H^0(C\setminus\{p_1,\ldots,p_n\})$
constructed in the above proof is often not optimal since it uses more generators than needed. 
%One can instead
%for each $i=1,\ldots,n$, $j=0,\ldots,a_i$, choose $\wt{h}_i(j)\in H^0(C,(a_i+1+j)p_i+\sum_{i'\neq i}a_{i'}p_{i'})$
%with the polar part at $p_i$ starting with $t_i^{-a_i-1-j}$. 
%Then a similar procedure to that of Section ??? leads
%to a Groebner basis in terms of these generators (see some examples in ???).
In the examples with $g=1$ considered in Section \ref{examples-sec} we will 
%consider some examples of moduli spaces $\wt{\UU}^{ns}_{g,n}(\ba)$, and in most cases will 
get presentations with a smaller number of generators.

\noindent 2. Let us denote by $\wt{\UU}^{ns}_{g,n}(\ba)'$ the stack of $(C,p_1,\ldots,p_n,v_1,\ldots,v_n)$ such that
$H^1(C,\OO(a_1p_1+\ldots+a_np_n))=0$ (but $\OO(p_1+\ldots+p_n)$ is not necessarily ample).
Let us consider the composition
\begin{equation}\label{bc-map-eq}
\bc: \wt{\UU}^{ns}_{g,n}(\ba)'\to ASG^\ba/\fG\rTo{\sim} \wt{\UU}^{ns}_{g,n},
\end{equation}
where the first map is induced by the Krichever map \eqref{Krichever-map}, while the second is the isomorphism
of Theorem B. The map $\bc$ sends a curve $(C,p_\bullet,v_\bullet)$ to the curve
$(\ov{C},p_\bullet,v_\bullet)$, where
$$\ov{C}=\Proj \bigl(\bigoplus_N H^0(C,\OO(N(p_1+\ldots+p_n)))\bigr).$$
Note that the natural map $C\to\ov{C}$ is an isomorphism near $p_1,\ldots,p_n$, so we have the induced marked points
on $\ov{C}$ and the tangent vectors at them. 
One can check that there is a similarly defined morphism
$$\UU^{r,(\infty)}_{g,n}\to \UU^{r,(\infty)}_{g,n}: (C,p_\bullet,t_\bullet)\mapsto (\ov{C},p_\bullet,t_\bullet),$$
such that $\OO(p_1+\ldots+p_n)$ is ample on $\ov{C}$ and
$$\Kr(C,p_\bullet,t_\bullet)=\Kr(\ov{C},p_\bullet,t_\bullet).$$

\noindent 3. If $(C,p_1,\ldots,p_n,v_1,\ldots,v_n)$ is a curve in $\wt{\UU}^{ns}_{g,n}(\ba)$ and
$(C',q_1,\ldots,q_m,w_1,\ldots,w_m)$ is a curve in $\wt{\UU}^{ns}_{g',m}(\ba')$, such that
the weights of the last marked points $a_n$ and $a'_m$ are both zero, then we can glue
$C$ and $C'$ by identifying $p_n$ and $q_m$ in such a way that we get a node on the glued curve.
It is easy to see that the glued curve $\wt{C}$ with the marked points $p_1,\ldots,p_{n-1},q_1,\ldots,q_{m-1}$ 
satisfies $h^1(\OO_{\wt{C}}(p_1+\ldots+p_{n-1}+q_1+\ldots+q_{m-1}))=0$, however, it may happen that it has
irreducible components without marked points. We can contract them by applying 
the map \eqref{bc-map-eq}.
The obtained curve, equipped with the tangent vectors induced by $v_1,\ldots,v_{n-1},w_1,\ldots,w_{m-1}$ defines a point
of $\wt{\UU}^{ns}_{g+g',n+m-2}(\bb)$, where $\bb=(a_1,\ldots,a_{n-1},a'_1,\ldots,a'_{m-1})$. 
Applying this procedure to the universal curves, we get a morphism
$$\wt{\UU}^{ns}_{g,n}(\ba)\times \wt{\UU}^{ns}_{g',m}(\ba')\to \wt{\UU}^{ns}_{g+g',n+m-2}(\bb).$$
\end{rems}
%Should also study the relation between the coordinates on $\Sigma$ and the tau-functions/theta-functions.

\section{More on the moduli schemes $\wt{\UU}^{ns}_{g,n}(\ba)$}\label{more-moduli-sec}

\subsection{Some special curves}\label{special-curves-subsec}

Everywhere in Section \ref{more-moduli-sec} we work over $\Q$.

First, we are going to construct some special singular curves corresponding to points of $\wt{\UU}^{ns}_{g,n}(\ba)$.
%In the case when all $a_i$'s are either $0$ or $1$ we will use these curves in Section \ref{big-moduli-sec}.
The cuspidal curves from the following definitions will appear as irreducible components of some of our special curves.

\begin{defi} For each $g\ge 0$ let us denote by $C^{\cusp}(g)$ the following rational cuspidal curve.
The underlying topological space is $\P^1$, and the structure sheaf of $C^{\cusp}(g)$ is the subsheaf of
$\OO_{\P^1}$ consisting of all functions $f$ such that near $0\in\P^1$ one has
$f-f(0)\in\mg^{g+1}$, where $\mg$ is the maximal ideal in the local ring of $0$. Note that for $g=0$ we
have $C^{\cusp}(0)=\P^1$.
\end{defi}

Without loss of generality let us assume that the weights $\ba=(a_1,\ldots,a_n)$ satisfy
$a_i>0$ for $i=1,\ldots,r$, and $a_i=0$ for $r<i\le n$.

\begin{defi} Let $x_1,\ldots,x_n$ be independent variables, and let $k$ be a field. First, we define a subalgebra in
$\bigoplus_{i=1}^rk[x_i]$ by
$$B=B(a_1,\ldots,a_r):=k\cdot 1+\bigoplus_{i=1}^r x_i^{a_i+1}k[x_i].$$
Next, let $(\ov{h}_j)_{j=r+1,\ldots,n}$ be a collection of elements in $\bigoplus_{i=1}^r x_ik[x_i]/(x_i^{a_i+1})$
with the property $\ov{h}_j\ov{h}_{j'}=0$ for $j\neq j'$.
Now we define $A(\ov{h}_\bullet)$ as the $B$-subalgebra in $\bigoplus_{i=1}^n k[x_i]$ 
generated by the elements $h_j=x_j+\ov{h}_j$, $j=r+1,\ldots,n$, 
where we lift each $\ov{h}_j$ to a sum of polynomials of degree $\le a_i$ in $x_i$
with no constant terms.
%Let 
%$$\Ga\sub\bigoplus_{i=1}^r x_ik[x_i]/(x_i^{a_i+1})\oplus\bigoplus_{j=r+1}^nx_ik[x_j]/x_j^2k[x_j]$$
%be the graph of the linear map given by $x_j\mapsto \ov{h}_j$. Note that $\Ga\cdot \Ga=0$.
%Then $A$ is the preimage
%of $1+\Ga$ under the natural projection
%$$\bigoplus_{i=1}^n k[x_i]\to \bigoplus_{i=1}^r x_ik[x_i]/(x_i^{a_i+1})\oplus\bigoplus_{j=r+1}^nx_ik[x_j]/x_j^2k[x_j].$$
%In the case $r=n$ no additional data is necessary and we just set
%$$A=1+\bigoplus_{i=1}^n x_i^{a_i+1}k[x_i].$$
\end{defi}

We equip the above algebra $A(\ov{h}_\bullet)$ with the $M_\ba$-valued filtration
$$F_{\bm}A(\ov{h}_\bullet)=A(\ov{h}_\bullet)\cap \bigoplus_{i=1}^n k[x_i]_{\le m_i},$$
where $k[x_i]_{\le d}$ denotes the space of polynomials of degree $\le d$.
We denote by $\RR(A(\ov{h}_\bullet))$
the Rees algebra associated with the filtration $F_m A(\ov{h}_\bullet)=F_{m\bD_1}A(\ov{h}_\bullet)$.

%increasing filtration $F_\bullet A(\ov{h}_\bullet)$ given by the maximum of
%degrees in $x_i$. 

\begin{prop}\label{special-curves-prop}
(i) The algebra $A(\ov{h}_\bullet)$ with the above filtration extends to a structure of a marked algebra of type $\ba$.
%Basis of $A(\ov{h}_\bullet)$???
Hence, the curve 
$$C(\ov{h}_\bullet):=\Proj(\RR(A(\ov{h}_\bullet)))$$ 
defines a point of $\wt{\UU}^{ns}_{g,n}(\ba)$.

\noindent
(ii) The curve $C(\ov{h}_\bullet)$
is the union of $n$ irreducible components $C_i$, 
joined in a single point $q$, which is the only singular point
of $C(\ov{h}_\bullet)$ 
(with $p_i\in C_i\setminus \{q\}$). The components $C_i$ corresponding to $i>r$ are all isomorphic to $\P^1$.
The component $C_i$ associated with $i\le r$ is the irreducible rational curve given as $\Proj$ of the Rees algebra of
the subalgebra in $k[x_i]$ generated by $\bigoplus x_i^{a_i+1}k[x_i]$ and by the $i$th components of $\ov{h}_j$,
$(\ov{h}_j)_i=\sum_{p=1}^{a_i}c_{ji}(p)x_i^p$, for $j=r+1,\ldots,n$. 
%where
%$\ov{h}_j=\sum_{i=1}^r\sum_{p=1}^{a_i} c_{ji}(p)x_i^p$.
Furthermore, for $i=1,\ldots,n$
the canonical formal parameter at $p_i\in C_i$ (see Corollary \ref{curves-sec-cor}) is given by $x_i^{-1}$
(for any choice of $i_0\le r$).

\noindent
(iii) The algebra $A(\ov{h}_\bullet)$ has the following description by generators and relations: generators
$h_i(m)=x_i^{a_i+1+m}$, for $i=1,\ldots,r$, $m=0,\ldots,a_i$; $h_j$, for $j=r+1,\ldots,n$, and relations
\begin{align} 
\begin{split}\label{special-gen-curve-eq}
&h_i(m)h_{i'}(m')=0, \\ 
&h_i(m)h_i(m')=h_i(m+m')h_i(0), \\
&h_i(m)h_j=\sum_{p=1}^{a_i}c_{ji}(p)h_i(m+p), \\
&h_jh_{j'}=\sum_{i=1}^r\sum_{m\ge 0}\sum_{p,p'\ge 1, p+p'=a_i+1+m}c_{ji}(p)c_{j',i}(p')h_i(m)
\ \text{ for } j\neq j',
\end{split}
\end{align} 
where $i\neq i'$, $j\neq j'$, $i,i'\le r$, $j,j'>r$, and we set $h_i(m)=h_i(m-a_i-1)h_i(0)$ for $a_i+1\le m\le 2a_i$.
\end{prop}

\Pf . (i) Let us consider the ideal $I=\bigoplus_{i=1}^r x_i^{a_i+1}k[x_i]\sub A(\ov{h}_\bullet)$.
It is easy to see that we have an exact sequence
$$0\to I\to A(\ov{h}_\bullet)\rTo{p_{>r}} B'\to 0,$$
where $B'=k\cdot 1+\bigoplus_{j=r+1}^n x_jk_[x_j]\sub \bigoplus_{j=r+1}^n k[x_j]$, the map
$p_{>r}$ is induced by the natural projection $\bigoplus_{i=1}^n k[x_i]\to \bigoplus_{j=r+1}^n k[x_j]$.
%This immediately implies that $F_{(a_1,\ldots,a_r,0,\ldots,0)}A(\ov{h}_\bullet)=k\cdot 1$.
Furthermore, since $p_{>r}(h_j^m)=x_j^m$ for $m>0$, it follows that
for any $\bm\ge \ba$ we have an exact sequence 
\begin{equation}\label{A-h-filtration-exact-sequence}
0\to I_{\le (m_1,\ldots,m_r)}\to F_\bm A(\ov{h}_\bullet)\to B'_{\le (m_{r+1},\ldots,m_n)}\to 0,
\end{equation}
where $B'_{\le (m_{r+1},\ldots,m_n)}$ consists of $(f_j)\in B'$ such that $\deg(f_j)\le m_j$,
and similarly for $I_{\le (m_1,\ldots,m_r)}$.
%Indeed, the surjectivity of the last map follows from the fact that for any $p\ge 1$ one has
%$$h_j^p\equiv \sum_{i=1}^r g_i(x_i)+x_j^p \mod I$$
%for some polynomials $g_i(x_i)$ such that $\deg(g_i)\le a_i$. 
The exact sequences \eqref{A-h-filtration-exact-sequence} easily imply that we indeed have a structure
of a marked algebra on $A(\ov{h}_\bullet)$.

\noindent
(ii) Let us consider the affine curve $C^{aff}=\Spec A(\ov{h}_\bullet)$, which is dense in $C(\ov{h}_\bullet)$.
By definition, we have a surjective finite morphism $\sqcup_{i=1}^n\A^1\to C^{aff}$
corresponding to the embedding of $A(\ov{h}_\bullet)$ into $\bigoplus_{i=1}^n k[x_i]$.
Hence, $C^{aff}$ is the union of $n$ irreducible components $C^{aff}_i$, where the ring of functions
on $C^{aff}_i$ is the image of the projection $A(\ov{h}_\bullet)\to k[x_i]$. This implies the assertion about the
irreducible components $C_i$. For $i=1,\ldots,r$, 
the fact that the canonical formal parameter at $p_i$ is induced by $x_i^{-1}$, follows from the fact that 
for each $m\ge a_i+1$ 
the element $x_i^m$ extends to a regular function on $C^{aff}$. For $j>r$ the identification of the canonical
formal parameter at $p_j$ follows from the fact that
$h_j$ is a regular function on $C^{aff}$, such that $h_j|_{C^{aff}_j}=x_j$ and the expansion of $h_j$ in $x_i$
has no constant term for any $i=1,\ldots,r$.

\noindent
(iii) It is easy to check that equations \eqref{special-gen-curve-eq} 
hold in $A(\ov{h}_\bullet)$. The fact that these are defining relations
follows from the fact that $A(\ov{h}_\bullet)$ is a marked algebra of type $\ba$.
\ed

\begin{exs}\label{special-curves-ex}
1. If we take all $\ov{h}_j$ to be zero then the corresponding algebra is
$$B(a_1,\ldots,a_n)=k\cdot 1+\bigoplus_{i=1}^n x_i^{a_i+1}k[x_i] \sub \bigoplus_{i=1}^n k[x_i].$$
The corresponding curve is the pointed transversal union of the cuspidal curves $C^{\cusp}(a_i)$:
\begin{equation}\label{a-cuspidal-curve-eq}
C^{\cusp}(\ba):=\cup_{i=1}^n C^{\cusp}(a_i),
\end{equation}
where we identify the origins in all $C^{\cusp}(a_i)$, and take the infinity on each component to be the marked point.
Note that by Proposition \ref{special-curves-prop}(ii), all the coordinates $\a_{ij}[p,q]$ vanish at the corresponding
point of $\wt{\UU}^{ns}_{g,n}(\ba)$.

\noindent
2. In the case $a_1=\ldots=a_r=1$ the elements $\ov{h}_{r+1},\ldots,\ov{h}_n$ can be arbitrary
linear combinations of $x_1,\ldots,x_r$. 
For $i=1,\ldots,r$, the irreducible component $C_i$ is isomorphic to the cuspidal curve $C^{\cusp}(1)$ if and only if
$(\ov{h}_j)_i=0$ for all $j>r$. Otherwise, $C_i\simeq \P^1$.
It is easy to see that the points in the moduli space $\wt{\UU}^{ns}_{r,n}(1,\ldots,1,0,\ldots,0)$, corresponding
to these curves, are invariant with respect
to the action of the diagonal subgroup $\G_m\sub\G_m^n$.
% (see Theorem \ref{big-moduli-thm}).

\noindent
3. In the case $g=1$, $r=1$, $a_1=1$, $n\ge 3$, if we take
$\ov{h}_j=x_1$ for $j=2,\ldots,n$, then the equations \eqref{special-gen-curve-eq} give 
$h_jh_{j'}=h_1(0)$ for $2\le j<j'\le n$ and $h_1(0)h_j=h_1(1)$, hence,
$h_2,\ldots,h_n$ are generators of the corresponding algebra, 
and the only equations on them are that the products $h_jh_{j'}$, where $j\neq j'$, do not depend on the pair $(j,j')$, and that
$h_jh_{j'}^2=h_j^2h_{j'}$ (the latter equations are superfluous for $n>3$).
In other words, we get precisely the elliptic $n$-fold curve (see \cite[Sec.\ 2]{SmythI}, \cite[Sec.\ 1.5]{LP}).
Similarly, in the case $g=1$, $r=1$, $a_1=1$, $n=2$, and $\ov{h}_2=x_1$ we get $h_1(0)h_2=h_1(1)$,
so $h_1(0)$ and $h_2$ are generators with the only relation $h_1(0)h_2^2=h_1(0)^2$, which is the equation
of the tacnode, i.e., the elliptic $2$-fold singularity.
%h_1(0)h_2=h_1(1), h_1(1)h_2=h_1(0)^2
\end{exs}

%Consider coordinates $\a_{ij}[-1,-p]$ for $i>r$, $j\le r$, $1\le p\le a_j$.
%Should have a curve depending only on these???

\subsection{The $\G_m^n$-action on $\wt{\UU}^{ns}_{g,n}(\ba)$}\label{subgroup-sec}

Below we freely use the identification of
$\wt{\UU}^{ns}_{g,n}(\ba)$ with the closed subset $\Si^{\ba,i_0}\cap ASG$ of the open
cell $SG^\ba$ in the Sato Grassmannian (depending on a choice of $i_0$ such that $a_{i_0}>0$). 
In particular, we view the coordinates 
$\a_{ij}[p,q]$ (see \eqref{coordinates-a-ij-eq}) as functions on $\wt{\UU}^{ns}_{g,n}(\ba)$. 

Recall (see \cite{P-ainf}) that the moduli space $\wt{\UU}^{ns}_{g,g}(1,\ldots,1)$ can be viewed as the deformation space of the
singular curve $C^{\cusp}(1,\ldots,1)$.
that corresponds to the unique $\G_m^g$-invariant point of the moduli space. 
%namely, the union of $g$ cuspidal curves of arithmetic genus $1$, joined transversally at the singular point.
Furthermore, this curve is the limit of the $\G_m$-orbit through every other point of $\wt{\UU}^{ns}_{g,g}(1,\ldots,1)$
with respect to the diagonal subgroup $\G_m\sub\G_m^g$.

In this section we will show that a similar picture holds for arbitrary $\ba\in X(g,n)$ with an appropriate choice of a subgroup
$\G_m$ in $\G_m^n$.

The functions $\a_{ij}[p,q]$ on $\wt{\UU}^{ns}_{g,n}(\ba)$ are semi-invariant with respect
to the $\G_m^n$-action, with the weights
%\begin{equation}\label{a-ij-wts-eq}
$$wt(\a_{ij}[p,q])=-p\be_i+q\be_j, \ \ p\le -a_i-1, q\ge -a_j.
$$
%\end{equation}
Let us denote by $\Om_\ba\sub\Z^n$ the set of these weights.

\begin{defi}\label{cone-defi} For $n>1$ we define $\bC_\ba\sub\R^n$ to be the closed cone 
generated by all the vectors $\om_{ij}=(a_i+1)\be_i-a_j\be_j$, with $i\neq j$, i.e.,
the set of all linear combination of $\om_{ij}$ with coefficients in $\R_{\ge 0}$. In the case $n=1$ (and
$\ba=(g)$) we set $\bC_\ba=\R_{\ge 0}$.
\end{defi}

\begin{lem}\label{cone-generation-lem}
The cone $\bC_\ba$ contains all the basis vectors $\be_i$. Hence, we have $\Om_\ba\sub\bC_\ba$.
In the case when $a_i>0$ for at least two indices $i$ 
the cone $\bC_\ba$ is generated by the vectors $\om_{ij}$ such that $i\neq j$
and $a_j>0$. In the case when $a_j=0$ for all $j\neq 1$, i.e., $\ba=g\be_1$, the cone $\bC_{\be_1}$ 
is generated by $\be_1$ and by the vectors $\om_{i1}$, for $i=2,\ldots,n$.
\end{lem}

\Pf . The first assertion follows from the equality 
\begin{equation}\label{om-e-eq}
(a_j+1)\om_{ij}+a_j\om_{ji}=(1+a_i+a_j)\be_i.
\end{equation}
For the second assertion let us denote by $\bC'_\ba$ the cone generated by $\om_{ij}$ with $i\neq j$ and $a_j>0$.
To see that $\bC'_\ba=\bC_\ba$ it is enough to check that $\bC'_\ba$ contains all the basis vectors $\be_i$.
In the case when $a_i>0$ for at least two $i$, \eqref{om-e-eq} implies that this is true
for each $i$ such that $a_i>0$. Hence, applying \eqref{om-e-eq} to a pair $(i,j)$ such that $a_i=0$ and $a_j>0$
we see that it is also true for $i$ such that $a_i=0$. The assertion about generators of the cone in the case $\ba=\be_1$ is
straightforward. 
\ed
%it contains the weights \eqref{a-ij-wts-eq}. 

%More generally, let us consider any collection $\ba(1),\ldots,\ba(r)\in X(g,n)$ with non-intersecting supports
%(i.e., for each $i=1,\ldots,n$ there is at most one $j=1,\ldots,r$ with $a_i(j)>0$).
%Then the affine scheme
%$\wt{\UU}^{ns}_{g,n}(\ba(1),\ldots,\ba(r))$ is the open subset in $\wt{\UU}^{ns}_{g,n}(\ba(1))$,
%where the determinants of the morphisms
%$$\pi_{\ba(1)}^{-1}\circ \pi_{\ba(i)}:\HH(\ba(i))/\HH_{\ge 0}\ot\OO\to \HH(\ba(1))/\HH_{\ga 0}\ot \OO,$$
%$i=2,\ldots,r$, do not vanish (see \eqref{double-det-eq}). Let us set 
%$$D_{\ba(1),\ba(i)}:=\det(\pi_{\ba(1)}^{-1}\circ \pi_{\ba(i)}).$$
%Note that the weights of these determinants are given by
%$$wt(D_{\ba(1),\ba(i)})=\sum_j \ba_j(i)e_j - \sum_j \ba_j(1)e_j.$$

%Need rational numbers 
%$w_1,\ldots,w_n$ such that for any coordinate on the moduli space ??? with the $\G_m^n$-weight $(m_1,\ldots,m_n)$
%one has $\sum w_im_i>0$. Note: if all $w_i\ge 0$ then need only coordinates with minimal $\G_m^n$-weights, which
%have form $(a_i+1)e_i-a_je_j$. 
For a given $\ba\in X(g,n)$, let us set
$$N=\max_{1\le i\le n}(a_i),$$
$$w_i=\begin{cases} \frac{1}{a_i}, & a_i>0\\ 1+\frac{1}{N}, & a_i=0.\end{cases},$$
$i=1,\ldots,n$.
%$\eps$ is a positive rational number, and let 
Let also $\ell$ be the linear function
on $\R^n$ given by $\ell(\be_i)=w_i$. 

%if both $a_i$ and $a_j$ are nonzero, we get $\frac{a_i+1}{a_i}-1>0$.
%If $a_i\neq 0$ and $a_j=0$ then we get $\frac{a_i+1}{a_i}>0$.
%If $a_j\neq 0$ and $a_i=0$, we get $2-1=1>0$. 

\begin{lem}\label{global-functions-lem} 
(i) For any $i,j$ one has $\ell(\om_{ij})\ge\frac{1}{N}$. Hence, the set $\Om_\ba$ lies in the half-space $\ell>0$.

\noindent (ii) There are no global $\G_m^n$-invariant functions on $\wt{\UU}^{ns}_{g,n}(\ba)$.
\end{lem}

\Pf . (i) This is immediate from the definitions (note that $\frac{1}{N}\le 1$ since at least one $a_i$ is $\ge 1$).

\noindent
(ii) This follows from the fact that any function on $\wt{\UU}^{ns}_{g,n}(\ba)$
is a linear combination of monomials in coordinates $\a_{ij}[p,q]$, and that any nonconstant monomial 
in these coordinates has $\G_m^n$-weight in the half-space $\ell>0$ by part (i).
\ed

%Recall Groebner basis picture.
%PROBLEM: MAY GET SOME NEGATIVE WEIGHTS

Let us consider the embedding
$$r\bw:\G_m\to \G_m^n:\la\mapsto (\la^{rw_1},\ldots,\la^{rw_n})$$
given by the multiples of the weights $w_i$, where we
choose a rational number $r>0$ so that $(rw_1,\ldots,rw_n)$ are coprime integers.

\begin{prop}\label{rw-invariant-prop}
The moduli scheme $\wt{\UU}^{ns}_{g,n}(\ba)$ has a unique $r\bw(\G_m)$-invariant point, namely, the curve
$C^{\cusp}(\ba)$ (see \eqref{a-cuspidal-curve-eq}).
The closure of the $r\bw(\G_m)$-orbit through every point in $\wt{\UU}^{ns}_{g,n}(\ba)$ contains the invariant point.
Thus, the corresponding curve can be degenerated to $C^{\cusp}(\ba)$.
%where $b_i$??? 
%(we identify the points $0\in C^{\cusp}(b_i)$).
%Should be similar to Prop. 2.1.5
\end{prop}

\Pf . By Lemma \ref{global-functions-lem}(i), all the coordinates $\a_{ij}[p,q]$ have positive weight with respect
to the subgroup $r\bw(\G_m)$. Hence, there is a unique $r\bw(\G_m)$-invariant point, namely, the point where
all these coordinates vanish, and it belongs to the closure of every orbit.
\ed
%For the proof need to analyze our Groebner bases.???

%For example, if $n=2$, $a_1=1$, $a_2=0$, then we get the ordinary cuspidal curve joined with $\P^1$ at the cusp.

%In the case when some $a_i$'s vanish, it makes sense to consider 

\subsection{Generators of the ring of functions on $\wt{\UU}^{ns}_{g,n}(\ba)$}\label{generators-sec}

The set of coordinates $(\a_{ij}[p,q])$ is superfluous: one can express some in terms of the others.
In this section we find a more convenient 
(infinite) subset of these coordinates that still generate the ring $\OO(\wt{\UU}^{ns}_{g,n}(\ba))$. 
We will later use this generating subset in analyzing the GIT stability for the $\G_m^n$-action on $\wt{\UU}^{ns}_{g,n}(\ba)$.

%Fix $\ba\in X(g,n)$.
%Recall (see ???) that we have an isomorphism of $\wt{\UU}^{ns}_{g,n}(\ba)$ with the closed subset
%$\wt{\UU}^{ns}_{g,n}(\ba)\simeq ASG\cap \Si^{\ba,i_0}$, so we can view the
%coordinates $(\a_{ij}[p,q])$ as functions on $\wt{\UU}^{ns}_{g,n}(\ba)$.

In the next result we normalize $(\a_{ij}[p,0])$ by $\a_{ii}[p,0]=0$.

\begin{prop}\label{generators-prop} 
The coordinates 
\begin{equation}\label{generators-wide-eq}
\{\a_{ii}[p,1] \ | \ p\le -a_i-1, a_i>0\} \  \text { and } \ \{\a_{ij}[p,q]) \ |\ p\le -a_i-1,-a_j\le q\le 0\}
% \ \text{ and } (\a_{ij}[p,0]-\a_{ij'}[p,0])_{p\le -a_i-1})
\end{equation}
(where $i$ and $j$ don't have to be distinct), 
generate the ring $\OO(\wt{\UU}^{ns}_{g,n}(\ba))$. 
\end{prop}

Let $A\sub \OO(\wt{\UU}^{ns}_{g,n}(\ba))$ be the subring generated by \eqref{generators-wide-eq}.

\begin{lem}\label{coefficients-lem}
For $m\ge 1$, $m'\ge 1$, $i\neq j$, one has
\begin{equation}\label{ff-expansion-C-eq}
f_i[-a_i-m]f_j[-a_j-m']=\sum_{k,q\le -a_k-1} c_{k,q} f_k[q]+C,
\end{equation}
where $C$ is a constant, and the coefficients $c_{k,q}$ for $k\neq i,j$ belong to $A$. 
Also, for $l\ge 1$ we have 
$$c_{i,-a_i-l}\equiv \a_{ji}[-a_j-m',m-l] \mod A, \ \ c_{j,-a_j-l}\equiv \a_{ij}[-a_i-m,m'-l] \mod A.$$
Similarly, for $m\ge 1$, $m'\ge 1$ one has
\begin{equation}\label{ff-ii-expansion-C-eq}
f_i[-a_i-m]f_i[-a_i-m']=\sum_{k,q\le -a_k-1} c_{k,q} f_k[q]+C,
\end{equation}
where $c_{k,q}\in A$ for $k\neq i$, and
$$c_{i,-a_i-l}\equiv \a_{ii}[-a_i-m',m-l]+\a_{ii}[-a_i-m,m'-l] \mod A.$$
%If $i=j$ the coefficients $c_{i,q}$ depend on functions in $\a_{ii}[-a_i-m,n']$ with $n'<n$,
%$\a_{ii}[-a_i-n,m']$ with $m'<m$ and on functions in $A$.
\end{lem}

\Pf . We will only consider the case $i\neq j$; the case $i=j$ is similar.
It is enough to choose coefficients $c_{k,q}$ in such a way that the difference
$$f_i[-a_i-m]f_j[-a_j-m']-\sum_{k,q\le -a_k-1} c_{k,q} f_k[q]$$
belongs to $H^0(C,\OO(\sum a_kp_k))$. Since $f_k[-a_k-l]$ has poles of order $\le a_{k'}$ for $k'\neq k$
and has expansion $t_k^{-a_k-l} \mod t_k^{-a_k+1}k[[t_k]]$ at $p_k$, we see that we can take
$c_{k,-a_k-l}$ to be the coefficient of $t_k^{-a_k-l}$ in the expansion of $f_i[-a_i-m]f_j[-a_j-m']$.
For $k\neq i,j$ we have
$$c_{k,-a_k-l}=\sum_{q,q'\ge -a_k;q+q'=-a_k-l}\a_{ik}[-a_i-m,q]\a_{jk}[-a_j-m',q'],$$
which is in $A$ since $q$ and $q'$ are in $[-a_k,-1]$.
For $k=i$ we have
$$c_{i,-a_i-l}=\a_{ji}[-a_j-m',m-l]+\sum_{q,q'\ge -a_i;q+q'=-a_i-l}\a_{ii}[-a_i-m,q]\a_{ji}[-a_j-m',q'],$$
and it remains to observe again that $q$ and $q'$ are in $[-a_i,-1]$. 
%$\a_{ik}[p,q_1]$ with $q_1<0$ and $\a_{jk}[p',q_2]$ with $q_2<0$
\ed

\noindent
{\it Proof of Proposition \ref{generators-prop}}.

\noindent

\noindent
{\bf Step 1}. For $a_i>1$ one has $\a_{ii}[p,q]\in A$. We will prove by induction on $m+l$ that for
$m\ge 1$, $l\ge 1$ one has $\a_{ii}[-a_i-m,l]\in A$. In the case when $l\le 1$ this holds by the definition.
Assume that the assertion is true for $m'+l'<m+l$. 
Consider the expression \eqref{ff-ii-expansion-C-eq} of $f_i[-a_i-m]f_i[-a_i-l]$. Then the coefficients of all nonconstant terms
belong to $A$ by Lemma \ref{coefficients-lem} and by the induction assumption. 
Now considering the coefficients of $t_i^{-a_i}$ and $t_i^{-a_i+1}$ we get
$$\a_{ii}[-a_i-m,l]\equiv -\a_{ii}[-a_i-l,m] \mod A,$$
$$\a_{ii}[-a_i-m,l+1]\equiv -\a_{ii}[-a_i-l,m+1] \mod A.$$
Combining these conditions we derive that for $l>1$ one has
$$\a_{ii}[-a_i-m,l]\equiv \a_{ii}[-a_i-m-1,l-1],$$
which implies by induction on $l$ that $\a_{ii}[-a_i-m,l]\in A$.

\noindent 
{\bf Step 2}. If for some $i\neq j$, such that $a_i>0$, and some
$m\ge 1$, one has $\a_{ij}[p,m']\in A$ for $m'<m$ and all $p\le -a_i-1$, then
$\a_{ji}[-a_j-m,q]\in A$ for all $q\ge -a_i$. 

Indeed, for $q\le 0$ this holds by definition, so we may assume that $q>0$.
Let us consider the expression \eqref{ff-expansion-C-eq} for the
 product $f_i[-a_i-q]f_j[-a_j-m]$. By Lemma \ref{coefficients-lem}, the coefficients of $f_k[p]$ for $k\neq i,j$
are in $A$, while the coefficient of $f_j[-a_j-l]$ is $\a_{ij}[-a_i-q,m-l]\mod A$, which is also in $A$ by assumption.
Since the elements $f_i[q']$ all have zero coefficient of $t_i^{-a_i}$ by the definition of $\Si^\ba$, we deduce that
the coefficient of $t_i^{-a_i}$ in $f_i[-a_i-q]f_j[-a_j-m]$ belongs to $A$. But this coefficient is equal to 
$$\a_{ji}[-a_j-m,q]+\sum_{q_1,q_2\ge -a_i;q_1+q_2=-a_i}\a_{ii}[-a_i-q,q_1]\a_{ji}[-a_j-m,q_2],$$
and our claim follows.

\noindent
{\bf Step 3}. Now let us prove by induction on $m\ge 0$ that for $i\neq j$, such that $a_i\neq 0$, one has
$\a_{ij}[p,m]\in A$ and $\a_{ji}[-a_j-1-m,q]\in A$ for all $p\le -a_i-1$,
$q\ge -a_i$. The base of induction follows from Step 2 (applied to $m=1$).
Suppose now that $m\ge 1$ and the assertion is true for $m'<m$. It suffices to check that $\a_{ij}[-a_i-m',m]\in A$.
The second assertion would follow by applying Step 2 again (to $m+1$). Let us consider the product
$f_i[-a_i-m']f_j[-a_j-m]$. By Lemma \ref{coefficients-lem} and by the induction assumption,
the coefficients of the nonconstant terms of the expansion \eqref{ff-expansion-C-eq} of this product 
depend only on elements of $A$. In the case when $a_j>0$ we deduce that the coefficient of 
$t_j^{-a_j}$ in the expansion of $f_i[-a_i-m']f_j[-a_j-m]$ belongs to $A$. But this coefficient
is equal to 
$$\a_{ij}[-a_i-m',m]+\sum_{k=1}^{a_j}\a_{ij}[-a_i-m',-k] \a_{jj}[-a_j-m,-a_j+k].$$
Hence, we derive that $\a_{ij}[-a_i-m',m]\in A$.
In the case $a_j=0$ we need the constant term $C$ in the expression \eqref{ff-expansion-C-eq} of
$f_i[-a_i-m']f_j[-m]$. Looking at the constant term in the expansion in $t_i$, we see that
$$C\equiv \a_{ji}[-m,a_i+m']+\sum_{-a_i\le q\le a_i}\a_{ii}[-a_i-m,q]\a_{ji}[-m,-q] \mod A.$$
Note that by Step 2 we have $\a_{ii}[-a_i-m,q]\in A$ for $q\le a_i$, while $\a_{ji}[-m,q']\in A$
by the induction assumption. Hence, $C\in A$. Then looking at the constant term in the expansion in $t_j$,
as above we deduce that $\a_{ij}[-a_i-m',m]\in A$.

\noindent
{\bf Step 4}.
Thus, we proved that whenever $i\neq j$, and either $a_i>0$ or $a_j>0$, then one has $\a_{ij}[p,q]\in A$.
Assume now that $a_i=a_j=0$ for some $i\neq j$. Let $i_0$ be such that $a_{i_0}\neq 0$.
We prove by induction on $m+l$, where $m\ge 1$, $l\ge 0$, that $\a_{ij}[-m,l]\in A$. 
The base case $m=1$, $l=0$ is clear. Assume the assertion holds for $m'+l'<m+l$.
We can assume that $l>0$, otherwise the assertion holds by definition.
Consider the expression \eqref{ff-expansion-C-eq} of $f_i[-m]f_j[-l]$. Then the coefficients of all nonconstant terms
belong to $A$
by Lemma \ref{coefficients-lem} and by the induction assumption. Looking at the expansion in $t_{i_0}$ and
using the previously proved case we see that the constant coefficient also belongs to $A$.
Now looking at the constant coefficient in the expansion in $t_j$ we derive that $\a_{ij}[-m,n]\in A$.

\noindent
{\bf Step 5}. Assume that $a_i=1$ and that for some $m\ge 1$, $l\ge 1$ one has $\a_{ii}[-1-m',l']\in A$ for $m'+l'<r$.
Then 
\begin{equation}\label{a-1-congruences-eq}
\begin{array}{l}
\a_{ii}[-1-m,l]\equiv -\a_{ii}[-1-l,m] \mod A,\\
\a_{ii}[-1-m,l+2]\equiv -\a_{ii}[-1-l,m+2] \mod A.
\end{array}
\end{equation}
Indeed, as in Step 1, we observe that the coefficients of nonconstant terms in the expression 
\eqref{ff-expansion-C-eq} of $f_i[-1-m]f_i[-1-l]$ belong to $A$ by Lemma \ref{coefficients-lem} and by the 
assumption. Now the assertion follows by considering the coefficients of $t_i^{-1}$ and $t_i$ in
$f_i[-1-m]f_i[-1-l]$ (and using the fact that $\a_{ji}[p,q]\in A$ for $j\neq i$).

\noindent
{\bf Step 6}. By the previous step, we immediately see by induction in $m\ge 1$ that for $a_i=1$ one has $\a_{ii}[-2,m]\in A$
(using the first of the congruences \eqref{a-1-congruences-eq}). Next, we will prove that for $a_i=1$ one has 
$\a_{ii}[-3,q]\in A$.
First, using the previous step we easily derive that $\a_{ii}[-3,3]$ and $\a_{ii}[-4,2]$ are in $A$.
Then applying Lemma \ref{coefficients-lem} we find that 
$$f_i[-2]^2=f_i[-4]+a+\ldots,$$
$$f_i[-2]f_i[-4]=f_i[-6]+\a_{ii}[-2,2]f_i[-2]+b+\ldots,$$
$$f_i[-3]^2=f_i[-6]+2\a_{ii}[-3,1]f_i[-2]+c+\ldots,$$
where the skipped terms are linear combinations of $f_j[p]$, $j\neq i$, with coefficients in $A$, and the constants $a,b,c$
satisfy
$$a\equiv 2\a_{ii}[-2,2]\mod A, \ b\equiv \a_{ii}[-2,4]+\a_{ii}[-4,2]\mod A, \ c\equiv 2\a_{ii}[-3,3]\mod A.$$
Thus, by the previous work, we have $a,b,c\in A$. Since the expansion of $f_i[-2]$ in $t_i$ has coefficients in $A$,
we derive that the same is true for $f_i[-4]$ and for $f_i[-6]$ (using that $\a_{ji}[p,q]\in A$ for $j\neq i$).
Therefore, the expansion of $f_i[-3]^2$ in $t_i$ has coefficients in $A$, and hence, the same is true for $f_i[-3]$.

\noindent
{\bf Step 7}. Applying Step 5, we see by induction on $m\ge 1$ that for $a_i=1$ one has $\a_{ii}[-m,2]\in A$ 
(using the first of the congruences \eqref{a-1-congruences-eq}). Next, let us prove by induction on $m+l$ that 
$\a_{ii}[-1-m,l]\in A$ for $m\ge 1$, $l\ge 1$. Assume this is true for $m'+l'<m+l$. Then applying 
Step 5, and combining the resulting congruences \eqref{a-1-congruences-eq} we derive that for $l>2$ one has
$$\a_{ii}[-1-m,l]\equiv \a_{ii}[-1-m-2,l-2].$$
Since we already know that $\a_{ii}[p,1],\a_{ii}[p,2]\in A$,
this implies by induction on $l$ that $\a_{ii}[-a_i-m,l]\in A$.

\noindent
{\bf Step 8}. Finally, assume that $a_i=0$ and let us prove that $\a_{ii}[-m,q]\in A$ by induction on $m\ge 1$.
In the base case $m=1$ we have $\a_{ii}[-1,q]=0$ for $q\ge 0$, by the definition of the section
$\Sigma^{\ba,i_0}$ and since we normalized the coordinates by $\a_{ii}[p,0]=0$. 
Assume that $\a_{ii}[-m',q]\in A$ for $m'<m$. By Lemma \ref{coefficients-lem}, we have
\begin{equation}\label{ff-a-0-expansion-eq}
f_i[-1]f_i[-m+1]=f_i[-m]+\sum_{l=1}^{m-1} c_{-l}f_i[-l]+C+\ldots,
\end{equation}
where $c_{-l}\in A$ and the skipped terms are linear combinations of $f_j[p]$, $j\neq i$, with coefficients in $A$.
Looking at the constant term in the expansion in $t_i$ and using the induction assumption we see that $C\in A$.
Hence, we can use \eqref{ff-a-0-expansion-eq} to find the entire expansion of $f_i[-m]$ in $t_i$.
\ed

%If $i=j$ then the coefficient of $t_i^{-a_i}$ in the expansion of $f_i[-a_i-m]f_i[a_i-n]$ is equal to
%$$\a_{ii}[-a_i-m,n]+\a_{ii}[-a_i-n,m]+\ldots$$
%Coef-t of $t_i^{-a_i+1}$ in $f_i[-a_i-m]f_i[-a_i-n+1]$ is
%$\a_{ii}[-a_i-m,n]+\a_{ii}[-a_i-n+1,m-1]+\ldots$. Takes care of the case $a_i>1$. Otherwise need constant term.

%\subsection{Examples of $\wt{\UU}^{ns}_{g,n}(\ba_1,\ba_2,\ldots)$ and of Brill-Noether loci in $\wt{\UU}^{ns}_{g,n}(\ba)$}

In the next Lemma,
generalizing \cite[Lem.\ 2.3.3(i),(ii)]{P-ainf}, 
we show how certain Brill-Noether loci in $\wt{\UU}^{ns}_{g,n}(\ba)$ are cut out by vanishing conditions for 
appropriate coordinates $\a_{ij}[p,q]$ with $-a_j\le q\le -1$ (which are among the coordinates \eqref{generators-wide-eq}).
%Here are some more examples of calculating the Brill-Noether loci using coordinates for $\ba=(1,\ldots,1)$.
%Let us use the following special notation for some of the coordinates on $\wt{\UU}^{ns}_{g,g}(1,\ldots,1)$:
%$$\a_{ij}:=\a_{ij}[-2,-1], \ \ \b_{ij}=\a_{ij}[-3,-1].$$

\begin{lem}\label{BN-ex-lem}
% (i) The locus $h^0(2p_1+\ldots+2p_a+p_{a+1}+\ldots+p_{a+b})\ge 2$ for $2a+b=g$ is given by the equation
Let $D=\sum_i a_ip_i$. Fix $j$ such that $a_j>0$ and $i$ such that $i\neq j$ (so $n>1$). 
Then for $m,m'\ge 1$, such that $m'\le a_j$, the locus in $\wt{\UU}^{ns}_{g,n}(\ba)$ given
by the condition $h^0(D+mp_i-m'p_j)\ge m+1$ (equivalently, $=m+1$) for $1\le m'\le m$, is cut out by the equations
$$\a_{ij}[p,q]=0 \ \text{ for } -a_i-m\le p\le -a_i-1, -a_j\le q\le -a_j+m'-1.$$
\end{lem}

\Pf . Since $h^0(D+mp_i)=m+1$, the equivalent condition describing this locus is 
$$H^0(D+mp_i-m'p_j)=H^0(D+mp_i).$$
In other words, every element of $H^0(D+mp_i-m'p_j)$ has to have a pole of order $\le a_j-m'$ at $p_j$.
Now the equations are obtained by writing this condition on poles for the basis $(f_i[p])_{-a_i-m\le p\le -a_i-1}$
of $H^0(D+mp_i)/\lan 1\ran$.
\ed

%Each BN-locus $h^1(n_1p_1+...+n_gp_g)\ge 0$, where one of $n_i$'s is zero, gets mapped by the Krichever map to
%the intersection with a linear subspace given by vanishing
%of some Pl\"ucker coordinates.

The next result (which is an analog of \cite[Lem.\ 2.3.3(iii)]{P-ainf}) will be useful in the analysis of GIT stability conditions
for the $\G_m^n$-action on $\wt{\UU}^{ns}_{g,n}(\ba)$ (see Section \ref{GIT-sec}).

\begin{lem}\label{smooth-nonvanishing-lem}
Let $(C,p_1,\ldots,p_n)\in \wt{\UU}^{ns}_{g,n}(\ba)$ be such that $C$ is smooth. Then for each $i\in [1,n]$, such that 
$a_i<g$, there
exists $j\in [1,n]$, $m\ge 1$ and $q\in [-a_j,-1]$ such that $\a_{ij}[-a_i-m,q]\neq 0$. If $a_i=g$ then there exists
$m\ge 1$ and $q\ge -g$, $q\neq 0$, such that $\a_{ii}[-a_i-m,q]\neq 0$.
\end{lem}

\Pf . Otherwise, we would have $h^0((a_i+m)p_i)=h^0((a_i+m-1)p_i)+1$ for every $m\ge 1$. In other words,
$h^1(a_ip_i)=h^1((a_i+1)p_i)=\ldots=0$. But this is possible only if $a_i\ge g$. 
In this case $a_i=g$ (and $a_j=0$ for $j\neq i$). Now let us consider 
the image $(C,p_i)$ of $(C,p_1,\ldots,p_n)$ under the forgetting map $\wt{\UU}^{ns}_{g,n}(\ba)\to\wt{\UU}^{ns}_{g,1}(g)$.
The point of  $\wt{\UU}^{ns}_{g,1}(g)$ where all coordinates vanish corresponds to the singular curve $C^{\cusp}(\ba)$
(see Example \ref{special-curves-ex}.1), hence, there exists a nonvanishing coordinate $\a_{ii}[-a_i-m,q]$.
\ed

\subsection{GIT quotients of $\wt{\UU}^{ns}_{g,n}(\ba)$}\label{GIT-sec}

In this section we work over an algebraically closed field $k$ of characteristic zero.

%??? Recall that a {\it modular birational model} of $M_{g,n}$ is a proper birational model $X$ such that there exists
%an open substack $\XX\sub\UU_{g,n}$ in the stack of pointed curves of genus $g$, such that the map
%$\XX\to X$ is the coarse moduli space. If instead the map $\XX\to X$ is only required to be the {\it good moduli space}
%(i.e., universal for maps to algebraic spaces and bijective on {\it closed} $k$-points), then $X$ is called a {\it weakly
%modular birational model} of $M_{g,n}$.

We can view any element $\chi\in\Z^n$ as a character of $\G_m^n$. The corresponding GIT quotient of
$\wt{\UU}^{ns}_{g,n}(\ba)$ is
$$\wt{\UU}^{ns}_{g,n}(\ba)\sslash_\chi \G_m^n:=\Proj \left(\bigoplus_{p\ge 0} 
H^0(\wt{\UU}^{ns}_{g,n}(\ba),\OO)_{\chi^p}\right),$$
where the subscript $\chi^p$ denotes the subset of functions $f$ such that $(\la^{-1})^*f=\chi(\la)^pf$
for $\la\in\G_m^n$.

As in \cite[Sec.\ 2.4]{P-ainf}, we are going to show that many of the GIT quotients of our affine moduli schemes
$\wt{\UU}^{ns}_{g,n}(\ba)$ by the natural $\G_m^n$-action provide birational models of $M_{g,n}$.
Furthermore, we will describe one chamber of characters $\chi$ where every point of $\wt{\UU}^{ns}_{g,n}(\ba)$ with a smooth curve $C$ is $\chi$-stable, and show that
the corresponding GIT quotient is the coarse moduli space of a modular compactification
of $\MM_{g,n}$ in the sense of \cite{Smyth} (see Corollary \ref{GIT-cor} below).

For a full-dimensional cone $\bC$ we denote by $\int(\bC)$ its interior.

\begin{thm}\label{GIT-thm}
(i) For any weight $\chi\in\Z^n$ the GIT quotient $\wt{\UU}^{ns}_{g,n}(\ba)\sslash_\chi \G_m^n$
is a projective scheme over $k$. For $\chi\not\in\bC_\ba$ this GIT quotient is empty.
Assume that $\chi\in\bC_\ba$ (resp., $\chi\in\int(\bC_\ba)$). 
Then for any smooth curve $C$ of genus $g$ and for generic points $p_1,\ldots,p_n$ the
point $(C,p_1,\ldots,p_n,v_1,\ldots,v_n)\in\wt{\UU}^{ns}_{g,n}(\ba)$ is 
$\chi$-semistable (resp., $\chi$-stable). 

\noindent
(ii) Let $\bC_0\sub \bC_\ba$ be the subcone generated by all the vectors $\be_i$.
Assume that $\chi\in\int(\bC_0)$. Then any $(C,p_1,\ldots,p_n,v_1,\ldots,v_n)\in \wt{\UU}^{ns}_{g,n}(\ba)$,
with $C$ smooth, is $\chi$-stable. Furthermore, 
for such $\chi$ every $\chi$-semistable point is $\chi$-stable and the notion of semistability does not depend on $\chi$, and
hence the GIT quotient $\wt{\UU}^{ns}_{g,n}(\ba)\sslash_\chi \G_m^n$
does not depend on $\chi$.

\noindent
(iii) For any $\ba,\ba'\in X(g,n)$, the geometric quotient
$\wt{\UU}^{ns}_{g,n}(\ba,\ba')/\G_m^n$ is isomorphic to a distinguished open affine subset of
$\wt{\UU}^{ns}_{g,n}(\ba)\sslash_\chi \G_m^n$ for 
$$\chi=\chi_{\ba,\ba'}=\sum_{i=1}^n\left({a'_i+1\choose 2}-{a_i+1\choose 2}\right)\be_i\in \bC_\ba\cap\Z^n.$$
\end{thm}

\Pf . (i) The argument is similar that of \cite[Prop.\ 2.4.1(i)]{P-ainf}. Lemma \ref{global-functions-lem}
implies projectivity of all the GIT quotients. The fact that the weights of all the coordinate functions
$\a_{ij}[p,q]$ belong to $\bC_\ba$ (see Lemma \ref{cone-generation-lem}) implies that
the GIT quotients are empty for $\chi\not\in\bC_\ba$.

%Also, the map from the $\chi$-semistable locus to the GIT quotient
%is known to be a good moduli space in the sense of Alper \cite{Alper}. 

Let $\Om'\sub\Om_\ba$ be the set of weights $\om_{ij}$, where $i\neq j$ and $a_j>0$, together with the weight
$\be_i$ for $a_i=g$ (when such $i$ exists).
By Lemma \ref{cone-generation-lem}, the weights from $\Om'$ generate $\bC_\ba$.
%In particular, they generate a lattice of maximal rank. Hence, every $\chi\in\bC_\ba$ can be written as
%$$\chi=\sum_{\om\in\Om'} x_\om \om,$$
%with $x_\om\in\Q_{\ge 0}$. In the case when $\chi$ is in the interior of $\bC_\ba$ we can further assume that
%all $x_\om$ are positive. Let $M>0$ be such that $M x_\om\in\Z$ for all $\om\in\Om'$.

Thus, by \cite[Lem.\ 2.4.1(i)(ii)]{P-ainf}, it
is enough to prove that for every $\om\in\Om'$ there exists
a function $f_\om$ of weight $\om$, which does not vanish at a generic smooth curve.
In the case when all $a_i<g$ we can take $f_{\om_{ij}}$ to be the coordinate $\a_{ij}[-a_i-1,-a_j]$.
Indeed, since $a_j>0$, 
the condition $\a_{ij}[-a_i-1,-a_j]\neq 0$ is equivalent to $h^0(C,D+p_i-p_j)=1$, where $D=\sum_i a_ip_i$,
(see Lemma \ref{BN-ex-lem}), so it holds generically. In the case when there exists $i$ such that $a_i=g$ and all other $a_j=0$
we use the fact that there exists $\a_{ii}[p,q]$ with $q\neq 0$, which does not vanish at a smooth curve (see Lemma
\ref{smooth-nonvanishing-lem}).

%then we will have that a generic smooth curve belongs to the open set $\prod f_\om^{Mx_\om}\neq 0$,
%and so is $\chi$-semistable. Furhtermore, if all $x_\om>0$ then $\G_m^n$ acts on this open set with finite stabilizers,
%so it consists of $\chi$-stable points.

\noindent
(ii) Again the argument is similar to that of \cite[Prop.\ 2.4.1(ii)]{P-ainf}. Let $(C,p_1,\ldots,p_n,v_1,\ldots,v_n)\in\wt{\UU}^{ns}_{g,n}(\ba)$
with $C$ smooth. By Lemma \ref{smooth-nonvanishing-lem}, for each $i\in [1,n]$ there exists a 
coordinate of weight $\xi_i\in\Z^n$, which is either a positive multiple of $\be_i$ or has form $p\be_i+q\be_j$ for some $i\neq j$,
where $p\ge a_i+1$, $ -a_j\le q\le -1$. Let us rescale the basis $\be_i$ as follows:
$$\be'_i=\frac{1}{w_i}\be_i=\begin{cases} a_i\be_i, & a_i>0, \\ \frac{N}{N+1}\be_i, & a_i=0.\end{cases}$$
Then in terms of the new basis we have that $\xi_i$ is either a positive multiple of $\be'_i$ or
has form $\xi_i=x\be'_i+y\be'_j$, where $x>1$, $-1\le y<0$.
By \cite[Lem.\ 2.4.3]{P-ainf}, this implies that the cone $\bC_\xi\sub\R^n$ generated by $\xi_1,\ldots,\xi_n$ contains
$\bC_0$. Now as in \cite[Prop.\ 2.4.2(ii)]{P-ainf} we conclude that for any $\chi\in\int(\bC_0)\sub\int(\bC_\xi)$ the point $(C,p_1,\ldots,p_n,v_1,\ldots,v_n)$ is $\chi$-stable.

Now let us prove that for $\chi$ in the interior of $\bC_0$ all $\chi$-semistable points are $\chi$-stable, and that
this notion does not depend on $\chi$. To this end we will use the generating coordinates of
Proposition \ref{generators-prop}.
%such that their weights are either positive multiples of $\be_i$ for some $i$, or have form $p\be_i+q\be_j$ for some $i\neq j$,
%$p\ge a_i+1$, $-a_j\le q\le -1$. 
Suppose we have a point $(C,p_\bullet,v_\bullet)\in\wt{\UU}^{ns}_{g,n}(\ba)$ which is
$\chi$-semistable for some $\chi$ in the interior of $\bC_0$. Then for each $i$ there should exist a generating coordinate
of weight $\xi_i$ that has a positive $i$th component. But such $\xi_i$ is either a positive multiple of $\be_i$ or has
form $p\be_i+q\be_j$ with $i\neq j$, $p\ge a_i+1$, $-a_j\le q\le -1$. By \cite[Lem.\ 2.4.3]{P-ainf}, this implies that
$\bC_0$ is contained in the cone generated by $\xi_1,\ldots,\xi_n$. As before, we deduce that
the point $(C,p_\bullet,v_\bullet)$ is stable with respect to any character in $\int(\bC_0)$. 

\noindent (iii) Recall that $\wt{\UU}^{ns}_{g,n}(\ba,\ba')$ is the distinguished affine open subset in $\wt{\UU}^{ns}_{g,n}(\ba)$
associated with the pull-back of the determinant of the morphism of the vector bundles
$$\pi_{\ba}^{-1}\circ\pi_{\ba'}:\HH(\ba')/\HH_{\ge 0}\ot\OO\to \HH_{\ge -\ba}/\HH_{\ge 0}\ot \OO$$
under the Krichever map (see Proposition \ref{Plucker-prop}(ii) and the proof of Theorems A and B).
Now the assertion follows immediately from the fact that $\chi_{\ba,\ba'}$ is the weight of 
$\det(\pi_{\ba}^{-1}\circ\pi_{\ba'})$.
\ed

Recall that a {\it modular compactification of } $\MM_{g,n}$ (over $k$) is an nonempty open substack $\XX\sub\NN_{g,n}$,
where $\NN_{g,n}$ is the stack of smoothable curves of arithmetic genus $g$ with $n$ smooth distinct
marked points, such that $\XX$ is proper. This term was introduced by Smyth \cite{Smyth} who considered substacks
defined over $\Z$, whereas we work over $k$.

Let $\wt{\UU}^{ns}_{g,n}(\ba)^{ss}_\chi\sub \wt{\UU}^{ns}_{g,n}(\ba)$ denote the set of $\chi$-semistable points.
As in \cite[Cor.\ 2.4.4]{P-ainf}), we deduce the following result.

\begin{cor}\label{GIT-cor} 
For $\chi\in\int(\bC_0)$ the quotient stack $\wt{\UU}^{ns}_{g,n}(\ba)^{ss}_\chi/\G_m^n$ is proper.
Hence, its irreducible component consisting of smoothable curves gives a modular compactification of $\MM_{g,n}$.
\end{cor}

%\Pf . By Theorem \ref{GIT-thm}, for such $\chi$ all semistable points are stable, hence, the map from
%$\XX=\wt{\UU}^{ns}_{g,n}(\ba)^{ss}_\chi/\G_m^n$
%to the GIT quotient
%is a coarse moduli map. Since, the GIT quotient is projective, this implies that $\XX$ is proper
%\ed

\begin{rems}\label{GIT-rems} 
1. As in \cite[Sec.\ 2.4]{P-ainf} we can use the weights of generating coordinates
\eqref{generators-wide-eq} to define a countable set of walls in $\bC_\ba$, such that the statement of Corollary
\ref{GIT-cor} also holds for characters that do not lie on these walls. In fact, this assertion holds for some finite subset
of these walls since we know that finitely many of the coordinates \eqref{generators-wide-eq} generate the ring
$\OO(\wt{\UU}^{ns}_{g,n}(\ba))$ (by Theorem A).

\noindent
2. For a pair $\ba\neq \ba'\in X(g,n)$ one often has nontrivial global $\G_m^n$-invariant functions on 
$\wt{\UU}^{ns}_{g,n}(\ba,\ba')$. For example, for $g=2$, $n=2$, $\ba=(2,0)$ and $\ba'=(0,2)$, we
have $\chi_{\ba,\ba'}=3(\be_2-\be_1)$ which is in the interior of the cone $\bC_{2,0}$, generated by $\be_1$ and
$\be_2-2\be_1$. Hence, by Theorem \ref{GIT-thm},
the geometric quotient $\wt{\UU}^{ns}_{2,2}(2\be_1,2\be_2)/\G_m^2$ is a (nonempty)
distinguished open affine in a projective birational model of $M_{2,2}$.

However, it may happen that the weight $\chi_{\ba,\ba'}$ belongs to the
boundary of $\bC_\ba$. For example, for $g=2$, $n=2$, $\ba=(1,1)$ and $\ba'=(2,0)$ we have
$\chi=2\be_1-\be_2$ which is one of the generators of $\bC_{(1,1)}$. In this case one can check
that there are no non-constant $\G_m^2$-invariant functions on $\wt{\UU}^{ns}_{2,2}((1,1),(2,0))$. 
%(e.g., using Corollary \ref{g-2-n-2-intersection-cor} below).
%in the case $g=1$ it is easy to see that
%$$H^0(\wt{\UU}^{ns}_{1,n}(e_1,\ldots,e_n),\OO)^{\G_m^n}=k,$$
%where $e_1=(1,0,\ldots,0$, etc. Furthermore, the corresponding stack
%\wt{\UU}^{ns}_{1,n}(e_1,\ldots,e_n
\end{rems}

\subsection{Extending the coordinates on $\wt{\UU}^{ns}_{g,n}(\ba)$ to $\ov{\MM}_{g,n}$}\label{extending-sec}

%As in the previous section we use the identification of $\wt{\UU}^{ns}_{g,n}(\ba)$ with
%Fix $\ba\in X(g,n)$.

In this section we will use a different normalization of the coordinates $\a_{ij}[p,0]$ than in
Section \ref{generators-sec} (see Theorem \ref{pole-thm} below).

Recall that the standard divisor classes $\psi_i$ on the moduli stack of pointed curves
correspond to the line bundles $L_i:=\OO(-p_i)|_{p_i}$ on $\UU_{g,n}$. The pull-back of these line bundles
to $\wt{\UU}^{ns}_{g,n}(\ba)$ is canonically trivialized. Now considering the semi-invariance
of $\a_{ij}[p,q]$ with respect to the $\G_m^n$-action we see that they descend to sections
\begin{equation}\label{ov-a-ij-pq-eq}
\ov{\a}_{ij}[p,q]\in H^0(\UU^{ns}_{g,n}(\ba),L_i^{-p}\ot L_j^q).
\end{equation}
Let $\ov{\MM}_{g,n}(\ba)$ denote the open substack of $\ov{\MM}_{g,n}$ consisting of 
$(C,p_1,\ldots,p_n)$ such that $H^1(C,\OO(p_1+\ldots+p_n))=0$.
We have a natural morphism 
$$\bc:\ov{\MM}_{g,n}(\ba)\to \UU^{ns}_{g,n}(\ba)$$
obtained by restricting the map $\bc$ (see \eqref{bc-map-eq}) to stable curves.
%Note that it sends ???
We have $\bc^*L_i=L_i$, and we denote still by $\ov{\a}_{ij}[p,q]$ the pull-back of $\ov{\a}_{ij}[p,q]$ to
$\ov{\MM}_{g,n}(\ba)$.

We are interested in the poles of $\ov{\a}_{ij}[p,q]$
along the complement $\ov{\MM}_{g,n}\setminus \ov{\MM}_{g,n}(\ba)$. In Theorem \ref{pole-thm} 
below we will achieve an estimate for these poles which ignores a possibly more subtle behavior
along the boundary of $\ov{\MM}_{g,n}$.

Let $(C,p_1,\ldots,p_n)$ be the universal curve over $\ov{\MM}_{g,n}$.
For a segment $[a,b]\sub\Z$ let us define the vector bundle on $\ov{\MM}_{g,n}$ by
$$\Polar_{[a,b],p_i}:=\pi_*(\OO_C(-ap_i)/\OO_C((-b-1)p_i)).$$
Let us set 
$$\Polar_\ba:=\bigoplus_{i=1}^n\Polar_{[-a_i,-1],p_i}=\bigoplus_{i=1}^n\pi_*(\OO_C(a_ip_i)/\OO_C).$$
As before, we denote by $\sV$ the bundle with the fibers $H^1(C,\OO)$ (the dual of the Hodge bundle).
We have a natural morphism of vector bundles of rank $g$,
\begin{equation}\label{ov-pi-a-eq}
\ov{\pi}_\ba:\Polar_\ba \to \sV,
\end{equation}
such that $\ov{\MM}_{g,n}(\ba)$ is precisely the locus
where $\ov{\pi}_\ba$ is an isomorphism. 

\begin{defi}\label{Z-ba-defi} 
Let us consider the section 
$$\ov{s}_{\ba}:=\det(\ov{\pi}_\ba)\in H^0(\det(\sV)\ot \prod \det(\Polar_\ba)^{-1}),$$
We define the effective Cartier divisor $Z_\ba\sub \ov{\MM}_{g,n}$ 
as the zero locus of $\ov{s}_\ba$. 
Note that $\ov{\MM}_{g,n}(\ba)=\ov{\MM}_{g,n}\setminus Z_\ba$.
\end{defi}

\begin{rem}
Let $\ov{\MM}^{(\infty)}_{g,n}\to \ov{\MM}_{g,n}$ be  the torsor
corresponding to choices of formal parameters at marked points.
Note that the pull-back of $\ov{\pi}_\ba$ 
agrees with the map $\Kr^*(\pi_\ba)$ over the open substack of $\ov{\MM}^{(\infty)}_{g,n}$
where the Krichever map is defined, where $\pi_\ba$ is given by \eqref{pi-ba-eq}.
Hence, $\ov{s}_\ba$ similarly agrees with $\Kr^*(s_\ba)$, where $s_\ba$ is the Pl\"ucker coordinate corresponding to $\ba$
(see Section \ref{cell-sec}).
\end{rem}

%Recall also that we defined effective divisors $Z_\ba$ on $SG_1(g)$ given as
%degeneration loci of morphisms of vector bundles 
%$$\pi_\ba:\HH_{\ba}/\HH_{\ge 0}\ot \OO\to \sV.$$

Below we use the standard divisor classes on $\ov{\MM}_{g,n}$, $\psi_i=c_1(L_i)$ and $\la=-c_1(\sV)$.

\begin{lem} One has the equality of divisor classes on $\ov{\MM}_{g,n}$
$$Z_\ba:=\Psi_\ba-\la,$$
where
\begin{equation}\label{Psi-eq}
\Psi_\ba:=\sum_{i=1}^n{a_i+1 \choose 2}\psi_i.
\end{equation}
\end{lem}

\Pf . Since $Z_\ba$ is the degeneration locus of the morphism $\ov{\pi}_\ba$, its class is given by
$$c_1(\sV)-c_1(\Polar_\ba).$$
It remains to use the isomorphism 
$$\det(\Polar_{[-m,-1],p_i})\simeq \bigotimes_{a=1}^m L_i^{-a}\simeq L_i^{-{m+1\choose 2}}$$
and recall that $c_1(\sV)=-\la$, $c_1(L_i)=\psi_i$.
\ed
%Then by ???, we have a Krichever morphism
%$$\Kr:\ov{M}^{(\infty)}_{g,n}\to SG_1(g).$$
%The class $Z_\ba$ is represented by an effective Cartier divisor on $\ov{M}_{g,n}$, such that
%its pull-back to $\ov{M}^{(\infty)}_{g,n}$ is equal to $\Kr^*Z_\ba$

\begin{rem}
Note that the effective divisor $Z_\ba\cap \MM_{g,n}$ is precisely the locus of $(C,p_1,\ldots,p_n)$
such that $h^1(\OO(p_1+\ldots+p_n))\neq 0$, or equivalently, $h^0(\OO(p_1+\ldots+p_n))\ge 2$.
Let $D_\ba$ be the closure of this locus in $\ov{\MM}_{g,n}$.
The divisors $(D_\ba)$ play an important role in the work of Logan \cite{Logan} on Kodaira dimension
of $\ov{\MM}_{g,n}$. In particular, he proves the following formula in $\Pic(\ov{\MM}_{g,n})$:
$$D_\ba=\Psi_\ba-\la-\sum m_{i,S}(\ba)\de_{i,S},$$
where $(\de_{i,S})$ are the boundary divisors on $\ov{\MM}_{g,n}$, and computes explicitly $m_{i,S}\ge 0$ 
(for example, he shows that $m_{0,\{i,j\}}=a_ia_j$).
On the other hand, we defined $Z_\ba=\Psi_\ba-\la$ as an effective divisor, given by zeros of the section $s_\ba$,
which vanishes with multiplicity $1$ along $D_\ba$. Since the boundary divisors are linearly independent
in $\Pic(\ov{\MM}_{g,n})$, it follows that $m_{i,S}(\ba)$ are precisely the multiplicities of zeros of $s_\ba$
along the boundary divisors. 
\end{rem}
%$$Z_\ba=D_\ba+\sum m_{i,S}(\ba)\de_{i,S}$$

For any $N\ge 1$ let us set 
$$\Polar_{\ba,N}:=\bigoplus_{i=1}^n\Polar_{[-a_i,N-1]}=\bigoplus_{i=1}^n\pi_*(\OO_C(a_ip_i)/\OO_C(-Np_i)).$$
Similarly to Section \ref{cell-sec}, it will be convenient to consider the natural morphism, induced by the coboundary homomorphism,
$$\ov{\pi}_{\ba,N}:\Polar_{\ba,N}/\OO \to \sV_N,$$
where $\sV_N$ is the vector bundle on $\ov{\MM}_{g,n}$ with the fiber $H^1(C,\OO(-N(p_1+\ldots+p_n)))$.
Note that the pull-back of this morphism to $\ov{\MM}^{(\infty)}_{g,n}$ agrees with $\Kr^*(\pi'_{\ba,N})$
on the open substack where $\Kr$ is defined (see the proof of Proposition \ref{Plucker-prop}).

\begin{lem}\label{pi-a-N-det-lem} 
We have a natural identification
$$\det(\sV_N)\ot\det(\Polar_{\ba,N})^{-1}\simeq \det(\sV)\ot\det(\Polar_\ba)^{-1},$$
so that $\det(\ov{\pi}_{\ba,N})=\det(\ov{\pi}_\ba)$.
\end{lem}

\Pf . This follows immediately 
from the morphism of exact sequences (similar to the one used in the proof of Proposition \ref{Plucker-prop})
\begin{diagram}
0\rTo{} &\Polar_{0,N}/\OO &\rTo{}&\Polar_{\ba,N}/\OO &\rTo{}& \Polar_{\ba} &
\rTo{} 0\\
&\dTo{\id }&&\dTo{\ov{\pi}_{\ba,N}}&&\dTo{\ov{\pi}_{\ba}}\\
0\rTo{} & \Polar_{0,N}/\OO &\rTo{}&\sV_{N}&\rTo{}&\sV&\rTo{} 0
\end{diagram}
\ed

We are going to work over $\UU^{ns}_{g,n}(\ba)$ for a while.
We consider the bundles $\Polar_{[a,b],p_i}$, $\Polar_{\ba,N}$ and the morphisms $\ov{\pi}_\ba$, $\ov{\pi}_{\ba,N}$
over $\UU^{ns}_{g,n}$, defined (and denoted) in the same way as over $\ov{\MM}_{g,n}$.

Recall that we have canonical formal parameters $t_i$ at the marked points $p_i$ on the universal curve
over $\wt{\UU}^{ns}_{g,n}(\ba)\simeq ASG\cap\Si^{\ba,i_0}$ (see Corollary \ref{curves-sec-cor}),
depending on a choice of $i_0$ such that $a_{i_0}>0$ in the case when some $a_i$ are zero.
These formal parameters are uniquely characterized by the following properties. If $a_i>0$ then for each
$m>a_i$ there is an element $f_i[-m]\in H^0(C,\OO(mp_i+\sum_{j\neq i}a_jp_j))$ with the expansion at $p_i$ satisfying
$f_i[-m]\equiv t_i^{-m}\mod t_i^{-a_i+1}k[[t_i]]$. If $a_i=0$ then there should exist an element 
$f_i[-1]\in H^0(C,\OO(p_i+\sum_{j\neq i}a_jp_j))$ which restricts to $t_i^{-1}$ in a formal neighborhood of $p_i$
and whose expansion in $t_{i_0}$ at $p_{i_0}$ has no constant term. 

For any $a\le b$, the splitting of the trivial vector bundle 
$(t_i^ak[t_i]/t_i^{b+1}k[t_i])\ot \OO$ associated with the basis $(t_i^j)$, over 
$\wt{\UU}^{ns}_{g,n}(\ba)$, descends to a splitting
\begin{equation}\label{W-splitting-eq}
\Polar_{[a,b],p_i}\simeq \bigoplus_{j=a}^b L_i^j
\end{equation}
over $\UU^{ns}_{g,n}(\ba)$.

Now we will give a recipe for computing
the sections $\ov{\a}_{ij}[p,q]$ of $L_i^{-p}L_j^q$, where $p\le -a_i-1$, $q\ge -a_j$ (see \eqref{ov-a-ij-pq-eq}),
on $\UU^{ns}_{g,n}(\ba)$, and hence on $\ov{\MM}_{g,n}(\ba)$, 
in terms of the above splitting. 

Let us set 
$$\Polar^i_{\ba,N}:=\bigoplus_{i'\neq i}^n\Polar_{[-a_{i'},N-1],p_{i'}}, \ \ \Polar^i_{\ba}:=\Polar^i_{\ba,0}.$$
Note that for each $p\le -a_i-1$
the isomorphism $\ov{\pi}_\ba$ over $\UU^{ns}_{g,n}(\ba)$ (see \eqref{ov-pi-a-eq}) induces an isomorphism
$$\Polar_{[p,-1],p_i}\oplus\Polar^i_\ba\rTo{\sim} \Polar_{[p,-a_i-1],p_i}\oplus\sV.$$
Hence, its inverse gives a morphism
$$F^i_p: \Polar_{[p,-a_i-1],p_i}\to \Polar_{[p,-1],p_i}\oplus\Polar^i_\ba$$
on $\UU^{ns}_{g,n}(\ba)$.
Similarly, for $N\ge 1$ the inverse of the isomorphism
$$\left(\Polar_{[p,N-1],p_i}\oplus\Polar^i_{\ba,N}\right)/\OO\rTo{\sim} \Polar_{[p,-a_i-1],p_i}\oplus\sV.$$
induces a morphism
$$F^i_{p,N}: \Polar_{[p,-a_i-1],p_i}\to \left(\Polar_{[p,N-1],p_i}\oplus\Polar^i_{\ba,N}\right)/\OO.$$

\begin{lem}\label{a-ij-composition-lem} 
In the case $q<0$ the section $\ov{\a}_{ij}[p,q]$ on $\UU^{ns}_{g,n}(\ba)$ is given by the composition
\begin{equation}\label{alpha-ij-composition-eq}
L_i^p\to \Polar_{[p,-a_i-1],p_i}\rTo{F^i_p}\Polar_{[p,-1],p_i}\oplus \Polar^i_\ba\to \Polar_{[-a_j,-1],p_j}\to L_j^q,
\end{equation}
where the first and the last arrows use the splitting \eqref{W-splitting-eq}.
If $q>0$ then $\ov{\a}_{ij}[p,q]$ is given by the composition
\begin{equation}\label{alpha-ij-composition-positive-q-eq}
L_i^p\to \Polar_{[p,-a_i-1],p_i}\rTo{F^i_{p,q+1}} \left(\Polar_{[p,q],p_i}\oplus\Polar^i_{\ba,q+1}\right)/\OO\to
\Polar_{[-a_j,q],p_j}/\OO\to L_j^q.
\end{equation}
Finally, in the case $q=0$, let us normalize $(\ov{\a}_{ij}[p,0])$ by the condition $\ov{\a}_{i,j_0}[p,0]=0$ for some fixed
$j_0\in [1,n]$. Then $\ov{\a}_{ij}[p,0]$, for $j\neq j_0$, 
is given by the composition obtained from \eqref{alpha-ij-composition-positive-q-eq} by replacing the last two arrows with
$$\left(\Polar_{[p,0],p_i}\oplus\Polar^i_{\ba,1}\right)/\OO\to \left(\Polar_{[-a_{j_0},0],p_{j_0}}\oplus \Polar_{[-a_j,0],p_j}\right)/\OO\to \OO,$$
where the last arrow is given by the difference of projections to $\OO$ from both factors in the direct sum.
\end{lem}

\Pf . It is enough to prove the similar assertions for the functions $\a_{ij}[p,q]$ on $\wt{\UU}^{ns}_{g,n}(\ba)$.
These follow essentially from the definition of $\a_{ij}[p,q]$. Let us consider the case $q<0$ first.
We can think of 
$\Polar_{[p,-1],p_i}\oplus\Polar^i_\ba$ as the bundle of polar parts at the marked points 
(up to order $-p$ at $p_i$ and up to order $a_j$ at $p_j$, $j\neq i$). The kernel of the morphism to $\sV$
corresponds to polar parts coming from $H^0(C\setminus\{p_1,\ldots,p_n\},\OO)$.
Hence, the map $F^i_p$ applied to a given polar part at $p_i$ gives its extension to polar parts at all points,
that come from $H^0(C\setminus\{p_1,\ldots,p_n\},\OO)$. Thus, the composition of the first two
arrows in \eqref{alpha-ij-composition-eq} sends $t_i^p$ to the polar parts of $f_i[p]$ at all points. The
composition of the two following arrows sends it to the term containing $t_j^q$ in the polar part at $p_j$,
which is by definition $\a_{ij}[p,q]t_j^q$.

The case $q>0$ is similar: the difference is that rather than considering polar parts we have to take into
account terms with higher order of the parameters. In the case $q=0$ we make a necessary change
to take into account that $\a_{ij}[p,0]$ is the difference of constant terms at $p_j$ and $p_{j_0}$.
\ed

Given a surjective map of vector bundles $p:V\to Q$, and an
effective Cartier divisor $D$, we say that a map $Q\to V(D)$ is a {\it splitting of $p$ with a pole} if the composition
$Q\to V(D)\to Q(D)$ is the map induced by the embedding $\OO\to \OO(D)$. Dually we talk about {\it retractions with a pole}
for embeddings of vector bundles.

\begin{lem}\label{split-lem} Let
$$0\to V_1\to V_2\to V_3\to 0$$
be an exact sequence of vector bundles, and let $D$ be an effective divisor, such that the sequence splits on the complement
$U$ of $D$. Then there is a splitting with a pole $V_3\to V_2(D)$ extending the given splitting over $U$,
$V_3|_U\to V_2|_U$,
if and only if there is a retraction with a pole $V_2\to V_1(D)$ extending the given retraction over $U$, $V_2|_U\to V_1|_U$.
\end{lem}

\Pf . Assume there exists $V_3\to V_2(D)$. Consider the global morphism 
$$f:V_1\oplus V_3(-D)\to V_2.$$ 
Then it is easy to see that $f$ is injective and $\coker(f)\simeq V_3/V_3(-D)$. In particular, $\coker(f)$ is killed by the
ideal sheaf of $D$, so we get an embedding $V_2(-D)\to V_1\oplus V_3(-D)$, which gives the required map
$V_2(-D)\to V_1$.

The proof of the converse is similar by looking at the morphism $V_2\to V_1(D)\oplus V_3$.
\ed

%We say that a morphism of vector bundles of equal rank $f:V\to W$ {\it degenerates on} an effective Cartier
%divisor $D$
%if the morphism $\det(f):\det(V)\to\det(W)$ factors through a morphism $\det(V)(D)\to \det(W)$.
We say that a morphism of vector bundles of equal rank $f:V\to W$ {\it degenerates exactly on} an effective Cartier
divisor $D$ if the morphism $\det(f):\det(V)\to\det(W)$ induces an isomorphism $\det(V)(D)\simeq\det(W)$.
Note that in this case the sheaf $W/f(V)$ is killed by the ideal sheaf of $D$, hence, we get a morphism $W(-D)\to V$
which is inverse to $f$ on the complement of $D$. 
%We say that $f$ {\it degenerates exactly on} $D$, if in addition the map $\det(V)(D)\to\det(W)$ is an isomorphism.

In the following lemma we analyze how to extend the splittings $L_i^a\to\Polar_{[a,b],p_i}$ of the natural projection
$\Polar_{[a,b],p_i}\to L_i^a$ (resp., retractions $\Polar_{[a,b],p_i}\to L_i^b$ of the embedding
$L_i^b\to\Polar_{[a,b],p_i}$) that we have over $\ov{\MM}_{g,n}(\ba)$ (see \eqref{W-splitting-eq}) to splittings 
(resp., retractions) with a pole
along the divisor $Z_a\sub\ov{\MM}_{g,n}$ (see Definition \ref{Z-ba-defi}).

\begin{lem}\label{W-embedding-splitting-lem}
(i) Let $Z$ be a divisor supported on $Z_\ba$.
Assume that for some $i\in [1,n]$ and some $m\ge 0$, $N\ge 1$, the splittings $L_i^{-m-j}\to \Polar_{[-m-j,-m],p_i}$
over $\ov{\MM}_{g,n}(\ba)$,  for $j=1,\ldots,N$, extend to splittings with a pole
$L_i^{-m-j}\to \Polar_{[-m-j,-m],p_i}(jZ)$ over $\ov{\MM}_{g,n}$. Then
we have a splitting with a pole 
$L_i^a\to\Polar_{[a,a+N],p_i}(NZ)$ for any $a\in\Z$, extending the splitting over $\ov{\MM}_{g,n}(\ba)$.

\noindent
(ii) For every $i\in [1,n]$ such that $a_i\ge 1$, and every $a\le b$, the splitting $L_i^a\to\Polar_{[a,b],p_i}$ (resp., the
retraction $\Polar_{[a,b],p_i}\to L_i^b$) over $\ov{\MM}_{g,n}(\ba)$
extends to a splitting (resp., retraction) with a pole
$$L_i^a\to\Polar_{[a,b],p_i}((b-a)Z_\ba)
 \ \text{ (resp., } \Polar_{[a,b],p_i}\to L_i^b((b-a)Z_\ba) \ )$$
over $\ov{\MM}_{g,n}$.

%\noindent
%(iii) Assume still that $a_i\ge 1$. In the case $[a,b]=[-a_i-m,-1]$, where $m\ge 1$, the embedding in (i) can be improved to
%$$L_i^{-a_i-m}\to \WW_{[-a_i-m,-1]}(p_i)(mZ_\ba).$$
%MAYBE DON'T NEED???

\noindent
(iii) In the case $a_i=0$ the above splittings and retractions acquire at most the following poles along $Z_\ba$:
$$L_i^a\to\Polar_{[a,b],p_i}((b-a)(a_{i_0}+1)Z_\ba),
 \ \ \Polar_{[a,b],p_i}\to L_i^b((b-a)(a_{i_0}+1)Z_\ba).$$
\end{lem}

\Pf . (i) Let us temporarily denote by $t_{i,can}$ the canonical formal parameter at $p_i$ over $\wt{\UU}^{ns}_{g,n}(\ba)$
obtained from Corollary \ref{curves-sec-cor}.
The question is local, so we can choose a formal parameter $t_i$ at $p_i$ locally over $\ov{\MM}_{g,n}$.
In particular, $t_i$ gives a nonzero relative tangent vector at $p_i$, and so we also have the corresponding
canonical parameter $t_{i,can}$ at $p_i$, which is defined away from
$Z_\ba$. We have
$$t_{i,can}=t_i+c_1t_i^2+\ldots+c_Nt_i^{N+1}+\ldots,$$
where $c_l$ are some functions, regular on the complement to $Z_\ba$.

We claim that $c_l\in\OO(lZ)$ for $l=1,\ldots,N$.
Using the induction on $N$ we can assume that we already know this for $c_l$ with $l<N$.
Now our assumption about the pole of the splitting $L_i^{-m-N}\to \Polar_{[-m-N,-m],p_i}$ means
that the coefficients of powers of $t_i$ in the expansion of $t_{i,can}^{-m-N}\mod t_i^{-m+1}k[[t_i]]$
belong to $\OO(NZ)$. Looking at the coefficient of $t_i^{-m}$ of this expansion we derive that $c_N\in \OO(NZ)$,
which proves our claim. 

Hence, for any $a\in\Z$ the coefficients in the expansion of $t_{i,can}^a\mod t_i^{a+N+1}k[[t_i]]$
belong to $\OO(NZ)$, which is precisely our assertion.

\noindent
(ii) Since, $Z_\ba$ is the degeneration locus of the morphism $\ov{\pi}_\ba$ (see \eqref{ov-pi-a-eq}), it follows
 that for any $m\ge 1$ the morphism of vector bundles of rank $g+m$ on $\ov{\MM}_{g,n}$,
$$\Polar_{[-a_i-m,-1],p_i}\oplus \Polar^i_\ba\to \Polar_{[-a_i-m,-a_i-1],p_i}\oplus \sV$$
degenerates exactly on $Z_\ba$.
Thus, we get a morphism 
\begin{equation}\label{a_i+m-morphism-eq}
\Polar_{[-a_i-m,-a_i-1],p_i}\to \Polar_{[-a_i-m,-1],p_i}(Z_\ba)\to \Polar_{[-a_i-m,-a_i],p_i}(Z_\ba).
\end{equation}
As in the proof of Lemma \ref{a-ij-composition-lem}, we see that over $\ov{\MM}_{g,n}(\ba)$  the first arrow in
\eqref{a_i+m-morphism-eq} sends $t_i^{-a_i-m'}$, for $1\le m'\le m$,
to the polar part the expansion of $f_i[-a_i-m']$ in $t_i$. By the definition of the canonical parameters, this
implies that the morphism \eqref{a_i+m-morphism-eq} is compatible with the splitting \eqref{W-splitting-eq} over
$\ov{\MM}_{g,n}(\ba)$.

Now, let us prove by induction that for any $m\ge 1$ we have the required splittings with poles
$$L_i^{a}\to\Polar_{[a,a+m],p_i}(mZ_\ba)$$
%, \ \ L_i^{-a_i-m}\to \WW_{[-a_i-m,-a_i]}(p_i)(mZ_\ba).$$
For $m=1$ the morphism \eqref{a_i+m-morphism-eq} gives a splitting with a pole
$$L_i^{-a_i-1}\to \Polar_{[-a_i-1,-a_i],p_i}(Z_\ba),$$
and the assertion follows from part (i).
Assume we already have such morphisms for $m'<m$. Then
we can construct the required splitting with a pole for $a=-a_i-m$ as the composition
$$L_i^{-a_i-m}\to \Polar_{[-a_i-m,-a_i-1],p_i}((m-1)Z_\ba)\to \Polar_{[-a_i-m,-a_i],p_i}(mZ_\ba),$$
where the first arrow exists by the induction assumption and the second arrow is given by \eqref{a_i+m-morphism-eq}.
It remains to apply part (i) again.

Next, combining our splittings with poles we get a splitting with a pole
$$\Polar_{[a,b-1],p_i}\to \Polar_{[a,b]}((b-a)Z_\ba).$$
Hence, by Lemma \ref{split-lem} we get the required retraction with a pole 
$$\Polar_{[a,b]}\to L_i^b((b-a)Z_\ba).$$

\noindent
(iii) Now consider the case $a_i=0$. Recall that we have fixed $i_0$ such that $a_{i_0}>0$.
It follows from Lemma \ref{pi-a-N-det-lem} that for each $m\ge 1$ the morphism of vector bundles of rank $g+n+m-1$
on $\ov{\MM}_{g,n}$,
$$\left(\Polar_{[-m,0],p_i}\oplus \Polar^i_{\ba,1}\right)/\OO\to \Polar_{[-m,-1],p_i}\oplus \sV_1$$
degenerates exactly on $Z_\ba$. Thus, we get a morphism
\begin{equation}\label{L-Z-a-i0-eq}
\Polar_{[-m,-1],p_i}(-Z_\ba)\to \left(\Polar_{[-m,0],p_i}\oplus \Polar_{[-a_{i_0},0],p_{i_0}}\right)/\OO.
\end{equation}

By part (ii), we also have a morphism $\Polar_{[-a_{i_0},0],p_{i_0}}\to \OO(a_{i_0}Z_\ba)$, 
which gives splitting with the pole of order $a_{i_0}$ along $Z_\ba$ of the exact sequence 
$$0\to \OO\to \Polar_{[-a_{i_0},0],p_{i_0}}\to \Polar_{[-a_{i_0},0],p_{i_0}}/\OO\to 0$$
Hence, the exact sequence
$$0\to \Polar_{[-m,0],p_i}\to (\Polar_{[-m,0],p_i}\oplus \Polar_{[-a_{i_0},0],p_{i_0}})/\OO\to 
\Polar_{[-a_{i_0},0],p_{i_0}}/\OO\to 0$$
also has a splitting with pole of order $a_{i_0}$ along $Z_\ba$, so we get a retraction with a pole
$$(\Polar_{[-m,0],p_i}\oplus \Polar_{[-a_{i_0},0],p_{i_0}})/\OO\to \Polar_{[-m,0],p_i}(a_{i_0}Z_\ba).$$
Composing it with \eqref{L-Z-a-i0-eq} we get a morphism
$$\Polar_{[-m,-1],p_i}\to \Polar_{[-m,0],p_i}((a_{i_0}+1)Z_\ba).$$
%$$L_i^{-1}\to \Polar_{[-1,0]}(p_i)((a_{i_0}+1)Z_\ba).$$
%By (i), this gives us a morphism $L_i^{a}\to \Polar_{[a,a+1]}(p_i)((a_{i_0}+1)Z_\ba)$ for any $a\in\Z$.
As in part (ii), we deduce from this by induction on $m$ that we have splittings with poles
$$L_i^{a}\to\Polar_{[a,a+m],p_i}(m(a_{i_0}+1)Z_\ba), \ \ L_i^{-m}\to \Polar_{[-m,0],p_i}(m(a_{i_0}+1)Z_\ba)$$
extending the ones associated with the canonical parameter $t_i$ over $\ov{\MM}_{g,n}(\ba)$. 
\ed

\begin{thm}\label{pole-thm} 
Let us choose $j_0\in [1,n]$ such that $a_{j_0}=\min(a_1,\ldots,a_n)$, and normalize $\ov{\a}_{ij}[p,0]$ by
$\ov{\a}_{i,j_0}[p,0]=0$. If there exists $i$ such that $a_i=0$ then we in addition choose $i_0$, such that $a_{i_0}>0$,
that is used to define $\ov{\a}_{i,j}[p,q]$ for $a_i=0$.
%\noindent
%(i) If $a_i>0$ and $a_j>0$ then 
Then the section $\ov{\a}_{ij}[p,q]$ of $L_i^{-p}\ot L_j^q$
on $\ov{\MM}_{g,n}(\ba)$, where $p\le -a_i-1$, $q\ge -a_j$, extends to an element of 
$$H^0(\ov{\MM}_{g,n},L_i^{-p}\ot L_j^q(1+d_j(q+a_j)-d_i(p+a_i+1))Z_\ba),$$ 
where
$$d_i=\begin{cases} 1 & a_i>0,\\ (a_{i_0}+1), & a_i=0.\end{cases}$$
In particular if $a_i>0$ and $a_j>0$ then $\ov{\a}_{ij}[p,q]$ extends to a global section of
$L_i^{-p}\ot L_j^q((q+a_j-p-a_i)Z_\ba)$.
\end{thm}

\Pf . We use the presentation of $\ov{\a}_{ij}[p,q]$ given in Lemma \ref{a-ij-composition-lem}.
%Note that because of our choice of $i_0$ Lemma \ref{W-embedding-splitting-lem} gives a uniform statement about
%the poles, independent on whether $a_i>0$ or $a_i=0$: we alway have morphisms
%$$L_i^a\to\WW_{[a,b]}(p_i)((b-a)Z_\ba), \ \ \WW_{[a,b]}(p_i)\to L_i^b((b-a)Z_\ba)$$
%on $\ov{\MM}_{g,n}$.
Assume first that $q<0$. Then using Lemma \ref{W-embedding-splitting-lem}, we can modify the composition 
\eqref{alpha-ij-composition-eq} as 
\begin{align*}
&L_i^p\left(d_i(p+a_i+1)Z_\ba\right)\to \Polar_{[p,-a_i-1],p_i}
\to\left(\Polar_{[p,-1],p_i}\oplus \Polar^i_\ba\right)(Z_\ba)\to \\
&\Polar_{[-a_j,-1],p_j}(Z_\ba)\to 
 \Polar_{[-a_j,q],p_j}(Z_\ba)\to L_j^q\left((1+d_j(q+a_j))Z_\ba\right),
\end{align*}
which gives the required pole estimate in this case.
Similarly, if $q>0$ then we modify \eqref{alpha-ij-composition-positive-q-eq} as
\begin{align*}
&L_i^p(d_i(p+a_i+1)Z_\ba)\to \Polar_{[p,-a_i-1],p_i}
\to\left(\left(\Polar_{[p,q],p_i}\oplus \Polar^i_{\ba,q+1}\right)/\OO\right)(Z_\ba)\to \\ 
& (\Polar_{[-a_j,q],p_j}/\OO)(Z_\ba)\to L_j^q\left((1+d_j(q+a_j))Z_\ba\right).
\end{align*}
Finally, in the case $q=0$ and $j\neq j_0$ we replace the last two arrows in this composition by
\begin{align*}
&\bigl(\left(\Polar_{[p,0],p_i}\oplus \Polar^i_{\ba,1}\right)/\OO\bigr)(Z_\ba)\to 
\bigl(\left(\Polar_{[-a_{j_0},0],p_{j_0}}\oplus \Polar_{[-a_j,0],p_j}\right)/\OO\bigr)(Z_\ba) \\
&\to L_j^q((1+d_ja_j)Z_\ba),
\end{align*}
where the last arrow exists since $a_{j_0}\le a_j$.
\ed

\begin{rem} In the case of $\ov{\a}_{ij}[-a_i-1,-a_j]$, where $a_j>0$, we can derive the result of the above theorem
much easier using the formula
$$\ov{\a}_{ij}[-a_i-1,-a_j]=\pm \frac{\ov{s}_{\ba+\be_i-\be_j}}{\ov{s}_\ba}$$
obtained similarly to Corollary \ref{a-simple-for-cor}. Indeed,  
the right hand side is a rational section of 
$\OO(\Psi_{\ba+\be_i-\be_j}-\Psi_\ba)=L_i^{a_i+1}\ot L_j^{-a_j}$ with the pole of order $1$ at $Z_\ba$, so we can view it as a regular
section of $L_i^{a_i+1}\ot L_j^{-a_j}(Z_\ba)$.
\end{rem}

%For $\ba=(a_1,\ldots,a_n)\in X(g,n)$.
In the next result we will use the notation from Section \ref{GIT-sec}. In particular, $N=\max(a_1,\ldots,a_n)$
and $\bC_\ba\sub \R^n$ is the cone generated
by all vectors $\om_{ij}=(a_i+1)\be_i-a_j\be_j$. We will also use the smaller subset of generators of this cone indicated
in Lemma \ref{cone-generation-lem}.
%Recall that in the case when $a_i>0$ for at least two indices $i$
%the cone $\bC_\ba$ is also generated
%by vectors $\om_{ij}$ with $i\neq j$ and $a_j>0$. It is also easy to see that

\begin{prop}\label{rational-classes-prop} 
(i) Assume that  $a_i>0$ for at least two indices $i$. Let $\chi=\sum_{i\neq j: a_j>0}x_{ij}\om_{ij}$ be in the cone $\bC_\ba$.
Then the rational map from $\ov{\MM}_{g,n}$ to $\wt{\UU}^{ns}_{g,n}(\ba)\sslash_\chi \G_m^n$ 
is given by a linear system in $|rZ|$, for some positive integer $r$, where 
$$Z=Z(\chi)=\sum_{i\neq j: a_i>0,a_j>0} x_{ij}\left(Z_{\ba+\be_i-\be_j}+(\frac{N}{a_i}-1)Z_\ba\right)+\sum_{i\neq j: a_i=0,a_j>0}
x_{ij}Z_{\ba+\be_i-\be_j}.$$

\noindent (ii) For $\ba=g\be_1$ and $\chi=x_1\be_1+\sum_{i=2}^n x_i\om_{i1}$
the rational map from $\ov{\MM}_{g,n}$ to $\wt{\UU}^{ns}_{g,n}(\ba)\sslash_{\chi}\G_m$
is given by a linear system in a positive multiple of the class
$$x_1\left(\psi_1+Z_{\be_1}\right)+\sum_{i=2}^n x_i Z_{(g-1)\be_1+\be_i}.$$
%$$\psi+Z_g=\frac{g^2+g+2}{2}\psi-\la.$$
%coordinate $\a_{ij}[p,q]$ descends to a section of ???
\end{prop}

\Pf . Set 
$$\wt{\psi}_i=\psi_i+N\ell(\be_i)Z_\ba=
\begin{cases} \psi_i+(\frac{N}{a_i})Z_\ba, & a_i\neq 0,\\ \psi_i+(N+1)Z_\ba, & a_i=0.\end{cases}$$
We claim that $\ov{\a}_{ij}[p,q]$ extends to a global section of the line bundle corresponding to the class
$-p\wt{\psi}_i+q\wt{\psi}_j$, where for a rational divisor class $\sum r_i [V_i]$, by global sections we mean global sections of
$\sum \lfloor r_i \rfloor [V_i]$. Indeed, 
$$-p\wt{\psi}_i+q\wt{\psi}_j=-p\psi_i+q\psi_j+N\ell(-p\be_i+q\be_j)Z_\ba.$$
But the linear functional $\ell$ has the property that $N\ell(\om_{ij})\ge 1$ and $N\ell(\be_i)\ge d_i$,
which implies that
$$N\ell(-p\be_i+q\be_j)\ge 1+d_j(q+a_j)-d_i(p+a_i+1)$$
whenever $p\le -a_i-1$, $q\ge -a_j$. Thus, our claim follows from Theorem \ref{pole-thm}.

Next, let us define the homomorphism
$$\wt{\psi}:\Q^n\to \Pic(\ov{\MM}_{g,n})_\Q: \be_i\mapsto \wt{\psi}_i.$$
Then we see that any monomial expression in coordinates $(\ov{\a}_{ij}[p,q])$ extends to a global section of $\wt{\psi}(v)$,
where $v$ is the corresponding $\G_m^n$-weight. 
Now the assertion (i) follows from the formula
$$\wt{\psi}(\om_{ij})=(a_i+1)\wt{\psi}_i-a_j\wt{\psi}_i=
\begin{cases} Z_{\ba+\be_i-\be_j}+(\frac{N}{a_i}-1)Z_\ba, & a_i>0, a_j>0; \\ Z_{\ba+\be_i-\be_j}, & a_i=0, a_j>0.\end{cases}$$
Indeed, in the case $a_i>0$, $a_j>0$ this follows from the identity
$$(a_i+1)\be_i-a_j\be_j=\Psi_{\ba+\be_i-\be_j}-\Psi_\ba,$$
where $\Psi_\ba$ are given by \eqref{Psi-eq}.
Similarly, in the case $a_i=0$, $a_j>0$ this follows from
$$\psi_i-a_j\psi_j+Z_\ba=\Psi_{\ba+\be_i-\be_j}-\la.$$
For the assertion (ii) we use in addition the formula for $\wt{\psi}_1$, which in the case $N=a_1=g$ becomes 
$\wt{\psi}_1=\psi_1+Z_{\be_1}$.
\ed

\begin{ex} In the case when $g$ is divisible by $n$ and $a_1=\ldots=a_n=g/n$, the class $Z(\chi)$ from
Proposition \ref{rational-classes-prop} 
is the image of $\chi$ under the linear map sending the generating vectors $\om_{ij}$ to
$Z_{\ba+\be_i-\be_j}$.
\end{ex}

\begin{rem} It follows easily from the results of Logan \cite{Logan} that if $\bC_\psi\sub \R\psi_1+\ldots+\R\psi_n$
is the convex hull of the classes $\Psi_\ba$, where $\ba$ ranges over $X(g,n)$, then for
any rational class $Z$ in the interior of $\bC_\psi+\R_{\ge 0}\psi_1+\ldots+\R_{\ge 0}\psi_n$, the $\Q$-divisor
$Z-\la$ is big. All the divisor classes appearing in Proposition \ref{rational-classes-prop} are multiples of classes of this form.
\end{rem}

%Example: $n=1$, $a_1=g$. Then $N=g$. We get that $\psi_1+Z=(1+{g+1\choose 2})\psi_1-\la$ is big on $\ov{M}_{g,1}$.

%$g=1$, $n>1$. For each $i$ have a cone generated by $e_i$ and $e_j-e_i$, $j\neq i$.
%Together they generate the entire region $x_1+\ldots+x_n\ge 0$.

%$n=g$, $a_1=\ldots=a_g=1$. Then $N=1$. We get that ??? are big on $\ov{M}_{g,g}$.

\section{Examples with $g=1$}\label{examples-sec}

\subsection{Case $g=1$, $n=2$, $\ba=(1,0)$}\label{g-1-n-2-sec}

In this section we work over $\Z[1/6]$.

The simplest case $g=n=1$, $a_1=1$ corresponds to the standard family of Weierstrass cubics (see e.g., 
\cite[Sec.\ 1.3]{LP}), so we start with the next simplest case $n=2$, $\ba=(1,0)$.

\begin{prop}\label{g-1-n-2-a-1-0-prop} 
There is an isomorphism $\wt{\UU}^{ns}_{1,2}(1,0)\simeq \A^4$ with coordinates $a,b,e,\pi$,
such that the universal affine curve $C\setminus\{p_1,p_2\}$ over $\wt{\UU}^{ns}_{1,2}(1,0)$
is given by the equations
\begin{equation}\label{g-1-n-2-1-0-curve-eq}
\begin{array}{l}
h_1^2=f_1^3+\pi f_1+s,\\
f_1h_{12}=ah_1+bh_{12}+ae,\\
h_1h_{12}=af_1^2+eh_{12}+abf_1+a(\pi+b^2),
\end{array}
\end{equation}
where 
$$s=e^2-b(\pi+b^2).$$
Here the functions $h_{12}\in H^0(C,\OO(p_1+p_2))$, $f_1\in H^0(C,\OO(2p_1))$ and $h_1\in H^0(C,\OO(3p_1))$  
have the expansions
$$h_{12}\equiv 1/t_2+\ldots \ \text{ at } p_2,$$ 
$$h_{12}\equiv a/t_1+\ldots, \  f_1\equiv 1/t_1^2+\ldots, \ h_1\equiv 1/t_1^3+\ldots \ \text{ at } p_1,$$
where $t_1$ and $t_2$ are formal parameters at $p_1$ and $p_2$ compatible with the fixed tangent vectors at $p_1$ and $p_2$. Note also that $b=f_1(p_2)$, $e=h_1(p_2)$.
The weights of the coordinates with respect to the $\G_m^2$-actions are:
$$wt(a)=\be_2-\be_1, \ wt(b)=2\be_1, \ wt(e)=3\be_1, \ wt(\pi)=4\be_1. 
$$
\end{prop}

\Pf . The proof is similar to that of \cite[Thm.\ 1.2.4]{P-ainf}. First, assuming that $6$ is invertible, 
we can construct canonical generators in any marked algebra $A$ of type $(1,0)$ 
over $R$ (see \eqref{marked-alg-def}). Namely,
we claim that there is a unique choice of $f_1\in F_{2\be_1}$, $h_1\in F_{3\be_1}$ and $h_{12}\in F_{\be_1+\be_2}$,
such that 
$$f_1\equiv u_1^2\mod F_{\be_1}, \ \ h_1\equiv u_1^3\mod F_{2\be_1}, \ \ h_{12}\equiv u_2\mod F_{\be_1},$$
and the following relations hold in $A$, with $\pi,s,a,b,d,e,p,q,r\in R$:
$$h_1^2=f_1^3+\pi f_1+s,$$
$$f_1h_{12}=ah_1+bh_{12}+d,$$
$$h_1h_{12}=af_1^2+eh_{12}+ph_1+qf_1+r.$$
Indeed, we have the following ambiguity in the choice of $f_1,h_1,h_{12}$: we can add constants to $f_1$ and $h_{12}$,
and can change $h_1$ to $h_1+c_1f_1+c_0$. The first relation fixes $f_1$ and $h_1$ uniquely, while adding
a constant to $h_{12}$ we make sure that there is no term with $f_1$ in the second equation.

The Buchberger's algorithm gives the following relations between the coefficients that are necessary and sufficient
for the elements
$$f_1^n, h_1f_1^m, h_{12}^k, \ \ n\ge 0, m\ge 0, k\ge 1,$$
to form a basis of $A$:
\begin{equation}\label{1-0-relations}
s=e^2-b(\pi+b^2), \ d=ae, \ p=0, \ q=ab, \ r=a(\pi+b^2).
\end{equation}

The rest of the proof is parallel to that of \cite[Thm.\ 1.2.4]{P-ainf}, considering the chain of the natural maps 
$$\wt{\UU}^{ns}_{1,2}(1,0)\to \MA_{(2,1)}\to \SS_{(1,0)}\to \wt{\UU}^{ns}_{1,2}(1,0),$$
where $\SS_{(1,0)}$ is the scheme defined by the equations \eqref{1-0-relations}, so $\SS_{(1,0)}\simeq \A^4$.
%defining equations of $H^0(C\setminus\{p_1,p_2\},\OO)$ take form
\ed

\begin{cor}
The forgetting morphism $\forg_2:\wt{\UU}^{ns}_{1,2}(1,0)\to \wt{\UU}^{ns}_1(1)$ identifies 
$\wt{\UU}^{ns}_{1,2}(1,0)$ with the total space $\Tot(T_{p_1})$ of the line bundle over the universal affine curve over $\wt{\UU}^{ns}_1(1)$ given by the tangent line at the marked point. This identification is compatible with the $\G_m^2$-action, where
the action of the first factor on $\Tot(T_{p_1})$ is induced by its action on $\wt{\UU}^{ns}_1(1)$, while the second factor
acts by rescalings in the fibers and trivially on the base.
\end{cor}

\Pf . The morphism $\forg_2$ is given by the functions $(\pi,s)$. Now we can rewrite the formula for $s$ as
the equation
$$e^2=b^3+\pi b+s,$$
which is exactly the equation of the universal affine curve over $\wt{\UU}^{ns}_1(1)$. The extra coordinate $a$ can be interpreted as the coordinate on the tangent line at the point $p_1$ via the identification of the latter with
$\OO(p_1)/\OO$, by recalling that $h_{12}\equiv a/t_1$ in $\OO(p_1)/\OO$.
\ed

\subsection{Case $g=1$, $n=3$, $\ba=(1,0,0)$}

\begin{prop}\label{g-1-n-3-a-1-0-0-prop} 
There is an isomorphism of $\wt{\UU}^{ns}_{1,3}(1,0,0)\simeq \A^1\times\RR_1(2,4)$,
where $\RR_1(2,4)$ is the variety of $2\times 4$-matrices $M=(m_{ij})$
such that $\rk M\le 1$ (its ideal is generated by the $2\times 2$-minors).
The weights of the $\G_m^3$-action on the coordinates $m_{ij}$ are given by
\begin{align*}
&w(m_{11})=\be_2-\be_1, w(m_{12})=\be_3-\be_1, w(m_{13})=2\be_1, w(m_{14})=3\be_1,\\
&w(m_{21})=\be_2, \  \ \ \ \ \ w(m_{22})=\be_3,  \ \ \  \ \ \   w(m_{23})=3\be_1, w(m_{24})=4\be_1,
\end{align*}
while the coordinate $t$ on the extra factor $\A^1$ has weight $2\be_1$.
%The universal affine curve $C\setminus\{p_1,p_2,p_3\}$ is given by ???
The scheme $\wt{\UU}^{ns}_{1,3}(1,0,0)$ is irreducible of dimension $6$, Cohen-Macauley, normal and nonsingular in codimension $4$, but not Gorenstein.
\end{prop}

\Pf . Again we argue as in \cite[Thm.\ 1.2.4]{P-ainf}.
Using the case $n=2$ (see Proposition \ref{g-1-n-2-a-1-0-prop})
we see that any marked algebra $A$ of type $(1,0,0)$ has unique generators $f_1\in F_{2\be_1}$,
$h_1\in F_{3\be_1}$, $h_{12}\in F_{\be_1+\be_2}$ and $h_{13}\in F_{\be_1+\be_3}$, such that
$$f_1\equiv u_1^2\mod F_{\be_1}, \ \ h_1\equiv u_1^3\mod F_{2\be_1}, \ \ h_{12}\equiv u_2\mod F_{\be_1},
\ \ h_{13}\equiv u_3\mod F_{\be_1},$$
and the following relations hold in $A$:
%the equations of the universal affine curve $C\setminus\{p_1,p_2,p_3\}$ have form
\begin{align*}
&h_1^2=f_1^3+\pi_1f_1+s_1,\\
&f_1h_{1i}=a_{1i}h_1+b_{1i}h_{1i}+a_{1i}e_{1i}, \\
&h_1h_{1i}=a_{1i}f_1^2+e_{1i}h_{1i}+a_{1i}b_{1i}f_1+a_{1i}(\pi_1+b_{1i}^2), \\
&h_{12}h_{13}=c_{32}h_{12}+c_{23}h_{13}+a_{12}a_{13}f_1+d,
\end{align*}
where $i=2,3$, and 
$$s_1=e_{12}^2-b_{12}^3-\pi_1b_{12}=e_{13}^2-b_{13}^3-\pi_1b_{13}.$$
%Here $h_{1i}\in H^0(C,\OO(p_1+p_i))$ has the polar part $\frac{1}{t_i}$ at $p_i$ (for $i=2,3$).
Note that the form of the last equation is dictated by the fact that
$h_{12}h_{13}-a_{12}a_{13}f_1\in F_{\be_1+\be_2+\be_3}$.
Using the Buchberger's algorithm, we see that the condition that 
$$(f_1^n,f_1^nh_1,h_{12}^k,h_{13}^k)_{n\ge 0, k\ge 1}$$ 
form a basis of $A$
%$H^0(C\setminus\{p_1,p_2,p_3\},\OO)$
is equivalent to the equation 
$$d=a_{12}a_{13}(b_{12}+b_{13}),$$
together the following equations on the
%an isomorphism of $\wt{\UU}^{ns}_{1,3}(1,0,0)$ with the subscheme in $\A^9$ with
cooefficients $a_{12}, a_{12}, b_{12}, b_{13}, e_{12}, e_{12}, c_{23}, c_{32}, \pi_1$:
\begin{align*}
&e_{12}^2-e_{12}^2=b_{12}^3-b_{13}^3+\pi_1(b_{12}-b_{12}),\\
&c_{23}(b_{13}-b_{12})=a_{12}(e_{12}+e_{13}),\\
&c_{32}(b_{12}-b_{13})=a_{13}(e_{12}+e_{13}),\\
&c_{23}(e_{13}-e_{12})=a_{12}(\pi_1+b_{13}^2+b_{12}b_{13}+b_{12}^2),\\
&c_{32}(e_{12}-e_{13})=a_{13}(\pi_1+b_{13}^2+b_{12}b_{13}+b_{12}^2),\\
&a_{12}c_{32}=-a_{13}c_{23}.
\end{align*}
These are exactly the equations given by the $2\times 2$-minors of the matrix 
$$\left(\begin{matrix} a_{12} & a_{13} & b_{13}-b_{12} & e_{13}-e_{12} \\
c_{23} & - c_{32} & e_{12}+e_{13} & \pi_1+b_{12}^2+b_{12}b_{13}+b_{13}^2\end{matrix}\right).$$
%$(m_{ij})_{(i=1,2;j=1,\ldots,4)}$ given by
%m_{11}=a_{12}, \ m_{12}=a_{13}, \ m_{13}

It remains to observe that the entries of this matrix complemented by $t=b_{12}$ form a change of variables from
the original $9$ coordinates (recall that we assume $2$ to be invertible).
The geometric properties of our space follow immediately from the well-known properties of the affine cone
over the Segre embedding of $\P^1\times\P^3$.
\ed

%Discuss morphism to $\wt{\UU}^{ns}_{1,2}(1,0)$, morphism of universal curves
%$\CC_{1,3}(1,0,0)\setminus p_3\to \CC_{1,2}(1,0)$.

\subsection{The moduli space $\wt{\UU}^{ns}_{1,n}(\be_1,\be_2,\ldots,\be_n)$}

An important example of the moduli scheme of the form 
$\wt{\UU}^{ns}_{g,n}(\ba_1,\ba_2,\ldots)$ is the scheme $\wt{\UU}^{ns}_{1,n}(\be_1,\be_2,\ldots,\be_n)$
because of its connection with the moduli scheme $\wt{\UU}^{sns}_{1,n}$
of genus $1$ curves $(C,p_1,\ldots,p_n)$ such that $h^1(\OO(p_i))=0$ and $\OO(p_1+\ldots+p_n)$ is ample,
with a nonzero global section of the dualizing sheaf $\om_C$.\footnote{``sns" stands for ``strongly non-special",
since each $p_i$ defines a non-special divisor. In \cite{LP} this scheme is denoted by $\wt{\UU}^{ns}_{1,n}$.} 
The scheme 
$\wt{\UU}^{sns}_{1,n}$ was constructed and studied in \cite{LP}, where we showed that its GIT quotient by $\G_m$
coincides with the moduli space of $(n-1)$-stable curves of genus $1$ defined and studied by Smyth in
\cite{SmythI}, \cite{SmythII}.

Let $T_1\sub \G_m^n$ be the kernel of the projection to the first component $\G_m^n\to \G_m$.

\begin{prop}\label{g-1-vns-prop} 
The action of $T_1$ on $\wt{\UU}^{ns}_{1,n}(\be_1,\be_2,\ldots,\be_n)$ is free and 
there is an isomorphism of $\G_m$-schemes
$$\wt{\UU}^{ns}_{1,n}(\be_1,\be_2,\ldots,\be_n)/T_1\simeq \wt{\UU}^{sns}_{1,n}$$
\end{prop}

\Pf . The restriction map $H^0(C,\om_C)\to \om_C|_{p_i}$ is an isomorphism for every $i$ (see \cite[Lem.\ 1.1.1]{LP}).
Hence, we have a natural embedding
$$\wt{\UU}^{sns}_{1,n}\to \wt{\UU}^{ns}_{1,n}(\be_1,\be_2,\ldots,\be_n),$$
associating to a trivialization of $H^0(C,\om_C)$ the corresponding 
trivializations of the tangent lines at each marked point. It is easy to see that its image is a section for the
$T_1$-action on $\wt{\UU}^{sns}_{1,n}$, which implies the assertion.
\ed

\section{Moduli of curves with chains of divisors}\label{special-curves-sec}

\subsection{The setup}\label{special-divisors-setup-sec}

The combinatorical data defining our moduli spaces
is a collection of formal divisors (formal integer combinations of $\bp_1,\ldots,\bp_n$):
$$0=\bD_0^-\le \bD_0^+\le \bD_1^-\le \bD_1^+\le \ldots \le \bD_s^-\le \bD_s^+,$$
such that $\supp(\bD_s^+)=\supp(\bD_s^+ -\bD_s^-)$, and for each $q=1,\ldots,s$ there
exists an index $i_q\in\{1,\ldots,n\}$ such that $\bD_q^-=\bD_{q-1}^+ + \bp_{i_q}.$
We set 
\begin{equation}\label{g-sum-D-p-eq}
g=\sum_{q=0}^s \deg(\bD_q^+ - \bD_q^-)=\deg(\bD_s^+)-s.
\end{equation}
We always assume that $g>0$.

We will construct the moduli space parametrizing curves $C$ of arithmetic genus $g$ with smooth marked points
$p_1,\ldots,p_n$, such that $\OO_C(p_1+\ldots+p_n)$ is ample, 
with restrictions on the number of sections of various divisors supported at $p_1,\ldots,p_n$.
We say that a divisor $D=\sum a_ip_i$ is of type $\bD=\sum a_i\bp_i$. Let 
$D_q^-, D_q^+$, for $q=0,\ldots,s$, be the divisors on $C$ of the types $\bD_q^-,\bD_q^+$.
Then we impose the conditions $h^1(C,D_s^+)=0$, and
\begin{equation}\label{h-0-rules-eq}
%\begin{array}{l}
h^0(C,D_q^-)=h^0(C,D_q^+)=q+1, \ \text{ for } q=0,\ldots,s,
%\\
%h^0(C,D_q^-)=h^0(C,D_{q-1}^+)+1 \ \text{ for } q=1,\ldots,p.
%\end{array}
\end{equation}
Note that since $\deg(D_q^-)=\deg(D_{q-1}^+)+1$, these conditions imply that
$h^1(C,D_q^-)=h^1(C,D_{q-1}^+)$. Thus, $h^1$ drops by $\deg(D_s^+ - D_s^-)$ only when going from
$D_s^-$ to $D_s^+$. The total drop from $h^1(D_0^-)=h^1(C,\OO)$ to $h^1(D_s^+)$ is $g$, which explains the formula
\eqref{g-sum-D-p-eq}. Let us set formally
\begin{equation}\label{formal-h1-eq}
%\begin{array}{l}
h^1(\bD_q^\pm)=g+q-\deg(\bD_q^\pm),
%g-\sum_{j=0}^{q}\deg(\bD_j^+ - \bD_j^-), \ q=0,\ldots,s,\\
%h^1(\bD_q^-)=h^1(\bD_{q-1}^+), q=1,\ldots,s,
%\end{array}
\end{equation}
so that in the above situation these are the values of $h^1(D_q^+)$ and $h^1(D_q^-)$.
%Note that since $h^0(C,D_0^-)=1$, from \eqref{h-0-rules-eq} we get that
%$$r_q^\pm:=h^0(C,D_q^\pm)$$

\begin{ex}
In the case $s=0$ the above data reduces to a single formal divisor $\bD_0^+=\sum_{i=1}^n a_i\bp_i$
of degree $g$. Thus, the corresponding moduli space, parametrizing $(C,p_1,\ldots,p_n)$, such that
$\OO(p_1+\ldots+p_n)$ is ample and the divisor $D_0^+=\sum a_ip_i$ satisfies $h^1(D_0^+)=0$, 
is exactly the stack $\UU^{ns}_{g,n}(\ba)$.
\end{ex}

\begin{ex}\label{gap-example} 
For $n=1$ our data is equivalent to the collection of numbers 
$$0=d_0^-\le d_0^+\le d_1^-\le d_1^+\le\ldots\le d_s^-<d_s^+,$$
where $d_q^-=d_{q-1}^+ +1$, so that $D_q^\pm=d_q^\pm p_1$. 
Then our condition on divisors supported at $p_1$ on a curve $C$ is that $h^1(d_s^+p_1)=0$ and
that the set
$$[1,d_s^+]\setminus \{d_q^- \ |\ q=1,\ldots, s\}$$ 
is precisely the set of Weierstrass gaps of $(C,p_1)$, i.e., the set of $m\ge 1$ such that
$h^0((m-1)p_1)=h^0(mp_1)$. Conversely, if $1=\ell_1<\ldots<\ell_g$ is the gap sequence of some 
subsemigroup $S\sub\Z_{\ge 0}$, i.e., $S=\Z_{\ge 0}\setminus \{\ell_1,\ldots,\ell_g\}$, then we can reconstruct our data by 
\begin{equation}\label{gap-divisors-eq}
d_s^+=\ell_g, \ \{d_q^- \ |\ q=1,\ldots, s\}=[1,\ell_g]\setminus \{\ell_1,\ldots,\ell_g\}.
\end{equation}
%For example, in the case when $C$ is smooth and $p_1$ is generic we have $\ell_i=i$, $i=1\ldots,g$ 
%so this corresponds to the case $s=0$, $d_0^+=g$.
\end{ex}
%in the union of segments
%$\sqcup_{q=0}^s [d_q^-+1,d_q^^+]$ form precisely the Weierstrass gap sequence of $(C,p_1)$.

%Let us fix a sequence $\bll_1,\ldots,\bll_r$ ???

The next Lemma records a well known construction of the loci over which the rank of a given map of vector bundles is constant.

\begin{lem}\label{rank-locus-lem} 
Let $\phi:\VV\to \WW$ be a morphism of vector bundles over a scheme $X$.
We define $Z_{\le r}(\phi)\sub X$ to be the closed subscheme given by the equation $\we^{r+1}(\phi)=0$,
and we call the locally closed subscheme in $X$,
$$Z_r(\phi):=Z_{\le r}(\phi)\setminus Z_{\le r-1}(\phi),$$ 
the locus where the rank of $\phi$ is equal to $r$. Then for a morphism $f:S\to X$ the sheaf
$\coker(f^*\phi)$ is locally free of rank $\rk\WW-r$ if and only if $f$ factors through $Z_r(\phi)$.
For such $f$ the sheaves $\im(f^*\phi)$ and $\ker(f^*\phi)$ are locally free of ranks $r$ and $\rk\VV-r$, respectively.
\end{lem}

\Pf . The question is local so we can assume that the bundles $\VV$ and $\WW$ are trivial. Note that the complement
of $Z_{\le r-1}(\phi)$ is the open subset of $X$ on which one of the $r\times r$-minors of the matrix of $\phi$ is
invertible. Equivalently, it is the union of open subsets on which there exists a decomposition 
$\WW=\WW_1\oplus \WW_2$ such that $\rk\WW_1=r$ and the component $\phi_1:\VV\to \WW_1$ of $\phi$
is surjective. Hence, locally we have a decomposition $\VV=\VV_1\oplus \VV_2$ such that $\phi_1(\VV_2)=0$ and
$\phi_1|_{\VV_1}:\VV_1\to\WW_1$ is an isomorphism. Thus, when we further
restrict to $Z_r(\phi)$, the condition that $\we^{r+1}(\phi)=0$ will imply that $\phi$ factors as a composition of
the projection $\VV\to\VV_1$ followed by an embedding as a subbundle $\VV_1\to\WW$. Thus,
$\im(\phi|_{Z_r(\phi)})$ is a subbundle of $\WW|_{Z_r(\phi)}$ of rank $r$.
 
Conversely, assume that $\coker(f^*\phi)$ is locally free of rank $\rk\WW-r$. Then $\im(f^*\phi)$
is a subbundle of $f^*\WW$ of rank $r$, so $\we^{r+1}(\phi)=0$ and
locally there exists an $r\times r$-minor of $f^*\phi$ which is invertible,
hence $f$ factors through $Z_r(\phi)$.
\ed

\begin{defi}\label{SG-D-defi}
With the data $\DD=(\bD_\bullet^\pm)$ we associate a locally closed subscheme $SG_\DD\sub SG_1(g)$
as follows. We start with an open subset $SG(\bD_s^+)\sub SG_1(g)$ consisting of $W$ such that $H^1(W(\bD_s^+))=0$.
Note that over $SG(\bD_s^+)$  
the spaces $H^0(W(\bD_s^+))/k$ are the fibers of the vector bundle $\KK(\bD_s^+)$ (see \eqref{K-f-C-f-def-eq}). 
Now for each $q=0,\ldots,s$ we have morphisms of vector bundles
$$\varphi_q^\pm: \KK(\bD_s^+)\to \HH_{\ge -\bD_s^+}/\HH_{\ge -\bD_q^\pm}\ot \OO,$$
and we define $SG_\DD\sub SG(\bD_s^+)$ as the intersection of the loci $Z_{s-q}(\varphi_q^\pm)$, for $q=0,\ldots,s$,
where the rank of $\varphi_q^\pm$ is equal to $s-q$
(see Lemma \ref{rank-locus-lem}). 
\end{defi}

Recall (see Lemma \ref{H-0-H-1-exact-seq-lem}) that for any $f:X\to SG(\bD_s^+)$
the morphism $f^*\varphi_q^\pm$ fits into an exact sequence
\begin{equation}\label{f*-phi-q-ex-seq}
0\to \KK(f,\bD_q^\pm)\to \KK(f,\bD_s^+)\rTo{f^*\varphi_q^\pm} \HH(\bD_s^+)/\HH(\bD_q^\pm)\ot\OO_X
\rTo{} \CC(f,\bD_q^\pm)\to 0.
\end{equation}

\begin{lem}\label{special-div-seq-lem} 
(i) For the embedding $i:SG_\DD\to SG(\bD_s^+)$ the sheaves
$\KK(i,\bD_q^\pm)$ and $\CC(i,\bD_q^\pm)$ are locally free of ranks
$q$ and $h^1(\bD_q^\pm)$, respectively.

\noindent (ii) For each $q$ the natural morphism $\KK(i,\bD_q^-)\to \KK(i,\bD_q^+)$ is an isomorphism.

\noindent (iii) For each $q\ge 1$ the natural morphism $\CC(i,\bD_{q-1}^+)\to \CC(i,\bD_q^-)$ is an isomorphism.
\end{lem}

\Pf . (i) This follows from the exact sequence \eqref{f*-phi-q-ex-seq} for $f=i$ together
with Lemma \ref{rank-locus-lem} and
the fact that the morphism $i^*\varphi_q^\pm$ has the constant rank $s-q$.

\noindent
(ii) By Lemma \ref{H-0-H-1-exact-seq-lem},
we have an exact sequence
\begin{equation}\label{K-V-seq-eq}
0\to \KK(i,\bD_q^-)\to \KK(i,\bD_q^+)\to \HH(\bD_q^+)/\HH(\bD_q^-)\ot\OO\to \VV\to 0,
\end{equation}
where $\VV$ fits into an exact sequence
$$0\to \VV\to \CC(i,\bD_q^-)\to \CC(i,\bD_q^+)\to 0.$$
Hence, $\VV$ is a vector bundle of rank $\deg(\bD_q^+-\bD_q^-)$, so the surjective morphism
$\HH(\bD_q^+)/\HH(\bD_q^-)\ot\OO\to \VV$ is in fact an isomorphism. Now
the assertion follows from the exact sequence \eqref{K-V-seq-eq}.

\noindent
(iii) By Lemma \ref{H-0-H-1-exact-seq-lem}, this morphism is surjective. It remains to observe that by part (i),
both sheaves are locally free of the same rank. Hence, this map is an isomorphism.
\ed

For a divisor $\bD$ let us denote by $m_i(\bD)$ the coefficient of $\bp_i$ in $\bD$. We also set
\begin{equation}\label{new-a-i-eq}
\begin{array}{l}
a_i=m_i(\bD_s^+), \ i=1,\ldots,n,\\
m_q=m_{i_q}(\bD_q^-), \ q=1,\ldots, s
\end{array}
\end{equation}
(recall that $i_q$ is determined by $\bD_q^-=\bD_{q-1}^+ + \bp_{i_q}$).
For each $i=1,\ldots, n$, let 
$$T_i=\{-m_q \ |\ q\in [1,s], i_q=i \}\cup (-\infty,-a_i-1],$$
and set 
$$S_i=\Z\setminus T_i, \ \bS=\sqcup_i S_i.$$
Recall that the open cell $U_{\bS,1}\sub SG_1(g)$ associated with $\bS$ is given by \eqref{U-S-1-eq}.

\begin{prop}\label{special-locus-open-cell-prop} 
$SG_\DD$ is a closed subscheme of the open cell $U_{\bS,1}$, which
is isomorphic to the infinite-dimensional affine space. 
\end{prop}

\Pf . By Lemma \ref{H-0-H-1-exact-seq-lem} and Lemma \ref{special-div-seq-lem}(iii), 
we have exact sequences of vector bundles over $SG_\DD$,
$$0\to \KK(\bD_{q-1}^+)\to \KK(\bD_q^-)\to \HH(\bD_q^-)/\HH(\bD_{q-1}^+)\ot \OO\to 0,$$
for $q=1,\ldots,s$, where the last vector bundle is canonically isomorphic to the trivial bundle
of rank $1$. On the other hand, for every $\bD\ge \bD_s^+$ we have an exact sequence
$$0\to \KK(\bD_s^+)\to \KK(\bD)\to \HH(\bD)/\HH(\bD_s^+)\ot\OO\to 0.$$
Hence, we can choose splittings of all of these sequences over $SG_\DD$,
and represent the universal subspace over $SG_\DD$ by a locally free subsheaf
$W\sub \HH\ot\OO$ (with locally free quotient), such that $W/\lan 1\ran$ has the basis 
\begin{equation}\label{special-divisor-basis-eq}
(g_q)_{q=1,\ldots,s}, (f_i[-m])_{i=1,\ldots,n,m>a_i}), 
\end{equation}
such that
\begin{equation}\label{special-divisor-basis-eq2}
g_q\equiv t_{i_q}^{-m_q}\mod \HH(\bD_{q-1}^+),
\end{equation}
$$f_i[-m]\equiv t_i^{-m}\mod \HH(\bD_s^+).$$
Adding to $f_i[-m]$ a linear combination of $g_q$, and to each $g_q$ a linear combination of $g_{q'}$ with $q'<q$, we 
arrive at the basis of the form \eqref{cell-basis-eq}, so we obtain the inclusion $SG_\DD\sub U_{\bS,1}$. 
In addition, the property \eqref{special-divisor-basis-eq2} implies that
for every $q=1,\ldots,s,$ the coefficient of $t_i^{-m}$ in $g_q$, where $-m\in S_i$, is zero unless
$m\le m_i(\bD_{q-1}^+)$. Note that these coefficients are among the coordinates on $U_{\bS,1}$, so their vanishing
gives a closed subscheme $Z\sub U_{\bS,1}$, isomorphic to the infinite-dimensional affine space.
Finally, it is easy to see that over $Z$ the image of the map
$\varphi_q^\pm$ is a subbundle of rank $s-q$ (for $q=1,\ldots,s$), which shows that $Z=SG_\DD$.
\ed

The following definition mimics Def.\ \ref{SG-D-defi} for the moduli of curves.

\begin{defi}\label{special-divisors-moduli-def} 
(i) For the data $\DD=(\bD_\bullet^\pm)$ as above let $\UU_{g,n}(\bD_s^+)$
denote the moduli stack of projective curves $C$ of arithmetic genus $g$ with $n$ marked points (distinct, smooth) $p_1,\ldots,p_n$
such that $\OO_C(p_1+\ldots+p_n)$ is ample and
$h^1(C,D_s^+)=0$, where $D_s^+$ is the divisor of type $\bD_s^+$ supported on $p_1,\ldots,p_n$.
%choose $N$ such that $\bD_s^+\le N(\bp_1+\ldots+\bp_n)$.
Then the moduli space
$\UU_{g,n}[\DD]$ is the locally closed locus of $\UU_{g,n}(\bD_s^+)$,
 defined as the intersection of the loci $Z_{s-q}(\phi_q^\pm)$, over $q\le s$,
on which the natural morphism of vector bundles 
\begin{equation}\label{R-0-R-0-morphism-eq}
\phi_q^\pm: R^0\pi_*(\OO(D_s^+))\to R^0\pi_*(\OO(D_s^+)/\OO(D_q^\pm)).
\end{equation}
has rank equal to $s-q$.
Here $\pi:C\to \UU_{g,n}(\bD_s^+)$ is the universal curve,
and $(D_q^\pm)$ is the divisor of type $(\bD_q^\pm)$ supported at $p_1,\ldots,p_n$. 
As before, $\wt{\UU}_{g,n}[\DD]$ (resp., $\UU^{(\infty)}_{g,n}[\DD]$) denote
the similar moduli space with the choice of nonzero tangent vectors (resp., formal parameters) at the marked points.

\noindent
(ii) In the case $n=1$, for a sequence $1\le \ell_1<\ldots<\ell_g$ we set
$$\UU_{g,1}[\ell_1,\ldots,\ell_g]=\UU_{g,n}[\DD], \ \wt{\UU}_{g,1}[\ell_1,\ldots,\ell_g]=\wt{\UU}_{g,n}[\DD],$$
where $\DD=(d_q^\pm\bp_1)_{q=0,\ldots,s}$ is defined by
\eqref{gap-divisors-eq}.
\end{defi}

\begin{lem} The stack $\UU_{g,n}[\DD]$ classifies families of pointed curves $(\pi:C\to S,p_1,\ldots,p_n)$ such that
$\OO(p_1+\ldots+p_n)$ is relatively ample and $R^1\pi_*(\OO_C(D_q^\pm))$ is locally free of rank $h^1(\bD_q^\pm)$ for each
$q=0,\ldots,s$, where $D_q^\pm$ is the divisor of type $\bD_q^\pm$.
For such a family and for any $\bD$, $\bD_q^-\le \bD\le \bD_q^+$, $q=0,\ldots,s$, the sheaf
$R^1\pi_*(\OO_C(D))$, where $D$ is the divisor of type $\bD$,
 is locally free of rank $h^1(\bD_q^+)+\deg(\bD_q^+ - \bD)$.
\end{lem}

\Pf . We have $R^1\pi_*(\OO(D_s^+))=0$ on $\UU_{g,n}(\bD_s^+)$. Hence, for $q=0,\ldots,s$,
the sheaf $R^1\pi_*(\OO_C(D_q^\pm))$
is isomorphic to the cokernel of \eqref{R-0-R-0-morphism-eq}. Therefore, the first assertion follows
from Lemma \ref{rank-locus-lem}. To prove the second assertion we use the commutative diagram
with exact rows and columns
\begin{diagram}
\pi_*(\OO(D)/\OO(D_q^-))&\rTo{\id}&\pi_*(\OO(D)/\OO(D_q^-))\\
\dTo{}&&\dTo{}\\
\pi_*(\OO(D_q^+)/\OO(D_q^-))&\rTo{}&R^1\pi_*(\OO(D_q^-))&\rTo{}&R^1\pi_*(\OO(D_q^+))\to 0\\
\dTo{}&&\dTo{}&&\dTo{\id}\\
\pi_*(\OO(D_q^+)/\OO(D))&\rTo{}& R^1\pi_*(\OO(D))&\rTo{}&  R^1\pi_*(\OO(D_q^+))\to 0.
\end{diagram}
Since by assumption $R^1\pi_*(\OO_C(D_q^-))$ and $R^1\pi_*(\OO_C(D_q^+))$ are locally free of prescribed
ranks, whose difference is $\deg(\bD_q^+ - \bD_q^-)$, we derive that the first arrow in the middle row is injective.
Since the left column is a short exact sequence, we deduce that the first arrow in the bottom row is injective,
so the bottom row is a short exact sequence. Hence, $R^1\pi_*(\OO(D))$ is locally free of prescribed rank.
\ed
%The base change theorem implies that 

%\begin{lem}
%\end{lem}

%\begin{ex}\label{gap-moduli-example}

Using this Lemma we can rephrase Definition \ref{special-divisors-moduli-def}(ii) 
as follows.
 
\begin{cor}\label{gap-moduli-cor}
For every gap sequence $1=\ell_1<\ldots<\ell_g$ the stack $\UU_{g,1}[\ell_1,\ldots,\ell_g]$
classifies families of pointed curves $(\pi:C\to S,p_1)$ such that
$\OO(p_1)$ is relatively ample, $R^1\pi_*(\OO_C(\ell_gp_1))=0$
and $R^1\pi_*(\OO_C(mp_1))$ is locally free of rank $g-i$ for $\ell_i\le m<\ell_{i+1}$, $i=0,\ldots,g-1$ (where $\ell_0=0$).
\end{cor}

\subsection{Curves with special divisors and the Krichever map}

It is clear from the definitions that the Krichever map defines a map
\begin{equation}\label{Kr-special-map}
\Kr:\UU^{(\infty)}_{g,n}[\DD]\to ASG_\DD:=ASG\cap SG_\DD,
\end{equation}
equivariant  with respect to the action of the group of changes of parameters $\fG$ (see Section \ref{changes-parameter-sec}).
%such that $R^0\pi_*(???)\simeq \Kr^*???$, etc.

%What does functor $\wt{\UU}_{g,n}[\DD]$ represent???

\begin{lem}\label{special-section-lem}
The action of $\fG$ on $SG_\DD$ is free and admits a section
$\Sigma_\DD\sub SG_\DD$, isomorphic to the infinite-dimensional affine space.
\end{lem}

\Pf . Renumbering the indices we can assume that 
$$\bD_s^+=a_1\bp_1+\ldots+a_r\bp_r,$$
with $a_i>0$, $r\le n$. By assumption, each of the formal points $\bp_1,\ldots,\bp_r$ is in the support of
$\bD_s^+-\bD_s^-$. In other words, 
$$\bD_s^-\le (a_1-1)\bp_1+\ldots+(a_r-1)\bp_r.$$
We know that $\CC(i,\bD_s^+)=0$, where $i:SG_\DD\to SG_1(g)$ is the embedding.
%The fact that $H^1(W(\bD_p^+))=0$ on $SG_\DD$ implies that the morphism
%$$\pi_{\bD_p^+}=\pi_{(m_1,\ldots,m_r)}:\HH_{\ge -\bD_p^+}/\HH_{\ge 0}\to \sV$$
%is surjective when restricted to $SG_\DD$. 
On the other hand, we claim that for each $i=1,\ldots,r$, the sheaf
%the restriction to $SG_\DD$ of the morphism
%$$\pi_{\bD_p^+-\bp_i}:\HH_{\ge -\bD_p^+ +\bp_i}/\HH_{\ge 0}\to \sV$$
$\CC(i,\bD_s^+-\bp_i)$ is locally free of rank $1$.
Indeed, by Lemma \ref{H-0-H-1-exact-seq-lem}, we have an exact sequence 
%we observe that for each $i=1,\ldots,r$, and for $v\in H^0(W(\bD_p^+ + jp_i))$, where $j>0$,
$$0\to \CC(i,\bD_s^+-\bp_i)\to \CC(i,\bD_s^+)\to \HH(\bD_s^+)/\HH(\bD_s^+-\bp_i)\ot\OO\to \CC(i,\bD_s^+-\bp_i)\to 0.$$
Since $\bD_s^-\le \bD_s^+ -\bp_i$, by Lemma \ref{special-div-seq-lem}(ii), the first arrow is an isomorphism,
so the map
$\HH(\bD_s^+)/\HH(\bD_s^+-\bp_i)\ot\OO\to \CC(i,\bD_s^+-\bp_i)$ is also an isomorphism, which proves our claim.

Therefore, Proposition \ref{section-prop} is applicable to the action of $\fG_i$ on $SG_\DD$. Arguing as in Corollary
\ref{section-open-cor}, we get a section for the action of $\prod_{i=1}^r \fG_i$ on $SG_\DD$ given by the locus of $W$
with the following property: each element $f_i[-m]\in H^0(W(\bD_s^+ +(m-a_i)\bp_i))$, for $i=1,\ldots,r$, $m>a_i$, 
from the basis \eqref{special-divisor-basis-eq}, has the expansion in $t_i$ of the form
$$f_i[-m]\equiv t_i^{-m} \mod t_i^{-a_i+1}k[[t_i]].$$
In other words, we require the vanishing of the coefficient of $t_i^{-a_i}$ in $f_i[-m]$. 
Next, as in Proposition \ref{section-open-prop}, for a fixed $i_0\le r$, 
we get a section for the action of the full group $\fG$ on $SG_\DD$ by requiring in addition
for each $i>r$ the expansion of $f_i[-1]\in H^0(W(\bD_s^+ +\bp_i)$ in $t_i$ to be of the form
$$f_i[-1]=t_i^{-1}+const,$$
and the expansion of $f_i[-1]$ in $t_{i_0}$ to have no constant term.

Thus, our section $\Sigma_\DD$ is the vanishing locus of some
set of coordinates in the infinite-dimensional affine space $SG_\DD$, hence, $\Sigma_\DD$ itself is an infinite-dimensional affine space.
\ed

\begin{thm}\label{special-divisors-thm} Let us work over $\Q$.
The action of $\fG$ on $ASG_\DD$ is free and admits a section.
The Krichever map induces an isomorphism
$$\wt{\UU}_{g,n}[\DD]\simeq ASG_\DD/\fG,$$
and $\wt{\UU}_{g,n}[\DD]$ is an affine scheme of finite type over $\Q$.
\end{thm}

\Pf . 
%The proof follows the steps similar to those in Theorems A and B. 
Let $\Sigma_\DD\sub SG_\DD$ be the section for the $\fG$-action constructed in Lemma \ref{special-section-lem}.
Then $\Sigma_\DD\cap ASG_\DD$ is an affine scheme, isomorphic to $ASG_\DD/\fG$.

The proof follows the same outline as the proof of Theorems A and B (see Section \ref{proof-A-B-sec}), so we will omit
some details. The key part of the proof is constructing the sequence of morphisms 
$$\wt{\UU}_{g,n}[\DD]\rTo{\ov{\Kr}} \Sigma_\DD\cap ASG_\DD \rTo{i} S_{GB}\rTo{r}  \wt{\UU}_{g,n}[\DD],$$
where $S_{GB}$ is a certain affine scheme of finite type classifying algebras with Gr\"obner bases of given form, 
and proving that $r\circ i\circ \ov{\Kr}=\id$, $\ov{\Kr}\circ r\circ i=\id$.

The morphism $\ov{\Kr}$ is simply induced by \eqref{Kr-special-map} by passing to quotients by $\fG$
and using the identification $\Sigma_\DD\cap ASG_\DD\simeq ASG_\DD/\fG$.

Let $W_R\sub \HH_R$ be
the universal subspace over $\Spec(R):=\Sigma_\DD\cap ASG_\DD$, and
let $N$ be such that $\bD_s^+\le N(\bp_1+\ldots+\bp_n)$.
Then one has $H^1(W_R(N(p_1+\ldots+p_n)))=0$, and
part of the basis of $W_R/R\cdot 1$, constructed in Proposition \ref{special-locus-open-cell-prop},
gives a basis of $H^0(W_R(N(p_1+\ldots+p_n)))/R\cdot 1$, namely, the elements
\begin{equation}\label{special-basis-N-eq}
(g_q)_{q=1,\ldots,s}, (f_i[-p])_{i=1,\ldots,n, a_i<p\le N}
\end{equation}
(see \eqref{special-divisor-basis-eq}). 
Thus, the condition ($\star$) of Section \ref{Groebner-sec} is satisfied, and we can apply Proposition \ref{Groebner-prop},
setting as before $h_i(j)=f_i[-N-j-1]$, $f_i=h_i(0)$. This gives us
generators and the basis of normal monomials in $W_R$ of the form \eqref{normal-monomials-eq}, with
the elements \eqref{special-basis-N-eq} playing the role of $(g_k)$.
As in Section \ref{proof-A-B-sec}, we define $S_{GB}$ to be the closed subscheme in the corresponding affine scheme of
Gr\"obner relations of the form \eqref{Groebner-leading-eq} by requiring in addition 
the relations with the leading terms $f_i[-p]f_i$ and $f_i[-p]f_{i'}$ to have form \eqref{S-GB-eq}, and those with the leading
terms $g_qf_i$ to have form
\begin{align*}
&g_qf_i=h_i(m_q)+A_i(\le N+m_q)+\sum_{k\neq i} A_k(\le 2N)+\text{terms}(\deg_1\le N) \ \text{ for } i=i_q,\\
&g_qf_i=B_i(\le N+1+m_i(\bD_q^-))+\sum_{k\neq i} B_k(\le 2N)+\text{terms}(\deg_1\le N)  \ \text{ for } i\neq i_q,
\end{align*}
where $A_i(\le a)$ and $B_i(\le a)$ are some linear combinations of the elements
$h_i(l)$ and $f_i$ of $\deg_1\le a$.
The relations do have this form for $W_R$, so we get a natural morphism $i:\Spec(R)\to S_{GB}$.

The morphism $r$ is associated with the family of curves over $S_{GB}$ given by $C_{GB}=\Proj(\RR(A_{GB}))$,
where $A_{GB}$ is the universal algebra defined by Gr\"obner relations over $S_{GB}$, $\RR(A_{GB})$ is
the Rees algebra associated with the increasing filtration induced by $\deg_1$.
Exactly as in Section \ref{proof-A-B-sec}, we equip $C_{GB}$ with smooth marked points $p_1,\ldots,p_n$ 
and show that it is a flat
family over $S_{GB}$. Furthermore, as in Section \ref{proof-A-B-sec}, we use the special form of the relations
to estimate the poles of the generators $h_i(j)$, $f_i[-p]$ and $g_q$ at the marked points.
Namely, we show that $h_i(j)$ has a pole of order $N+j+1$ at $p_i$ and a pole of order $\le N$ at $p_{i'}$ for $i'\neq i$;
$f_i[-p]$ has a pole of order $p$ at $p_i$ and a pole of order $\le a_{i'}$ at $p_{i'}$ for $i'\neq i$; and
$g_q$ has a pole of order $m_q$ at $p_{i_q}$ and a pole of order $\le m_i(\bD_q^-)$ at $p_i$ for 
$i\neq i_q$. This gives the following bases for the relevant divisors on $C_{GB}$:
\begin{equation}
\begin{array}{l}
H^0(C_{GB}, \OO(mD)) \ \text{ for } m\ge N: \ \ (g_q)_{q=1,\ldots,s}, (f_i[-p])_{i=1,\ldots,n, a_i<p\le N}, 
(h_i(j))_{i=1,\ldots,n, 0\le j\le m-N-1}, \\
H^0(C_{GB}, \OO(D_q^+))=H^0(C_{GB}, \OO(D_q^-)): \ \ 1, (g_{q'})_{q'=1,\ldots,q},
\end{array}
\end{equation}
where $D_q^\pm$ denotes the divisor supported at $p_1,\ldots,p_n$ of type $\bD_q^\pm$.
Hence, we derive that $C_{GB}\to S_{GB}$ is a family of genus $g$ curves, and from Riemann-Roch we get
that $h^1(C_s,\OO(D_s^+))=0$ for each curve in this family. Furthermore, the above explicit bases
are compatible with the maps \eqref{R-0-R-0-morphism-eq}, so we derive that each sheaf
$R^1\pi_*(\OO(D_q^\pm))$ over $S_{GB}$ is locally of rank $h^1(\bD_q^\pm)$. 

Thus, we get a well-defined morphism $r:S_{GB}\to  \wt{\UU}_{g,n}[\DD]$.
The rest of the proof is similar to that of Sec.\ \ref{proof-A-B-sec}. Note only that in proving that 
$\ov{\Kr}\circ r\circ i=\id$ we have to establish an isomorphism of two $R$-points of
$\Sigma_\DD\cap ASG_\DD\simeq ASG_\DD/\fG$. For this we again can use Proposition \ref{marked-alg-prop}(i)
which reduces this to checking that the corresponding weakly marked algebras are isomorphic.
\ed

Note that Theorem E is a special case of Theorem \ref{special-divisors-thm} when $n=1$, where we
set $ASG[\ell_1,\ldots,\ell_g]=ASG_\DD$ for $\DD$ given by \eqref{gap-divisors-eq}.

It is straightforward to see that in the case $n=1$ all the coordinates on the affine space $\Si_\DD$ 
have positive $\G_m$-weights. Thus, similarly to Proposition \ref{rw-invariant-prop} we derive the following result.

\begin{cor} Let us work over an algebraically closed field $k$ of characteristic zero. 
Let $1=\ell_1<\ldots<\ell_g$ be a gap sequence, i.e.,
$$S=\Z_{\ge 0}\setminus\{\ell_1,\ldots,\ell_g\}$$
is a subsemigroup in $\Z_{\ge 0}$.
The action of $\G_m$ on the ring of functions on the
scheme $\wt{\UU}_{g,1}[\ell_1,\ldots,\ell_g]$ has non-negative weights, with zero weight subspace given by constants.
The unique $\G_m$-invariant point of $\wt{\UU}_{g,1}[\ell_1,\ldots,\ell_g]$ corresponds
to the projective curve $C^S$ compactifying the curve $\Spec(k[S])$, where $k[S]$ is the semigroup ring of $S$.
The geometric quotient 
\begin{equation}\label{Uprime-eq}
U'_{g,1}[\ell_1,\ldots,\ell_g]:=\wt{\UU}_{g,1}[\ell_1,\ldots,\ell_g]\setminus\{C^S\}/\G_m
\end{equation}
is a projective scheme.
\end{cor}

Recall that the maximal gap $\ell_g$ of a smooth point $p$ on a curve of arithmetic genus $g$ is always $\le 2g-1$
(this follows from the fact that $S=\Z_{\ge 0}\setminus\{\ell_1,\ldots,\ell_g\}$ is a subsemigroup, see e.g., 
\cite{Kunz}). Using the $\G_m$-action on $\wt{\UU}_{g,1}[\ell_1,\ldots,\ell_g]$ we deduce that any curve that
has a smooth point for which $\ell_g=2g-1$ is
necessarily Gorenstein (provided it is irreducible).

\begin{cor}\label{Gorenstein-cor} 
Let $C$ be a projective irreducible and reduced curve of arithmetic genus $g\ge 1$ with a smooth point $p$ such that
$h^1((2g-2)p)=1$, or equivalently the maximal gap at $p$ is $\ell_g=2g-1$. 
Then $C$ is Gorenstein and $\om_C\simeq \OO((2g-2)p)$.
\end{cor}

\Pf . Let $1=\ell_1<\ldots<\ell_g$ be the gap sequence of $(C,p)$.  The fact that $\ell_g=2g-1$ 
implies that the subsemigroup $S=\Z_{\ge 0}\setminus\{\ell_1,\ldots,\ell_g\}$ is symmetric,
 i.e., for an integer $m$, one has $m\in S$ if and only if
$2g-1-m\not\in S$ (see \cite[Lemma]{Kunz}). 
Hence, by the main theorem of \cite{Kunz}, the semigroup ring $k[S]$ is Gorenstein.
But the $\G_m$-action on $\wt{\UU}_{g,1}[\ell_1,\ldots,\ell_g]$ gives a deformation over $\A^1$ with the special fiber
$\Spec(k[S])$ and every other fiber isomorphic to $C\setminus \{p\}$. Hence, $C\setminus \{p\}$ is Gorenstein as well,
and so $C$ is Gorenstein.

For the last assertion we observe that $L=\om_C(-(2g-2)p)$ is a line bundle of degree $0$ on $C$. 
Furthermore, by Serre duality, we get $h^0(L)=1$. Since $C$ is irreducible and reduced,
this implies that $L$ is trivial.
%(see \cite[Ch.\ IV]{Serre}).???
% $L$ is numerically trivial (see ???). By Serre duality we have
%$h^0(L)=1$, so we derive that $L$ is trivial (see ???). 
\ed

\begin{rems}\label{Stohr-rem}
1. Let $1=\ell_1<\ldots<\ell_g$ be a gap sequence with $\ell_g=2g-1$. 
Assuming in addition that $\ell_2=2$ and $\ell_{g-1}\ge g$
St\"ohr constructed in \cite{Stohr} a projective scheme $S[\ell_1,\ldots,\ell_g]$ whose points are in bijection with
Gorenstein projective irreducible 
curves $C$ of arithmetic genus $g$ with a smooth point $p$ such that $(C,p)$ has finite automorphism
group and the gap sequence of $(C,p)$ is $(\ell_1,\ldots,\ell_g)$. His method is to study the equations of $C$ in the
canonical embedding. It is easy to see that by construction $S[\ell_1,\ldots,\ell_g]$ carries a family satisfying the conditions
of Corollary \ref{gap-moduli-cor}, so we get a morphism
\begin{equation}\label{S-U'-map}
S[\ell_1,\ldots,\ell_g]\to U'_{g,1}[\ell_1,\ldots,\ell_g]
\end{equation}
where the scheme on the right is given by \eqref{Uprime-eq}. 
Corollary \ref{Gorenstein-cor} implies that this morphism is a bijection on $k$-points.
We conjecture that in fact, it is an isomorphism. 

The construction of \cite{Stohr} has been extended in \cite{OS} to
gap sequences with $\ell_g=2g-2$, $\ell_3=3$, $\ell_{g-1}\ge g$, to give a quasi-projective variety $S[\ell_1,\ldots,\ell_g]$
parametrizing Gorenstein irreducible curves with such gap sequences. It seems plausible that in this
case the map \eqref{S-U'-map} is an open embedding.
%in this case it is not clear anymore whether the Gorenstein condition is automatic.???
%Assume that the gap sequence $1=\ell_1<\ldots<\ell_g$ satisfies
%\end{rem}
%\begin{rem}\label{Pinkham-rem} 

\noindent
2.
The construction of the moduli space of pointed curves $(C,p)$ with a given gap sequence, and such that $C$ is smooth, goes back to Pinkham's work \cite{Pinkham}. He also constructs a compactification of the space, using
the versal deformation space of the monomial curve associated with the gap sequence, viewed as a $\G_m$-curve.
It would be interesting to compare this compactification to our space $\wt{\UU}_{g,1}[\ell_1,\ldots,\ell_g]$.
In the case $\ell_g=2g-1$ St\"ohr shows that Pinkham's space for smooth curves is an open subset of the moduli
space considered in \cite{Stohr}.
\end{rems}

The following easy example of the moduli space $\wt{\UU}_{g,1}[\ell_1,\ldots,\ell_g]$ is related to hyperelliptic curves.

\begin{prop} Let us work over $\Spec\Z[1/2]$.
For each $g\ge 2$, the moduli scheme
$\wt{\UU}_{g,1}[1,3,\ldots,2g-1]$ is isomorphic to the affine space $\A^{2g}$, and the universal affine curve $C\setminus\{p\}$
over it is given by the equation
\begin{equation}\label{hyperell-equation}
y^2=x^{2g+1}+a_1x^{2g-1}+\ldots+a_{2g}.
\end{equation}
The weights of the $\G_m$-action on the coordinates are given by $wt(a_i)=2i+2$.
In other words, this is just the miniversal deformation space of the $A_{2g}$-singularity with its natural $\G_m$-action.
\end{prop}

\Pf . Given a curve $(C,p)$ in $\wt{\UU}_{g,1}[1,3,\ldots,2g-1]$, 
we can find elements $x\in H^0(C,2p)$ and $y\in H^0(C,(2g+1)p)$ with the expansions $x\equiv \frac{1}{t^2}+\ldots$,
$y\equiv \frac{1}{t^{2g+1}}+\ldots$ at $p$. Then there should be a relation of the form
$$y^2=P_{2g+1}(x)+yP_g(x),$$
where $P_{2g+1}$ and $P_g$ are polynomials of degree $2g+1$ and $g$, respectively.
We can use the ambiguity $y\mapsto y+Q_g(x)$, where $Q_g$ is a polynomial of
degree $g$, to get $P_g=0$, and the ambiguity $x\mapsto x+const$ to make the coefficient of $x^{2g}$ in $P_{2g+1}$ to be $0$. Since this is a single equation, it already gives the Gr\"obner basis, so that the monomials $(x^m, x^my)_{m\ge 0}$
form a basis of $H^0(C\setminus\{p\},\OO)$. Conversely, starting from the algebra $A$ with the defining equation 
\eqref{hyperell-equation} between the generators $x,y$, we can construct a curve $C$ as $\Proj(\RR(A))$,
where $\RR(A)$ is the Rees algebra of $A$ with respect to the filtration induced by $\deg(x)=2$, $\deg(y)=2g+1$.
Similarly to the proof of \cite[Thm.\ 1.2.4]{P-ainf} one checks that these constructions (that work in families)
are mutually inverse.
%he coefficients of $P_{2g+1}$ (except the top one) give the coordinates on the moduli space.
\ed
%\subsection{Example: $g=2$, ???}

%Another example: two fixed points $p_1$ and $p_2$, $h^0(p_1+p_2)=2$, $h^0(2p_1+2p_2)=3$, etc. ???

\begin{rem} Thus, the moduli stack $\UU_{g,1}[1,3,\ldots,2g-1]\setminus\{C_0\}/\G_m$, 
where $C_0$ is the curve with the affine part
$y^2=x^{2g+1}$, is isomorphic to the weighted projective stack 
$\P(4,6,\ldots,4g+2)$. Note that the same stack also appears in \cite{Fed-hyp} as the stack
$\HH_{2g}[2g-1]$ of quasi-admissible hyperelliptic covers of $\P^1$, and in \cite{ASvdW} as the stack
of cyclic covers of $\P^1$ of degree $2$ with a fixed simple ramification point.
\end{rem}

%Another example: 
%$g=2$, $n=2$, $h^0(p_1)=1$, $h^0(p_1+p_2)=h^0(p_1+2p_2)=2$, hence, $h^1(2p_1+p_2)=0$. 
%Choose $f\in H^0(p_1+p_2)$ such that $f\equiv \frac{1}{t_2}+\ldots$ at $p_2$, and 
%$h\in H^0(2p_1+2p_2)$ and $k\in H^0(3p_1+2p_2)$ with expansions $\frac{1}{t_1^2}+\ldots$, $\frac{1}{t_1^3}+\ldots$.
%Basis: $(f^n, h^m, h^mk)$, $n\ge 1$, $m\ge 0$.???

%Version: $h^0(p_1)=1$, $h^0(p_1+p_2)=h^0(2p_1+p_2)=2$.
%Then choose $f\in H^0(p_1+p_2)$ as before,
%$h_1\in H^0(3p_1+p_2)$, $h_2\in H^0(4p_1+p_2)$, $h_3\in H^0(5p_1+p_2)$ such that ???

\end{document}